\documentclass[11pt]{article}
\usepackage[a4paper,margin=25mm]{geometry}
\usepackage[T1]{fontenc}
\usepackage{lmodern,microtype}
\usepackage{amsmath,amssymb,amsthm,mathtools,mathrsfs,esint}
\usepackage{xcolor}
\usepackage{tikz}
\usetikzlibrary{arrows.meta,calc,decorations.pathreplacing}
\usepackage[font=small,labelfont=bf]{caption}
\usepackage[section]{placeins}
\usepackage{needspace}
\usepackage[hidelinks]{hyperref}
\hypersetup{pdftitle={Subsonic Euler flows past an airfoil with an attached contact discontinuity in a semi-infinite nozzle},
	pdfsubject={Euler--Lagrange formulation, elliptic estimates with discontinuous coefficients, and semi-infinite nozzle flows},
	pdfauthor={Jun Chen, Xiaoguang You, Aibin Zang}}
\numberwithin{equation}{section}
\newtheorem{proposition}{Proposition}[section]
\newtheorem{lemma}[proposition]{Lemma}
\newtheorem{theorem}{Theorem}
\theoremstyle{definition}
\newtheorem{problem}{Problem}
\newtheorem*{mainproblem}{Main Problem}
\theoremstyle{remark}
\newtheorem*{remark}{Remark}
\theoremstyle{plain}
\allowdisplaybreaks
\title{Subsonic Euler flows past an airfoil with an attached contact discontinuity in a semi-infinite nozzle}
\author{Jun Chen\thanks{Three Gorges Mathematical Research Center,
		China Three Gorges University, No.~8 University Road,
		Yichang, Hubei 443002, China.}
	\and Xiaoguang You\thanks{Corresponding author. School of Mathematical Sciences,
		Jiangxi Science and Technology Normal University, No.~589 Xuefu Avenue,
		Nanchang, Jiangxi 330038, China.
		E-mail: \href{mailto:xgyou@jxstnu.edu.cn}{\texttt{xgyou@jxstnu.edu.cn}}.}
	\and Aibin Zang\thanks{College of Artificial Intelligence and Information Engineering,
		Yichun University, No.~576 Xuefu Road,
		Yichun, Jiangxi 336000, China.}}
\date{}
\begin{document}
	\maketitle
	\begin{abstract}
		For variable inlet data and prescribed downstream pressure, we prove the existence and uniqueness of two-dimensional
		steady subsonic Euler flows past an airfoil attached to the entrance of a semi-infinite
		nozzle, with a contact discontinuity issuing from its cusped trailing edge.
		An Euler--Lagrange coordinate transformation is employed to decouple the elliptic--hyperbolic system and
		reduce the free boundary problem to a nonlinear elliptic problem in a fixed strip
		with a finite slit. To handle the resulting boundary degeneracy, we develop a new approach in which tubular neighbourhoods of the slit are excised to obtain Lipschitz approximating domains. On these domains, barrier functions adapted to the contact discontinuity conditions are constructed in suitably chosen local comparison regions to obtain weighted H\"older estimates independent of the tube thickness. These estimates allow us to pass to the limit and construct the desired flow.
	\end{abstract}
	
	\section{Introduction}
	\label{sec:physical}
	
	In the early twentieth century, Prandtl developed the boundary-layer
	theory and the lifting-line theory \cite{Prandtl1904,Prandtl1918}. His work
	clarified the role of trailing vorticity in airfoil theory and motivated
	the study of flows past airfoils with attached vortex sheets. In
	two-dimensional steady flow, such a sheet is represented by a curve
	across which the tangential velocity may jump. This curve is a contact
	discontinuity for the Euler equations and forms a free boundary whose
	position should be determined together with the flow.
	
	To handle the free boundary, S.-X.~Chen \cite{ChenMach} introduced
	an Euler--Lagrange transformation in his study of a Mach configuration.
	The transformation maps the unknown contact line to a fixed straight
	line. Bae \cite{Bae} later used these coordinates to study subsonic
	contact discontinuities in an infinite nozzle. The resulting fixed-domain
	formulation, together with the elliptic estimates of Li and Vogelius
	\cite{LiVogelius} for discontinuous coefficients, yields stability under
	small wall perturbations. Bae and Park \cite{BaePark} developed a
	different formulation based on Helmholtz decomposition. Using
	compensated compactness, Chen, Huang, Wang and Xiang \cite{CHWX2019}
	obtained subsonic flows with large vorticity and characteristic
	discontinuities in general infinitely long nozzles. Chen, Deng and You
	\cite{CDY} subsequently treated general piecewise smooth subsonic
	upstream profiles without small-oscillation or monotonicity assumptions.
	Their solutions exist uniquely on a maximal interval of wall
	perturbations and have a $C^{2,\alpha}$ contact curve. These nozzle-flow
	results do not include an airfoil.
	
	For flows past airfoils, the classical work of Bers \cite{Bers1954,Bers1958} and
	Finn and Gilbarg \cite{FinnGilbarg1957,FinnGilbarg3D1957} concerned
	potential flows, in which vorticity is absent. These results therefore
	do not address an attached vortex sheet. For the full Euler equations,
	Chen, Xin and Zang \cite{CXZ} established existence, uniqueness and
	stability of subsonic flows past a thin airfoil in a half-plane. The
	vortex line issuing from the trailing edge joins the cusped profile
	with $C^{1,\alpha}$ regularity. Chen \cite{ChenAir} subsequently
	considered the whole-plane problem, where vortex lines are attached
	to both edges and the asymptotic flow direction is determined as part of the solution rather than an imposed condition.
	In both works, an
	Euler--Lagrange transformation reduces the free-boundary problem to a
	fixed-boundary problem. The transformation collapses the airfoil to a slit,
	whose two boundary faces coincide geometrically while retaining separate
	traces. To treat this degenerate geometry, they solve elliptic problems
	separately in the regions above and below the airfoil and the contact curves.
	The remaining contact discontinuity conditions are then matched by the implicit
	function theorem.
	
	In this paper, we consider the corresponding airfoil problem in a
	semi-infinite nozzle. This setting is motivated by wind-tunnel experiments
	used to assess airfoil performance, in which the surrounding walls
	influence the flow and the measured aerodynamic characteristics
	\cite[Sections~9.4 and~13.2]{BarlowRaePope}.
	The analyses in \cite{CXZ,ChenAir} concern piecewise constant background
	states in domains without confining walls. Their matching arguments use
	reflection and rescaling of the Lagrangian quadrants or half-planes.
	In the nozzle considered here, the two flow regions may have unequal widths, and
	the inlet density and horizontal velocity may vary transversely without a
	small-oscillation assumption. The resulting variable coefficients and wall
	conditions prevent a direct application of those matching arguments.
	We therefore treat the elliptic problem on the full slit domain,
	combining suitably chosen Lipschitz approximations with local comparison
	regions and barrier functions adapted to the slit geometry and the
	contact discontinuity conditions. These constructions yield weighted H\"older
	estimates independent of the thickness of the excised tubes, preserving
	$C^\alpha$ control at the trailing edge as the slit is recovered.
	Passing to the limit then gives the desired flow.
	
	The rest of this section states the problem and the main theorem.
	Section~\ref{sec:lagrangian} uses Euler--Lagrange coordinates to
	reformulate the problem as a nonlinear elliptic system. The linear estimates on the
	approximating domains and the recovery of the pressure are given in
	Section~\ref{sec:analysis}. Section~\ref{sec:finite-nonlinear} establishes
	existence and uniqueness on finite nozzles. Sections~\ref{sec:infinite}
	and \ref{sec:physical-return} treat the passage to the semi-infinite
	nozzle and the recovery of the flow in physical coordinates, respectively.
	
	\subsection{Statement of the problem}
	\label{subsec:physical-problem}
	
	Let $W^\pm:[0,\infty)\to\mathbb R$ describe the upper and lower
	walls, respectively, with
	\begin{equation}\label{eq:wall-ordering}
		W^-(x_1)<W^+(x_1),\qquad x_1\ge0.
	\end{equation}
	At the entrance, we assume $W^\pm(0)=H^\pm$, where
	$H^-<0<H^+$ are prescribed constants. The semi-infinite nozzle is
	then defined by
	\[
	\mathcal N=\{x=(x_1,x_2)\in\mathbb R^2:x_1>0,
	W^-(x_1)<x_2<W^+(x_1)\}.
	\]
	We denote its upper and lower walls by
	\[
	\Sigma_{\mathrm w}^\pm
	:=\{(x_1,W^\pm(x_1)):x_1>0\}.
	\]
	
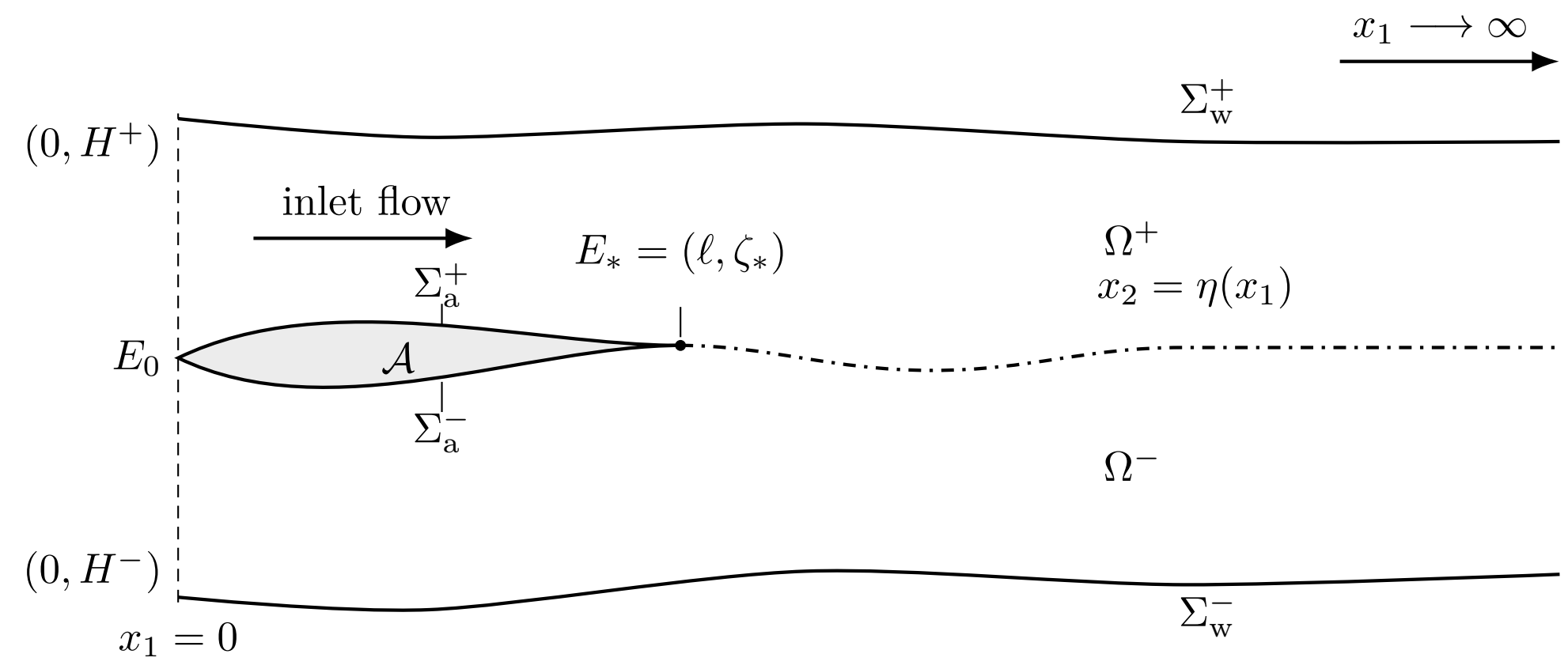
\begin{figure}[htbp]
		\centering
		\begin{tikzpicture}[x=1.05cm,y=0.87cm,>=Latex,font=\small]
			\draw[thick,smooth] plot coordinates {
				(0,2.30)
				(2,2.12)
				(5,2.25)
				(8,2.08)
				(11,2.08)
			};
			\draw[thick,smooth] plot coordinates {
				(0,-2.30)
				(2,-2.42)
				(5,-2.05)
				(8,-2.18)
				(11,-2.08)
			};
			\draw[densely dashed,smooth] (0,-2.35)--(0,2.35);
			
			\filldraw[fill=gray!15,draw=black,thick]
			(0,0) .. controls (1.2,0.70) and (2.8,0.12) .. (4,0.12)
			.. controls (2.8,0.12) and (1.2,-0.70) .. (0,0) --cycle;
			
			\draw[thick,dash dot] (4,0.12)
			.. controls (4.7,0.12) and (5.2,-0.12) .. (6,-0.12)
			.. controls (6.8,-0.12) and (7.2,0.10) .. (8,0.10)
			.. controls (9,0.10) and (10,0.10) .. (11,0.10);
			\fill (4,0.12) circle (1.3pt);
			\node[above] at (4,0.65) {$E_*=(\ell,\zeta_*)$};
			\draw[thin] (4,0.49)--(4,0.20);
			\node at (1.75,0) {$\mathcal A$};
			\node[left] at (0,0) {$E_0$};
			\node[left] at (0,2.05) {$(0,H^+)$};
			\node[left] at (0,-2.05) {$(0,H^-)$};
			\node at (2.1,0.70) {$\Sigma_{\mathrm a}^+$};
			\draw[thin] (2.1,0.52)--(2.1,0.32);
			\node at (2.1,-0.70) {$\Sigma_{\mathrm a}^-$};
			\draw[thin] (2.1,-0.52)--(2.1,-0.23);
			\node[below] at (0,-2.40) {$x_1=0$};
			\node at (7.6,1.15) {$\Omega^+$};
			\node at (7.6,-1.02) {$\Omega^-$};
			\node[above] at (8.1,0.32) {$x_2=\eta(x_1)$};
			\node[above] at (8.2,2.1) {$\Sigma_{\mathrm w}^+$};
			\node[below] at (8.2,-2.1) {$\Sigma_{\mathrm w}^-$};
			\draw[->,thick] (0.6,1.15)--(2.35,1.15);
			\node[above] at (1.5,1.22) {inlet flow};
			\draw[->,thick] (9.25,2.85)--(11,2.85);
			\node[above] at (10.05,2.87) {$x_1\longrightarrow\infty$};
		\end{tikzpicture}
		\caption{The nozzle, airfoil and contact discontinuity.}
		\label{fig:physical-nozzle-airfoil}
	\end{figure}

	An airfoil $\mathcal A$ of length $\ell>0$ is attached to the entrance.
	Its upper and lower surfaces are described by profiles
	$\zeta^\pm:[0,\ell]\to\mathbb R$:
	\[
	\Sigma_{\mathrm a}^\pm
	:=\{(x_1,\zeta^\pm(x_1)):0<x_1<\ell\}.
	\]
	We assume $\zeta^-<\zeta^+$ in $(0,\ell)$, so the airfoil occupies
	\[
	\mathcal A=\{x\in\mathbb R^2:0<x_1<\ell,
	\zeta^-(x_1)<x_2<\zeta^+(x_1)\}.
	\]
	The profiles meet at the leading edge $E_0=(0,0)$ and the trailing
	edge $E_*=(\ell,\zeta_*)$:
	\begin{equation}\label{eq:airfoil-endpoints}
		\zeta^+(0)=\zeta^-(0)=0,\qquad
		\zeta^+(\ell)=\zeta^-(\ell)=\zeta_*.
	\end{equation}
	At $E_*$, we impose the cusp condition
	\begin{equation}\label{eq:airfoil-cusp}
		(\zeta^+)'(\ell)=(\zeta^-)'(\ell)=k,
	\end{equation}
	where $k$ is a constant. The airfoil lies strictly between the nozzle
	walls:
	\begin{equation}\label{eq:airfoil-placement}
		W^-(x_1)<\zeta^-(x_1)\le\zeta^+(x_1)<W^+(x_1),
		\qquad 0\le x_1\le\ell.
	\end{equation}
	The fluid occupies
	\[
	\Omega=\mathcal N\setminus\overline{\mathcal A}.
	\]
	The contact discontinuity emanates from $E_*$ along the unknown curve
	\[
	\Gamma_\eta=\{(x_1,\eta(x_1)):x_1>\ell\},\qquad
	\eta(\ell)=\zeta_*.
	\]
	Together, the airfoil and $\Gamma_\eta$ separate $\Omega$ into an upper region
	$\Omega^+$ and a lower region $\Omega^-$, as shown in
	Figure~\ref{fig:physical-nozzle-airfoil}.
	
	\FloatBarrier
	
	The flow in $\Omega$ is governed by the two-dimensional steady
	compressible Euler equations
	\begin{equation}\label{eq:physical-Euler}
		\begin{cases}
			\operatorname{div}_x(\rho u)=0,\\
			\operatorname{div}_x(\rho u\otimes u)+\nabla_xp=0,\\
			\operatorname{div}_x\bigl(\rho u(E+p/\rho)\bigr)=0.
		\end{cases}
	\end{equation}
	Here $u=(u_1,u_2)$ is the velocity, and $\rho>0$ and $p>0$ are
	the density and pressure. For a polytropic gas with adiabatic exponent
	$\gamma>1$, the energy, sound speed and Mach number are
	\[
	E=\frac{|u|^2}{2}+\frac{p}{(\gamma-1)\rho},\qquad
	c=\sqrt{\frac{\gamma p}{\rho}},\qquad \mathrm{Ma}=\frac{|u|}{c}.
	\]
	We require the flow variables to extend continuously to the boundary
	of each $\Omega^\pm$ and denote their restrictions by
	$u^\pm,\rho^\pm,p^\pm$, respectively.
	The slip conditions on the nozzle walls and the airfoil surfaces read
	\begin{equation}\label{eq:physical-slip}
		u^\pm\cdot n=0\quad\text{on }
		\Sigma_{\mathrm w}^\pm\cup\Sigma_{\mathrm a}^\pm,
	\end{equation}
	where $n$ denotes a unit normal to the corresponding boundary curve.
	Boundary values of $u^\pm,\rho^\pm,p^\pm$ are taken from $\Omega^\pm$.
	Along the free boundary $\Gamma_\eta$, we impose the contact discontinuity
	conditions
	\begin{equation}\label{eq:physical-contact}
		u^\pm\cdot n=0,\qquad [p]=0.
	\end{equation}
	Here $[v]=v^+-v^-$ denotes the jump of $v$ across the curve.
	The derivation of
	\eqref{eq:physical-contact} from the Rankine--Hugoniot conditions can
	be found in \cite{ChenAir,CXZ}.
	
	At the entrance, we distinguish the upper and lower intervals
	\[
	I_0^+=(0,H^+),\qquad I_0^-=(H^-,0),\qquad
	I_0=I_0^+\cup I_0^-.
	\]
	On these intervals, the density and velocity are prescribed by
	\begin{equation}\label{eq:physical-inlet}
		(\rho,u_1,u_2)(0,x_2)
		=\bigl(\rho_{\mathrm{in}},u_{1,\mathrm{in}},u_{2,\mathrm{in}}\bigr)(x_2),
		\qquad x_2\in I_0.
	\end{equation}
	At downstream infinity, the pressure approaches a prescribed constant
	$p_0>0$ uniformly across the nozzle:
	\begin{equation}\label{eq:physical-infinity}
		p(x_1,\cdot)\longrightarrow p_0\qquad(x_1\to\infty).
	\end{equation}
	We assume that the inlet data and the downstream pressure $p_0$ satisfy
	\begin{equation}\label{eq:physical-background}
		\rho_{\mathrm{in}}\ge\rho_*>0,\qquad
		u_{1,\mathrm{in}}\ge u_*>0,\qquad
		\frac{\rho_{\mathrm{in}}u_{1,\mathrm{in}}^2}{\gamma p_0}\le1-\delta_b
	\end{equation}
	in $I_0$, where $\rho_*,u_*>0$ and $\delta_b\in(0,1)$ are fixed constants.
	The inlet density and horizontal velocity determine the positive mass
	fluxes in the two flow regions:
	\begin{equation}\label{eq:mass-inlet-maps}
		m_+=\int_0^{H^+}(\rho_{\mathrm{in}}u_{1,\mathrm{in}})(\xi)\,d\xi,
		\qquad
		m_-=\int_{H^-}^0(\rho_{\mathrm{in}}u_{1,\mathrm{in}})(\xi)\,d\xi.
	\end{equation}
	To match the inlet data with the adjoining solid boundaries, we require
	\begin{equation}\label{eq:physical-compatibility}
		\frac{u_{2,\mathrm{in}}(H^\pm)}{u_{1,\mathrm{in}}(H^\pm)}
		=(W^\pm)'(0),\qquad
		\frac{u_{2,\mathrm{in}}^\pm(0)}{u_{1,\mathrm{in}}^\pm(0)}
		=(\zeta^\pm)'(0).
	\end{equation}
	
	With these preparations, we formulate the free boundary problem as follows.
	
	\begin{mainproblem}[Subsonic flow past an airfoil]\label{prob:physical}
		Suppose that the wall profiles $W^\pm$ and the airfoil profiles $\zeta^\pm$
		satisfy \eqref{eq:wall-ordering} and
		\eqref{eq:airfoil-endpoints}--\eqref{eq:airfoil-placement},
		and that \eqref{eq:physical-background} and
		\eqref{eq:physical-compatibility} hold.
		Find a flow $U=(\rho,u_1,u_2,p)$ and a continuously differentiable
		contact graph $\eta$ such that:
		\begin{enumerate}
			\item $U$ solves \eqref{eq:physical-Euler} classically in
			$\Omega^+$ and $\Omega^-$,
			with the slip and contact discontinuity conditions
			\eqref{eq:physical-slip}--\eqref{eq:physical-contact}, the inlet
			condition \eqref{eq:physical-inlet}, and the downstream pressure
			condition \eqref{eq:physical-infinity}.
			\item The flow is uniformly subsonic and has uniformly positive horizontal mass flux:
			\begin{equation}\label{eq:main-flow-admissibility}
				\sup_{\Omega^+\cup\Omega^-}\mathrm{Ma}<1,\qquad
				\inf_{\Omega^+\cup\Omega^-}(\rho u_1)>0.
			\end{equation}
			\item The contact curve joins the trailing edge tangentially and
			lies strictly inside the nozzle:
			\begin{equation}\label{eq:main-contact-location}
				\begin{gathered}
					\eta(\ell)=\zeta_*,\qquad \eta'(\ell)=k,\\
					W^-(x_1)<\eta(x_1)<W^+(x_1),\qquad x_1\ge\ell.
				\end{gathered}
			\end{equation}
		\end{enumerate}
	\end{mainproblem}
	\setcounter{problem}{1}
	
	\subsection{Function spaces and the main result}
	\label{subsec:physical-main-result}
	
	To state the main result, we introduce weighted H\"older norms that
	describe the regularity at the trailing edge and the downstream
	behavior \cite{ChenAir,CXZ}.
	Let $G\subset\{z\in\mathbb R^n:z_1\ge0\}$ be a domain, with
	$n=1$ or $2$, and let $\mathcal E$ be a
	nonempty finite subset of $\mathbb R^n\setminus G$. Set
	\[
	d_{\mathcal E}(z)=\min\{1,\operatorname{dist}(z,\mathcal E)\},\qquad
	d_{\mathcal E}(z,z')=\min\{d_{\mathcal E}(z),d_{\mathcal E}(z')\}
	\]
	and
	\[
	\Delta_z=1+z_1,\qquad
	\Delta_{z,z'}=\min\{\Delta_z,\Delta_{z'}\}.
	\]
	For $m\in\mathbb N_0$, $0<\alpha<1$, $s\in\mathbb R$, and
	$\boldsymbol\beta=(\beta_0,\ldots,\beta_m)\in\mathbb R^{m+1}$, define
	\begin{equation}\label{eq:general-X-definition}
		\begin{split}
			\|v\|_{m,\alpha;(\boldsymbol\beta);G}^{(s;\mathcal E)}
			:={}&\sum_{j=0}^m\sup_{z\in G}
			d_{\mathcal E}(z)^{\max\{j+s,0\}}
			\Delta_z^{\beta_j}|D^jv(z)|\\
			&+\sup_{\substack{z,z'\in G\\z\ne z'}}
			d_{\mathcal E}(z,z')^{\max\{m+\alpha+s,0\}}
			\Delta_{z,z'}^{\beta_m}
			\frac{|D^mv(z)-D^mv(z')|}{|z-z'|^\alpha}.
		\end{split}
	\end{equation}
	Here $\mathbb N_0=\{0,1,2,\ldots\}$, and $D^jv$ denotes the array
	of derivatives of order $j$, with $D^0v=v$. The space
	$C^{m,\alpha;(\boldsymbol\beta)}_{(s;\mathcal E)}(G)$ consists of functions
	that are locally $C^{m,\alpha}$ up to the boundary away from
	$\mathcal E$ and have finite norm in \eqref{eq:general-X-definition}.
	A scalar weight $\beta$ means $\beta_0=\cdots=\beta_m=\beta$. Omitting
	$(\boldsymbol\beta)$ or $(s;\mathcal E)$ removes the corresponding
	weight factors. In the applications below, we omit $\mathcal E$
	from the norm notation; it consists of the trailing edge in the
	coordinates under consideration. Any other weight set is indicated
	explicitly. With both weights omitted, $\|v\|_{m,\alpha;G}$ is the
	usual H\"older norm. For $\alpha=0$, the difference-quotient term
	is omitted, so $\|v\|_{0,0;G}=\sup_G|v|$.
	
	For the flow variables and coefficients, we use
	\[
	\|v\|_{1,\alpha;G}^{(-\alpha)},\qquad
	\|v\|_{1,\alpha;(\beta);G}^{(-\alpha)},\qquad
	\|v\|_{1,\alpha;(\beta-1,\beta);G}^{(-\alpha)}.
	\]
	These norms have the same corner weights. The last two differ only
	in the weight on the function value. For general $\beta_0,\beta_1$, we have
	\[
	\|v\|_{1,\alpha;(\beta_0,\beta_1);G}^{(-\alpha)}
	\le\sup_G\Delta^{\beta_0}|v|
	+\|Dv\|_{0,\alpha;(\beta_1);G}^{(1-\alpha)}
	\le C\|v\|_{1,\alpha;(\beta_0,\beta_1);G}^{(-\alpha)}.
	\]
	Here $C$ is a positive constant depending only on the convention for
	derivative arrays.
	
	On a smooth boundary segment $J$, the same definition uses tangential
	derivatives and distance in the ambient space. With an arclength
	parametrization $\gamma_J(s)$, we set
	$D^jv=d^j(v\circ\gamma_J)/ds^j$. In particular,
	$\partial_{\mathrm{tan}}v=Dv\cdot\tau_J$, where
	$\tau_J=\gamma_J'$, whenever $v$ extends differentiably off $J$.
	On intervals in $\mathbb R$, derivatives and
	distances are one-dimensional. For example, the space
	$C^{1,\alpha}_{(-\alpha;\{\ell\})}((0,\ell))$ uses
	$d_{\{\ell\}}(x_1)=\min\{1,\ell-x_1\}$ and requires ordinary
	$C^{1,\alpha}$ regularity at $x_1=0$.
	On $(\ell,\infty)$, we use
	$d_{\{\ell\}}(x_1)=\min\{1,x_1-\ell\}$ and $\Delta_{x_1}=1+x_1$.
	For the contact graph, the norm
	$\|v\|_{2,\alpha;(\beta-1,\beta,\beta);(\ell,\infty)}^{(-1-\alpha;\{\ell\})}$
	measures $v$ with downstream weight $\beta-1$ and its derivatives
	with weight $\beta$, while allowing $v''$ to grow like
	$(x_1-\ell)^{\alpha-1}$ near the trailing edge.
	
	For vector- and matrix-valued functions, we sum the component norms.
	On two subdomains $G^+,G^-$, the subscript $\mathrm{pm}$ denotes the piecewise
	norm
	\[
	\|v\|_{m,\alpha;(\boldsymbol\beta);G,\mathrm{pm}}^{(s;\mathcal E)}
	=\sum_{\sigma\in\{+,-\}}
	\|v^\sigma\|_{m,\alpha;(\boldsymbol\beta);G^\sigma}^{(s;\mathcal E)}.
	\]
	The corresponding product space is denoted by
	$C_{(s;\mathcal E),\mathrm{pm}}^{m,\alpha;(\boldsymbol\beta)}(G)$.
	The same piecewise convention applies to unweighted H\"older and
	Sobolev spaces. In $C^0_{\mathrm{pm}}$, the restrictions to $G^+$ and $G^-$
	have bounded continuous extensions to $\overline{G^+}$ and
	$\overline{G^-}$, respectively. For the domains considered below, a piecewise $W^{1,q}$ function whose
	traces agree on the interior contact line belongs to $W^{1,q}(G)$,
	$1\le q<\infty$. The two faces of a slit retain separate traces. We use
	$H_0^1(G)=\overline{C_c^\infty(G)}^{\,H^1(G)}$ on all domains,
	including slit domains. In the original coordinates,
	$\mathcal E=\{E_*\}$. After a coordinate change, $\mathcal E$ denotes
	the image of this set.
	
	For local Schauder estimates we use the standard scaled norm
	\cite{GT,CD}:
	\begin{equation}\label{eq:GT-scaled-norm}
		\|v\|_{m,\alpha;G}'
		=\sum_{j=0}^m d_G^{\,j}\|D^jv\|_{L^\infty(G)}
		+d_G^{\,m+\alpha}[D^mv]_{\alpha;G},\qquad
		d_G=\operatorname{diam}G,
	\end{equation}
	where $[\,\cdot\,]_{\alpha;G}$ is the usual H\"older seminorm.
	We also write $B_r(z)=\{z':|z'-z|<r\}$ and $B_r=B_r(0)$.
	For a measurable set $G$ with $0<|G|<\infty$, $|G|$ denotes its
	Lebesgue measure, and the integral average is
	\[
	\fint_G v\,dz=\frac1{|G|}\int_G v\,dz.
	\]
	Throughout the paper, $C$ denotes a positive constant that may
	change from line to line.
	
	Fix arbitrary $\alpha\in(0,\tfrac12)$ and $\beta>1$.
	We assume that the inlet data and the solid-boundary profiles have
	the following regularity:
	\begin{equation}\label{eq:physical-data-regularity}
		\begin{gathered}
			\rho_{\mathrm{in}},u_{1,\mathrm{in}},u_{2,\mathrm{in}}
			\in C^{1,\alpha}_{\mathrm{pm}}(I_0),\\
			(W^\pm)'\in C^{1,\alpha;(\beta)}([0,\infty)),\qquad
			(\zeta^\pm)'\in C^{1,\alpha}_{(-\alpha;\{\ell\})}((0,\ell)).
		\end{gathered}
	\end{equation}
	Using these norms, we measure the transverse inlet velocity and the
	wall and airfoil slopes by
	\begin{equation}\label{eq:physical-kappa}
		\kappa:=\|u_{2,\mathrm{in}}\|_{1,\alpha;I_0,\mathrm{pm}}
		+\sum_{\sigma\in\{+,-\}}\left(
		\|(W^\sigma)'\|_{1,\alpha;(\beta);[0,\infty)}
		+\|(\zeta^\sigma)'\|_{1,\alpha;(0,\ell)}^{(-\alpha;\{\ell\})}\right).
	\end{equation}
	\Needspace{10\baselineskip}
	\begin{theorem}[Existence and uniqueness of subsonic flows]\label{thm:physical-main}
		There exist positive constants $\kappa_0,C$ such that, if
		$\kappa\le\kappa_0$, the \hyperref[prob:physical]{Main Problem} has
		a solution $(U,\eta)$ with
		\[
		U\in C^{1,\alpha}_{(-\alpha;\{E_*\}),\mathrm{pm}}(\Omega),\qquad
		\eta\in C^{2,\alpha}_{(-1-\alpha;\{\ell\})}((\ell,\infty)),
		\]
		and a constant $\eta_\infty$ such that
		\begin{equation}\label{eq:physical-main-estimate}
			\|u_2\|_{1,\alpha;(\beta);\Omega,\mathrm{pm}}^{(-\alpha)}
			+\|p-p_0\|_{1,\alpha;(\beta-1,\beta);\Omega,\mathrm{pm}}^{(-\alpha)}
			+\|\eta-\eta_\infty\|_{2,\alpha;(\beta-1,\beta,\beta);(\ell,\infty)}
			^{(-1-\alpha;\{\ell\})}\le C\kappa.
		\end{equation}
		The solution is unique among solutions of the
		\hyperref[prob:physical]{Main Problem} satisfying
		\eqref{eq:physical-main-estimate} with the same constant $C$.
		The constants depend only on the inlet profiles
		$\rho_{\mathrm{in}},u_{1,\mathrm{in}}$ and on
		$p_0,\ell,H^\pm,\gamma,\alpha,\beta$.
	\end{theorem}
	
	\begin{remark}
		In Lagrangian coordinates, $\rho$ and $u_1$ satisfy weighted H\"older estimates
		analogous to those for $p$, relative to their downstream limits and with the
		weights specified in \eqref{eq:downstream-state-estimate}. Their regularity and
		uniqueness follow from those of $(w,p)$ and the prescribed inlet data
		through the algebraic relations \eqref{eq:algebraic-flow-variables}.
	\end{remark}
	
	\section{Reformulation of the Main Problem}
	\label{sec:lagrangian}
	
	We reformulate the \hyperref[prob:physical]{Main Problem} in
	Euler--Lagrange coordinates, following S.-X.~Chen~\cite{ChenMach}
	and J.~Chen~\cite{ChenAir}.
	
	\subsection{The Euler--Lagrange transformation}
	\label{subsec:Euler-Lagrange}
	
	Let $U$ be a piecewise smooth subsonic solution of the
	\hyperref[prob:physical]{Main Problem}. In each simply connected
	region $\Omega^\pm$, the mass equation determines a stream function
	$\psi$, up to an additive constant, by
	\begin{equation}\label{eq:stream-function}
		\psi_{x_1}=-\rho u_2,\qquad \psi_{x_2}=\rho u_1.
	\end{equation}
	We choose the additive constants so that $\psi^\pm(E_0)=0$ and introduce the
	Euler--Lagrange transformation
	\begin{equation}\label{eq:Euler-Lagrange-map}
		y=\mathcal T(x),\qquad
		\begin{cases}
			y_1=x_1,\\
			y_2=\psi(x_1,x_2).
		\end{cases}
	\end{equation}
	By the slip and contact discontinuity conditions
	\eqref{eq:physical-slip}--\eqref{eq:physical-contact}, the nozzle walls,
	airfoil boundaries and contact curve are streamlines. The choice
	$\psi^\pm(E_0)=0$ and the mass fluxes defined in
	\eqref{eq:mass-inlet-maps} therefore determine their images:
	the nozzle walls are mapped to
	\[
	\mathcal T(\Sigma_{\mathrm w}^\pm)
	=\Gamma_{\mathrm w,\infty}^\pm:=(0,\infty)\times\{\pm m_\pm\},
	\]
	while the upper and lower airfoil boundaries are mapped to
	\[
	\mathcal T(\Sigma_{\mathrm a}^\pm)
	=\Gamma_{\mathrm a}^\pm:=(0,\ell)\times\{0\}.
	\]
	These two airfoil images form the opposite sides of the slit
	\[
	S_*=(0,\ell]\times\{0\}.
	\]
	Accordingly, we introduce the transverse intervals
	\[
	I^+=(0,m_+),\qquad I^-=(-m_-,0),\qquad I=I^+\cup I^-.
	\]
	The two flow regions are then mapped to
	\[
	\mathcal T(\Omega^\pm)=D_\infty^\pm:=(0,\infty)\times I^\pm.
	\]
	Including the contact line, we obtain the fixed slit domain
	\begin{equation}\label{eq:infinite-domain}
		\mathcal T(\Omega)=D_\infty
		=\bigl((0,\infty)\times(-m_-,m_+)\bigr)\setminus S_*.
	\end{equation}
	The contact curve itself is mapped to
	\[
	\mathcal T(\Gamma_\eta)=\Gamma_\infty^*
	=(\ell,\infty)\times\{0\},
	\]
	and the trailing edge becomes $e_*=(\ell,0)$.
	
	Since $y_1=x_1$ and $\psi_{x_2}=\rho u_1>0$, the transformation is
	invertible in $\Omega^+$ and $\Omega^-$. We retain the notation $U(y)$ for
	$U(\mathcal T^{-1}(y))$. Using \eqref{eq:stream-function} to express
	the derivatives in the $y$-coordinates, we rewrite
	\eqref{eq:physical-Euler} as
	\begin{subequations}\label{eq:lagrangian-conservation}
		\begin{align}
			\partial_{y_1}\!\left(\frac1{\rho u_1}\right)
			-\partial_{y_2}\!\left(\frac{u_2}{u_1}\right)&=0,
			\label{eq:lagrangian-mass}\\
			\partial_{y_1}\!\left(u_1+\frac p{\rho u_1}\right)
			-\partial_{y_2}\!\left(\frac{pu_2}{u_1}\right)&=0,
			\label{eq:lagrangian-horizontal-momentum}\\
			(u_2)_{y_1}+p_{y_2}&=0,
			\label{eq:lagrangian-vertical-momentum}\\
			\partial_{y_1}\!\left(\frac{|u|^2}{2}
			+\frac{\gamma p}{(\gamma-1)\rho}\right)&=0.
			\label{eq:lagrangian-energy}
		\end{align}
	\end{subequations}
	
	\subsection{Reduction to an elliptic problem}
	\label{subsec:elliptic-reduction}
	\label{subsec:lagrangian-problem}
	
	We formulate the elliptic problem in terms of the pressure $p$ and
	the flow direction
	\[
	w=\frac{u_2}{u_1}.
	\]
	In terms of $w$, the slip conditions on the nozzle walls and
	the airfoil give
	\begin{equation}\label{eq:lagrangian-slip-data}
		\begin{cases}
			w(y_1,m_+)=(W^+)'(y_1),&y_1\ge0,\\
			w(y_1,-m_-)=(W^-)'(y_1),&y_1\ge0,\\
			w^\pm(y_1,0)=(\zeta^\pm)'(y_1),&0<y_1<\ell.
		\end{cases}
	\end{equation}
	To express the inlet condition in the same coordinates, we use
	\begin{equation}\label{eq:inverse-inlet-map}
		y_2=\psi(0,x_2)=\int_0^{x_2}
		(\rho_{\mathrm{in}}u_{1,\mathrm{in}})(\xi)\,d\xi,
		\qquad x_2\in\overline{I_0^\pm}.
	\end{equation}
	The positivity of $\rho_{\mathrm{in}}u_{1,\mathrm{in}}$ makes this map
	strictly increasing on each $\overline{I_0^\pm}$. Let
	$s_\pm:\overline{I^\pm}\to\overline{I_0^\pm}$ denote its inverse
	on the corresponding interval, so that $x_2=s_\pm(y_2)$.
	Writing $s=s_\pm$ on $I^\pm$, we set
	\[
	\widehat\rho_{\mathrm{in}}=\rho_{\mathrm{in}}\circ s,\qquad
	\widehat u_{i,\mathrm{in}}=u_{i,\mathrm{in}}\circ s,\qquad i=1,2,
	\]
	so the prescribed inlet velocity gives
	\begin{equation}\label{eq:lagrangian-inlet-angle}
		w(0,y_2)=w_{\mathrm{in}}(y_2)
		:=\frac{\widehat u_{2,\mathrm{in}}(y_2)}
		{\widehat u_{1,\mathrm{in}}(y_2)},\qquad y_2\in I.
	\end{equation}
	
	On the fixed contact line, \eqref{eq:physical-contact} implies
	\begin{equation}\label{eq:lagrangian-contact-data}
		[w]=[p]=0\quad\text{on }\Gamma_\infty^*.
	\end{equation}
	Since $y_1=x_1$, the prescribed downstream pressure gives
	\begin{equation}\label{eq:lagrangian-pressure-infinity}
		p(y_1,\cdot)\longrightarrow p_0\qquad(y_1\to\infty).
	\end{equation}
	
	To derive the equations for $(w,p)$, we write the Lagrangian
	system~\eqref{eq:lagrangian-conservation} in nondivergence form:
	\begin{equation}\label{eq:characteristic-matrix-form}
		\mathsf M(U)U_{y_1}+\mathsf N(U)U_{y_2}=0,
	\end{equation}
	where
	\[
	\mathsf M=
	\begin{pmatrix}
		-\dfrac1{\rho^2u_1}&-\dfrac1{\rho u_1^2}&0&0\\[2mm]
		-\dfrac p{\rho^2u_1}&1-\dfrac p{\rho u_1^2}&0&\dfrac1{\rho u_1}\\[2mm]
		0&0&1&0\\[1mm]
		-\dfrac{\gamma p}{(\gamma-1)\rho^2}&u_1&u_2&
		\dfrac\gamma{(\gamma-1)\rho}
	\end{pmatrix},\qquad
	\mathsf N=
	\begin{pmatrix}
		0&\dfrac{u_2}{u_1^2}&-\dfrac1{u_1}&0\\[2mm]
		0&\dfrac{pu_2}{u_1^2}&-\dfrac p{u_1}&-\dfrac{u_2}{u_1}\\[2mm]
		0&0&0&1\\[1mm]
		0&0&0&0
	\end{pmatrix}.
	\]
	The characteristic equation $\det(\lambda\mathsf M-\mathsf N)=0$
	has four roots:
	\[
	\lambda_1=\lambda_2=0,\qquad
	\lambda_\pm=-\frac{c\rho}{c^2-u_1^2}
	\left(cu_2\mp i u_1\sqrt{c^2-q^2}\right),
	\]
	where $q=(u_1^2+u_2^2)^{1/2}<c$. Corresponding left eigenvectors are
	\[
	\begin{aligned}
		\boldsymbol l_1={}&(0,0,0,1),\\
		\boldsymbol l_2={}&(-pu_1,u_1,u_2,-1),\\
		\boldsymbol l_\pm={}&\Biggl(
		\left(\frac{\gamma p^2}{(\gamma-1)\rho u_1}
		-\frac{pu_1}{\gamma-1}\right)\lambda_\pm
		+\frac{\gamma p^2u_2}{(\gamma-1)u_1},\\
		&\quad-\left(u_1+\frac{\gamma p}{(\gamma-1)\rho u_1}\right)\lambda_\pm
		-\frac{\gamma pu_2}{(\gamma-1)u_1},\quad
		\frac{\gamma p}{\gamma-1}-u_2\lambda_\pm,\quad
		\lambda_\pm\Biggr).
	\end{aligned}
	\]
	Multiplying \eqref{eq:characteristic-matrix-form} by
	$\boldsymbol l_1$ recovers \eqref{eq:lagrangian-energy} and gives
	Bernoulli's law:
	\begin{subequations}\label{eq:transported-profiles}
		\begin{equation}\label{eq:Bernoulli-law}
			\frac{q^2}{2}+\frac{\gamma p}{(\gamma-1)\rho}=B(y_2).
		\end{equation}
		Multiplication by $\boldsymbol l_2$ gives
		$\partial_{y_1}(p/\rho^\gamma)=0$, and hence
		\begin{equation}\label{eq:entropy-law}
			\frac p{\rho^\gamma}=S(y_2).
		\end{equation}
	\end{subequations}
	These identities show that $S$ and $B$ are constant along each
	streamline $y_2=\mathrm{const}$. To obtain the equations for the flow
	direction and pressure, we multiply \eqref{eq:characteristic-matrix-form}
	by $\boldsymbol l_+$ and take real and imaginary parts:
	\begin{subequations}\label{eq:characteristic-wp-system}
		\begin{align}
			D_Rw+\mu D_Ip&=0,\label{eq:characteristic-wp-real}\\*
			D_Iw-\mu D_Rp&=0,\label{eq:characteristic-wp-imaginary}
		\end{align}
	\end{subequations}
	where
	\[
	D_R=\partial_{y_1}+\lambda_R\partial_{y_2},\qquad
	D_I=\lambda_I\partial_{y_2},
	\]
	and
	\begin{equation}\label{eq:characteristic-real-coefficients}
		\mu=\frac{\sqrt{c^2-q^2}}{c\rho u_1^2},\qquad
		\lambda_R=-\frac{c^2\rho u_2}{c^2-u_1^2},\qquad
		\lambda_I=\frac{c\rho u_1\sqrt{c^2-q^2}}{c^2-u_1^2}>0.
	\end{equation}
	
	Since $\mu\lambda_I>0$, \eqref{eq:characteristic-wp-system} is
	equivalent to
	\begin{equation}\label{eq:nonlinear-rotated-system}
		\nabla_y^\perp p=A\nabla_yw,
	\end{equation}
	where $\nabla_y^\perp v=(-v_{y_2},v_{y_1})^{\mathsf T}$, the superscript
	$\mathsf T$ denotes transpose, and
	\begin{equation}\label{eq:characteristic-coefficient-matrix}
		A=\frac1{\mu\lambda_I}
		\begin{pmatrix}1&\lambda_R\\
			\lambda_R&\lambda_R^2+\lambda_I^2\end{pmatrix}.
	\end{equation}
	The matrix $A$ is symmetric and, by \eqref{eq:main-flow-admissibility},
	positive definite. Hence \eqref{eq:nonlinear-rotated-system} is a
	first-order elliptic system.
	
	To determine the dependence of $A$ on the unknown pair $(w,p)$,
	we use \eqref{eq:transported-profiles} to express the remaining flow
	variables in terms of $(S,B,w,p)$:
	\begin{equation}\label{eq:algebraic-flow-variables}
		\rho=\left(\frac pS\right)^{1/\gamma},\qquad
		u_1=\frac{\sqrt{2B-\dfrac{2\gamma p}{(\gamma-1)\rho}}}
		{\sqrt{1+w^2}},\qquad u_2=wu_1.
	\end{equation}
	In particular, the sound speed and the flow speed satisfy
	\begin{equation}\label{eq:coefficient-functions-SB}
		c=\sqrt{\gamma S^{1/\gamma}p^{(\gamma-1)/\gamma}},\qquad
		q=\sqrt{2B-\frac{2c^2}{\gamma-1}}.
	\end{equation}
	Substituting $u_1=q/\sqrt{1+w^2}$ and $\rho c^2=\gamma p$ into
	\eqref{eq:characteristic-coefficient-matrix} gives
	\begin{equation}\label{eq:chen-matrix}
		A=(a^{ij})=\frac{q}{\sqrt{1+w^2}(c^2-q^2)}
		\begin{pmatrix}
			c^2-\dfrac{q^2}{1+w^2}
			&-\dfrac{\gamma p\,qw}{\sqrt{1+w^2}}\\[3mm]
			-\dfrac{\gamma p\,qw}{\sqrt{1+w^2}}
			&\dfrac{\gamma^2p^2q^2}{c^2}
		\end{pmatrix}.
	\end{equation}
	The matrix $A$ is thus determined by $(S,B,w,p)$, and we
	write $A=A(S,B,w,p)$. Since $\boldsymbol l_1$, $\boldsymbol l_2$,
	$\operatorname{Re}\boldsymbol l_+$ and $\operatorname{Im}\boldsymbol l_+$
	are linearly independent in the subsonic regime, the preceding reduction
	is reversible. Thus the Lagrangian system
	\eqref{eq:lagrangian-conservation} is equivalent to the nonlinear
	elliptic system \eqref{eq:nonlinear-rotated-system} together with
	\eqref{eq:transported-profiles}, with the remaining flow variables
	recovered by \eqref{eq:algebraic-flow-variables}.
	
	The prescribed inlet density and velocity relate the transported profiles $S$ and
	$B$ to the inlet pressure trace. Evaluating~\eqref{eq:transported-profiles} at
	$y_1=0$ gives
	\begin{equation}\label{eq:inlet-transported-profiles}
		\begin{aligned}
			S(y_2)&=\frac{p(0,y_2)}{\widehat\rho_{\mathrm{in}}(y_2)^\gamma},\\
			B(y_2)&=\frac{\widehat u_{1,\mathrm{in}}(y_2)^2+
				\widehat u_{2,\mathrm{in}}(y_2)^2}{2}
			+\frac{\gamma p(0,y_2)}{(\gamma-1)\widehat\rho_{\mathrm{in}}(y_2)}.
		\end{aligned}
	\end{equation}
	Consequently, with the inlet density and velocity fixed,
	\eqref{eq:coefficient-functions-SB} and \eqref{eq:chen-matrix} show that $A(y)$
	depends on $w(y)$, $p(y)$ and the inlet pressure trace $p(0,y_2)$.
	
	The transformed inlet data and boundary profiles satisfy the
	compatibility conditions inherited from the
	\hyperref[prob:physical]{Main Problem}:
	\begin{equation}\label{eq:lagrangian-compatibility}
		\begin{cases}
			w_{\mathrm{in}}(m_+)=(W^+)'(0),\quad
			w_{\mathrm{in}}(-m_-)=(W^-)'(0),\\
			w_{\mathrm{in}}^+(0)=(\zeta^+)'(0),\quad
			w_{\mathrm{in}}^-(0)=(\zeta^-)'(0),\\
			(\zeta^+)'(\ell)=(\zeta^-)'(\ell)=k.
		\end{cases}
	\end{equation}
	
	These equations and boundary conditions lead to the following
	reformulation of the \hyperref[prob:physical]{Main Problem}.
	\begin{problem}[Reformulated problem]\label{prob:lagrangian}
		Suppose that the data are obtained from the
		\hyperref[prob:physical]{Main Problem} by the above transformation,
		so that \eqref{eq:lagrangian-compatibility} holds.
		Find $w,p$ on $D_\infty$, with $p>0$, such that:
		\begin{enumerate}
			\item The pair solves \eqref{eq:nonlinear-rotated-system} in
			$D_\infty^+$ and $D_\infty^-$, with $A=A(S,B,w,p)$ given by
			\eqref{eq:chen-matrix}, where $S,B$ are determined by
			\eqref{eq:inlet-transported-profiles}.
			\item The slip and inlet conditions \eqref{eq:lagrangian-slip-data}
			and \eqref{eq:lagrangian-inlet-angle}, the contact discontinuity
			conditions \eqref{eq:lagrangian-contact-data},
			and the downstream pressure condition \eqref{eq:lagrangian-pressure-infinity}
			hold.
			\item The flow recovered from $(w,p)$ by
			\eqref{eq:algebraic-flow-variables} is uniformly subsonic and has
			uniformly positive horizontal mass flux:
			\begin{equation}\label{eq:wp-admissibility}
				\sup_{D_\infty}\frac{q^2}{c^2}<1,\qquad
				\inf_{D_\infty}(\rho u_1)>0.
			\end{equation}
		\end{enumerate}
	\end{problem}
	
	In the following sections, we use the weighted norms introduced in
	Section~\ref{sec:physical}, with the trailing edge located at $e_*$:
	\[
	\mathcal E=\{e_*\},\qquad d(y)=d_{\mathcal E}(y),\qquad
	d(y,z)=d_{\mathcal E}(y,z).
	\]
	
	In these coordinates, the parallel reference flow is
	\begin{equation}\label{eq:lagrangian-reference-state}
		\widehat U_0(y_2)=
		(\widehat\rho_{\mathrm{in}}(y_2),\widehat u_{1,\mathrm{in}}(y_2),0,p_0),
		\qquad y_2\in I.
	\end{equation}
	The inequalities in \eqref{eq:physical-background} remain valid with
	$\rho_{\mathrm{in}},u_{1,\mathrm{in}}$ replaced by
	$\widehat\rho_{\mathrm{in}},\widehat u_{1,\mathrm{in}}$.
	The positivity and $C^{1,\alpha}$ regularity of
	$\rho_{\mathrm{in}}u_{1,\mathrm{in}}$ imply
	$s_\pm\in C^{2,\alpha}(\overline{I^\pm})$; composition with $s_\pm$
	then gives $\widehat\rho_{\mathrm{in}},\widehat u_{1,\mathrm{in}},
	\widehat u_{2,\mathrm{in}}\in C^{1,\alpha}_{\mathrm{pm}}(I)$.
	The matrix associated with the reference state $\widehat U_0$ through
	\eqref{eq:chen-matrix} is diagonal:
	$A_b(y_2)=\operatorname{diag}(a_b^{11}(y_2),a_b^{22}(y_2))$.
	It satisfies
	\begin{equation}\label{eq:lagrangian-background-bounds}
		\lambda\le a_b^{11},a_b^{22}\le\Lambda,\qquad
		\|A_b\|_{1,\alpha;I,\mathrm{pm}}\le\Lambda,
	\end{equation}
	for fixed positive constants $\lambda,\Lambda$.
	
	We measure the size of the transformed data relative to this background by
	\begin{equation}\label{eq:lagrangian-data-size}
		\hat\kappa:=\|\widehat u_{2,\mathrm{in}}\|_{1,\alpha;I,\mathrm{pm}}
		+\sum_{\sigma\in\{+,-\}}\left(
		\|(W^\sigma)'\|_{1,\alpha;(\beta);[0,\infty)}
		+\|(\zeta^\sigma)'\|_{1,\alpha;(0,\ell)}^{(-\alpha;\{\ell\})}\right).
	\end{equation}
	
	\begin{theorem}[Solvability of the reformulated problem]\label{thm:infinite-existence}
		Suppose that $\alpha,\beta$ and the inlet data are as in
		Theorem~\ref{thm:physical-main} and that $\widehat U_0$ is given by
		\eqref{eq:lagrangian-reference-state}.
		Then there exist positive constants $\hat\kappa_0,\sigma_0,C$ such that,
		if $\hat\kappa\le\hat\kappa_0$, Problem~\ref{prob:lagrangian} has a solution
		$(w,p)$ in the piecewise weighted $C^{1,\alpha}$ class satisfying
		\begin{equation}\label{eq:lagrangian-main-estimate}
			\|w\|_{1,\alpha;(\beta);D_\infty,\mathrm{pm}}^{(-\alpha)}
			+\|p-p_0\|_{1,\alpha;(\beta-1,\beta);D_\infty,\mathrm{pm}}^{(-\alpha)}
			\le C\hat{\kappa}.
		\end{equation}
		The solution is unique among solutions of Problem~\ref{prob:lagrangian}
		for which the left-hand side of \eqref{eq:lagrangian-main-estimate} is
		at most $\sigma_0$. The constants depend only on
		$\widehat U_0,\ell,m_\pm,\gamma,\alpha,\beta$.
	\end{theorem}
	
	\section{The linear problem in a finite nozzle}
	\label{sec:analysis}
	
	Throughout this section, fix $0<\alpha<\tfrac12$ and $\beta>1$.
	For prescribed $\ell,m_+,m_->0$ and $L>\ell+1$, consider the
	finite slit domain
	\begin{equation}\label{eq:finite-original-domain}
		D_L=\bigl((0,L)\times(-m_-,m_+)\bigr)\setminus S_*.
	\end{equation}
	We distinguish the portions above and below $y_2=0$:
	\[
	D_L^+=(0,L)\times I^+,\qquad D_L^-=(0,L)\times I^-.
	\]
	The contact line in the truncated domain and its endpoint at the outlet are
	\[
	\Gamma_L^*=(\ell,L)\times\{0\},\qquad e_L=(L,0).
	\]
	Restricting the nozzle walls to $0<y_1<L$ gives
	\[
	\Gamma_{\mathrm w}^\pm
	=\Gamma_{\mathrm w,\infty}^\pm\cap\{0<y_1<L\}
	=(0,L)\times\{\pm m_\pm\}.
	\]
	The two sides of the slit retain the notation
	\[
	\Gamma_{\mathrm a}^+=(0,\ell)\times\{0\},\qquad
	\Gamma_{\mathrm a}^-=(0,\ell)\times\{0\}.
	\]
	We denote the inlet segments and the outlet by
	\[
	I_{\mathrm{in}}^\sigma=\{0\}\times I^\sigma\quad(\sigma\in\{+,-\}),\qquad
	O_L=\{L\}\times(-m_-,m_+).
	\]
	The labels $\Gamma_{\mathrm a}^\pm$ refer to the distinct boundary
	traces on the two sides of the slit. We keep $L$ fixed throughout
	this section and Section~\ref{sec:finite-nonlinear}, and omit it from
	the notation for coefficients, data and solutions when there is no ambiguity.
	Unless otherwise indicated, differential operators act on the variable $y$.
	
	For a prescribed symmetric, uniformly elliptic matrix $A=(a^{ij})$,
	we consider the following linear elliptic system:
	\begin{equation}\label{eq:linear-first-order}
		\nabla^\perp\widetilde p=A\nabla \widetilde w
		\quad\text{in }D_L^\pm,
	\end{equation}
	subject to the contact discontinuity conditions
	\begin{equation}\label{eq:finite-linear-contact}
		[\widetilde w]=[\widetilde p]=0\quad\text{on }\Gamma_L^*.
	\end{equation}
	
	We prescribe a boundary function $g$, with traces $g^\pm$
	on the two sides of the slit, and impose
	\begin{equation}\label{eq:original-Dirichlet-data}
		\widetilde w=g\quad\text{on }\partial D_L.
	\end{equation}
	For any two adjacent boundary segments $J_1,J_2$ meeting at an
	ordinary corner $Q$, the endpoint values of $g$ satisfy
	\begin{equation}\label{eq:compat-left}
		(g|_{J_1})(Q)=(g|_{J_2})(Q).
	\end{equation}
	Here endpoint values are understood as one-sided limits, and the
	upper and lower inlet corners at $(0,0)$ are treated separately.
	At the trailing edge, the two slit traces have a common value:
	\begin{equation}\label{eq:original-boundary-structure}
		g^+(\ell,0)=g^-(\ell,0)=k.
	\end{equation}
	Near the outlet point $e_L$, assume
	\begin{equation}\label{eq:linear-outlet-vanishing}
		g(L,y_2)=0\quad\text{for }|y_2|\le\delta_{\mathrm{out}},
		\qquad \delta_{\mathrm{out}}=\tfrac14\min\{m_+,m_-\}.
	\end{equation}
	At $(L,0)$, we prescribe the pressure condition
	\begin{equation}\label{eq:finite-linear-pressure-normalization}
		\widetilde p(L,0)=p_0,
	\end{equation}
	where $p_0$ is a given constant.
	
	\Needspace{6\baselineskip}
	We consider the following linear boundary value problem.
	\begin{problem}[Linear problem in a finite nozzle]\label{prob:finite-linear}
		Suppose that $g$ satisfies
		\eqref{eq:compat-left}--\eqref{eq:linear-outlet-vanishing}.
		Find $(\widetilde w,\widetilde p)$ on $D_L$ such that:
		\begin{enumerate}
			\item The pair solves \eqref{eq:linear-first-order} in $D_L^+$ and $D_L^-$.
			\item The contact discontinuity conditions \eqref{eq:finite-linear-contact}
			hold on $\Gamma_L^*$.
			\item The Dirichlet condition \eqref{eq:original-Dirichlet-data}
			and condition
			\eqref{eq:finite-linear-pressure-normalization} hold.
		\end{enumerate}
	\end{problem}
	
	To study this problem, we assume that $A$ is a small perturbation of
	a diagonal matrix
	$A_b(y_2)=\operatorname{diag}(a_b^{11}(y_2),a_b^{22}(y_2))$,
	whose entries are piecewise $C^{1,\alpha}$ on $I$ and satisfy
	\begin{equation}\label{eq:background-bounds}
		\lambda_b\le a_b^{11},a_b^{22}\le\Lambda_b,\qquad
		\|a_b^{11}\|_{1,\alpha;I,\mathrm{pm}}
		+\|a_b^{22}\|_{1,\alpha;I,\mathrm{pm}}\le K_b,
	\end{equation}
	where $\lambda_b,\Lambda_b>0$ and $K_b\ge1$.
	More precisely, for $\tau>0$, we require
	\begin{equation}\label{eq:coeff-small}
		\|A-A_b\|_{1,\alpha;D_L,\mathrm{pm}}^{(-\alpha)}\le\tau,\qquad
		a^{12}(L,y_2)=a^{21}(L,y_2)=0
		\quad(0<|y_2|\le\delta_{\mathrm{out}}).
	\end{equation}
	Restricting the coefficients to any open set $G\subset D_L$
	preserves the perturbation bound:
	\begin{equation}\label{eq:coeff-restriction}
		\|A-A_b\|_{1,\alpha;G,\mathrm{pm}}^{(-\alpha)}
		\le\|A-A_b\|_{1,\alpha;D_L,\mathrm{pm}}^{(-\alpha)}\le\tau,
	\end{equation}
	with the same weight set $\mathcal E$.
	
	We measure the size of $g$ by the boundary norm
	\begin{equation}\label{eq:original-data-bound}
		\kappa_g:=\|g\|_{1,\alpha;(\beta);\partial D_L,\mathrm{pm}}^{(-\alpha)}
		:=|k|+\sum_{J\in\mathscr J}
		\|g|_J\|_{1,\alpha;(\beta);J}^{(-\alpha)},
	\end{equation}
	where $\mathscr J=\{I_{\mathrm{in}}^\sigma,\Gamma_{\mathrm w}^\sigma,
	\Gamma_{\mathrm a}^\sigma:\sigma\in\{+,-\}\}\cup\{O_L\}$.
	
	The following theorem gives the estimates independent of $L$ needed
	for the nonlinear construction in Section~\ref{sec:finite-nonlinear}.
	
	\begin{theorem}[Solvability of the linear problem]\label{thm:finite-linear}
		Suppose that $A_b$ and $A$ satisfy
		\eqref{eq:background-bounds} and \eqref{eq:coeff-small}, respectively,
		and that the boundary norm \eqref{eq:original-data-bound} is finite.
		Then there exist positive constants $\tau_0,C$ such that, if
		$\tau\le\tau_0$, Problem~\ref{prob:finite-linear} has a unique solution
		$(\widetilde w,\widetilde p)\in H^1(D_L)^2\cap C^0_{\mathrm{pm}}(\overline{D_L})^2$.
		Moreover, $\widetilde w\in C^{1,\alpha;(\beta)}_{(-\alpha),\mathrm{pm}}(D_L)$ and
		$\widetilde p-p_0\in C^{1,\alpha;(\beta-1,\beta)}_{(-\alpha),\mathrm{pm}}(D_L)$, with
		\[
		\|\widetilde w\|_{1,\alpha;(\beta);D_L,\mathrm{pm}}^{(-\alpha)}
		+\|\widetilde p-p_0\|_{1,\alpha;(\beta-1,\beta);D_L,\mathrm{pm}}^{(-\alpha)}
		\le C\kappa_g
		\]
		and
		\[
		\|\widetilde w\|_{H^1(D_L)}
		+\|\nabla\widetilde p\|_{L^2(D_L)}\le C\kappa_g.
		\]
		The constants $\tau_0$ and $C$ depend only on
		$\alpha,\beta,A_b,\ell,m_\pm$ and are independent of $L$.
	\end{theorem}
	
	Taking the divergence of \eqref{eq:linear-first-order} reduces the
	construction in Theorem~\ref{thm:finite-linear} to the scalar Dirichlet problem
	\begin{equation}\label{eq:slit-problem}
		\mathscr L\widetilde w:=\operatorname{div}(A\nabla \widetilde w)=0
		\quad\text{in }D_L,\qquad
		\widetilde w=g\quad\text{on }\partial D_L.
	\end{equation}
	The slit has two geometrically coincident boundary faces, so a direct
	application of boundary estimates on Lipschitz domains is unavailable.
	Excising a tubular neighbourhood of the slit produces the Lipschitz
	approximating domains $D_\varepsilon$, on which the variational
	Dirichlet problem can be solved. The essential step is to obtain
	estimates that remain uniform as the tube thickness tends to zero.
	In Subsection~\ref{subsec:truncated-linear}, local comparison regions
	are chosen to accommodate the slit geometry, and barrier functions
	compatible with the contact discontinuity conditions give downstream
	decay and pointwise control near the trailing edge and the other
	corners. Combined with local Schauder estimates and an explicit
	boundary lifting for the energy bound, these comparisons yield
	weighted H\"older and energy estimates independent of $\varepsilon$
	and $L$. For fixed $L$, they permit passage to the slit limit in
	Subsection~\ref{subsec:epsilon-limit}, where uniqueness follows from
	the weak formulation. The pressure is then recovered from
	\eqref{eq:linear-first-order} in Subsection~\ref{subsec:pressure};
	condition~\eqref{eq:finite-linear-pressure-normalization} fixes its
	additive constant, and the resulting pressure estimates complete
	the proof of Theorem~\ref{thm:finite-linear}.
	
	\subsection{Uniform estimates on the approximating domains}
	\label{subsec:truncated-linear}
	
\begin{figure}[htbp]
		\centering
		\begin{tikzpicture}[x=1.06cm,y=0.96cm,>=Latex,font=\small]
			\draw[thick] (0,2)--(11,2)--(11,-2)--(0,-2);
			\draw[thick] (0,2)--(0,0.43);
			\draw[thick] (0,-0.43)--(0,-2);
			\fill[gray!15]
			(0,0.43)--(4,0.43)
			arc[start angle=90,end angle=-90,radius=0.43]
			--(0,-0.43)--cycle;
			\draw[thick]
			(0,0.43)--(4,0.43)
			arc[start angle=90,end angle=-90,radius=0.43]
			--(0,-0.43);
			\draw[densely dashed] (4,0) circle (0.86);
			\draw[thick,dash dot] (4.43,0)--(11,0);
			\fill (4,0) circle (1.2pt);
			\fill (11,0) circle (1.2pt);
			\node at (1.55,0) {$K_\varepsilon$};
			\node at (1.80,1.03) {$\Gamma_{\mathrm a,\varepsilon}^+$};
			\draw[thin] (1.80,0.83)--(1.80,0.46);
			\node at (1.80,-1.03) {$\Gamma_{\mathrm a,\varepsilon}^-$};
			\draw[thin] (1.80,-0.83)--(1.80,-0.46);
			\node at (4.18,1.40) {$B_{2\varepsilon}(e_*)$};
			\draw[thin] (4.18,1.19)--(4.18,0.85);
			\node at (3.12,-0.91) {$e_*$};
			\draw[thin] (3.35,-0.76)--(3.94,-0.07);
			\node at (5.62,-0.78) {$\Sigma_\varepsilon:\ r=\varepsilon$};
			\draw[thin] (5.14,-0.64)--(4.41,-0.14);
			\node[above] at (7.55,0.13) {$\Gamma_\varepsilon^*$};
			\node at (8.0,1.20) {$D_\varepsilon^+$};
			\node at (8.0,-1.15) {$D_\varepsilon^-$};
			\node[right] at (11,0) {$e_L$};
			\node[below] at (11,-2.17) {$y_1=L$};
			\node[below] at (0,-2.17) {$y_1=0$};
		\end{tikzpicture}
		\caption{The approximating domain $D_\varepsilon$.}
		\label{fig:approx-domain}
	\end{figure}
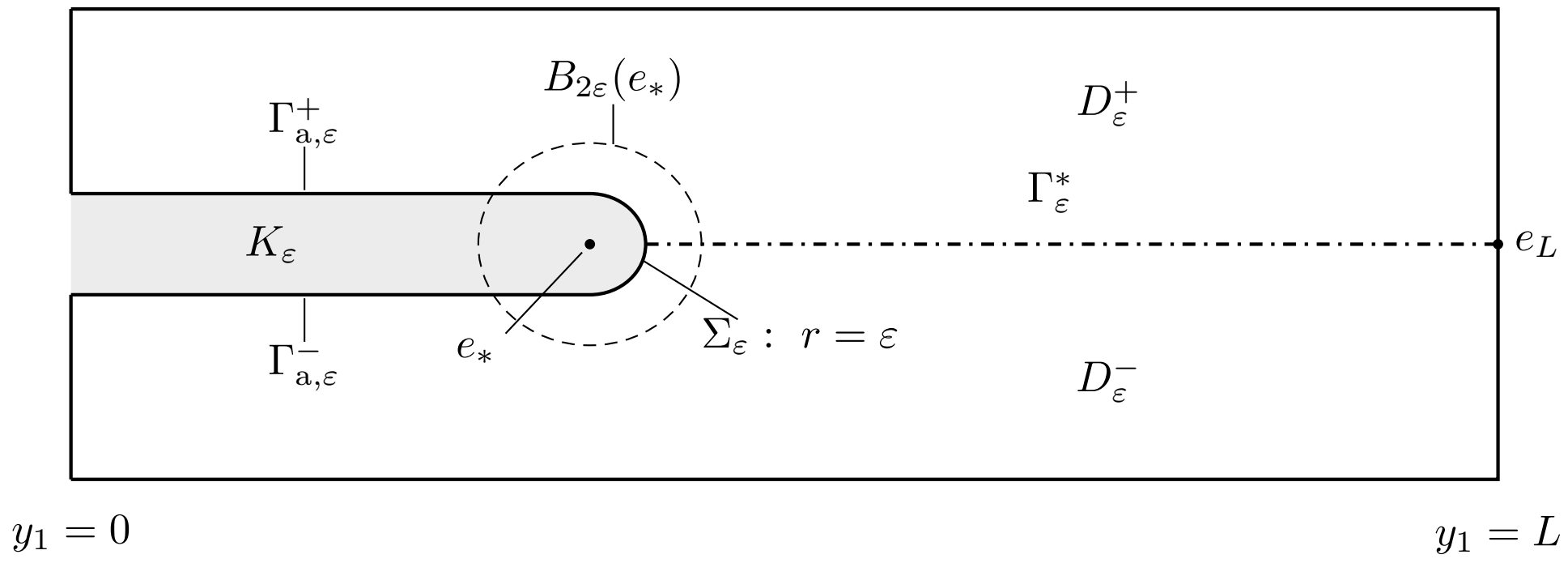

	For $0<\varepsilon<\min\{1,\ell/2,m_+/2,m_-/2\}$, we replace the slit
	with a strip of width $2\varepsilon$ capped by a semicircle on the right.
	More precisely, we remove
	\begin{equation}\label{eq:semicircle-obstacle}
		K_\varepsilon=([0,\ell]\times[-\varepsilon,\varepsilon])
		\cup\{(y_1,y_2):y_1\ge\ell,\ (y_1-\ell)^2+y_2^2\le\varepsilon^2\}
	\end{equation}
	and work on the approximating domain
	\begin{equation}\label{eq:approximate-domain}
		D_\varepsilon=\bigl((0,L)\times(-m_-,m_+)\bigr)\setminus K_\varepsilon.
	\end{equation}
	Intersecting $D_\varepsilon$ with the upper and lower half-planes gives
	\[
	D_\varepsilon^+=D_\varepsilon\cap\{y_2>0\},\qquad
	D_\varepsilon^-=D_\varepsilon\cap\{y_2<0\}.
	\]
	The shifted airfoil boundaries
	$\Gamma_{\mathrm a,\varepsilon}^+=(0,\ell)\times\{\varepsilon\}$ and
	$\Gamma_{\mathrm a,\varepsilon}^-=(0,\ell)\times\{-\varepsilon\}$
	meet the semicircular cap tangentially:
	\[
	\Sigma_\varepsilon
	=\{(\ell+\varepsilon\cos\theta,\varepsilon\sin\theta):
	-\pi/2\le\theta\le\pi/2\}.
	\]
	At the inlet, removing $K_\varepsilon$ leaves the segments
	\[
	I_\varepsilon^+=\{0\}\times(\varepsilon,m_+),\qquad
	I_\varepsilon^-=\{0\}\times(-m_-,-\varepsilon).
	\]
	The construction preserves the nozzle walls and the outlet
	$O:=O_L$. The closures of the inlet and wall segments,
	the outlet, the shifted airfoil boundaries and
	$\Sigma_\varepsilon$ together form $\partial D_\varepsilon$.
	The contact line now starts at the tip of the cap:
	\[
	\Gamma_\varepsilon^*=(\ell+\varepsilon,L)\times\{0\}.
	\]
	We establish the weighted derivative estimates on
	\begin{equation}\label{eq:estimated-domains}
		\widehat D_\varepsilon=D_\varepsilon\setminus\overline{B_{2\varepsilon}(e_*)},
		\qquad \widehat D_\varepsilon^\sigma=\widehat D_\varepsilon\cap D_\varepsilon^\sigma,
		\qquad \sigma\in\{+,-\}.
	\end{equation}
	On this set, the distance to the cap satisfies
	\begin{equation}\label{eq:semicircular-distance}
		r=|y-e_*|>2\varepsilon
		\quad\Longrightarrow\quad
		\operatorname{dist}(y,\Sigma_\varepsilon)
		\ge r-\varepsilon>\frac r2.
	\end{equation}
	Figure~\ref{fig:approx-domain} shows the geometry.
	
	\FloatBarrier
	
	To pose the Dirichlet problem on $D_\varepsilon$, we transfer $g$
	to its boundary. At the inlet, define the piecewise affine map
	\[
	h_\varepsilon(s)=
	\begin{cases}
		\displaystyle\frac{m_+}{m_+-\varepsilon}(s-\varepsilon),
		&\varepsilon\le s\le m_+,\\[4pt]
		\displaystyle\frac{m_-}{m_--\varepsilon}(s+\varepsilon),
		&-m_-\le s\le-\varepsilon.
	\end{cases}
	\]
	On each interval, $h_\varepsilon$ fixes the outer endpoint and sends
	the inner endpoint to zero. Its derivative and that of its inverse
	are uniformly bounded for small $\varepsilon$.
	
	On the shifted airfoil boundaries, we use a cutoff to join the
	traces $g^\pm$ on the slit to the constant value $k$ on $\Sigma_\varepsilon$. Choose
	$\chi\in C^\infty([0,\infty))$ with $0\le\chi\le1$,
	$\chi=0$ on $[0,1]$, and $\chi=1$ on $[2,\infty)$, and set
	\[
	\varphi_\varepsilon^\sigma(y_1)
	=k+\chi\!\left(\frac{\ell-y_1}{\varepsilon}\right)
	\bigl(g^\sigma(y_1,0)-k\bigr),\quad
	0\le y_1\le\ell,\quad \sigma\in\{+,-\}.
	\]
	These constructions give the boundary data
	\begin{equation}\label{eq:assembled-boundary-data}
		g_\varepsilon(y_1,y_2)=
		\begin{cases}
			g(0,h_\varepsilon(y_2)),
			&(y_1,y_2)\in I_\varepsilon^+\cup I_\varepsilon^-,\\
			g(y_1,y_2),       &(y_1,y_2)\in \Gamma_{\mathrm w}^+\cup \Gamma_{\mathrm w}^-\cup O,\\
			\varphi_\varepsilon^\sigma(y_1),
			&(y_1,y_2)\in \Gamma_{\mathrm a,\varepsilon}^\sigma,\\
			k,                 &(y_1,y_2)\in\Sigma_\varepsilon.
		\end{cases}
	\end{equation}
	The inlet map and the fact that the airfoil cutoff equals one near
	$y_1=0$ preserve the matching conditions \eqref{eq:compat-left} at
	the shifted inlet corners. Near $(\ell,\pm\varepsilon)$, the
	airfoil traces equal the cap value $k$. Thus $g_\varepsilon$ is
	continuous at every boundary junction.
	
	Using the inlet maps and rescaling the cutoff region
	$\ell-y_1\sim\varepsilon$, we deduce from
	\eqref{eq:original-data-bound} that
	\begin{equation}\label{eq:uniform-boundary-size}
		\begin{aligned}
			\|g_\varepsilon\|_{1,\alpha;(\beta);\partial D_\varepsilon,\mathrm{pm}}^{(-\alpha)}
			&:=|k|+\sum_{J\in\mathscr J_\varepsilon}
			\|g_\varepsilon|_J\|_{1,\alpha;(\beta);J}^{(-\alpha)}\le C\kappa_g,
		\end{aligned}
	\end{equation}
	where $C$ is independent of $L$ and $\varepsilon$, and
	\(
	\mathscr J_\varepsilon
	=\{I_\varepsilon^+,I_\varepsilon^-,\Gamma_{\mathrm w}^+,\Gamma_{\mathrm w}^-,
	\Gamma_{\mathrm a,\varepsilon}^+,\Gamma_{\mathrm a,\varepsilon}^-,O\}.
	\)
	
	With these boundary data, we consider the Dirichlet problem
	\begin{equation}\label{eq:bvp}
		\operatorname{div}(A\nabla \widetilde w_\varepsilon)=0\quad\text{in }D_\varepsilon,
		\qquad \widetilde w_\varepsilon=g_\varepsilon\quad\text{on }\partial D_\varepsilon,
	\end{equation}
	for $\widetilde w_\varepsilon\in H^1(D_\varepsilon)$, with $A$
	restricted from $D_L$.
	
	\Needspace{9\baselineskip}
	\begin{proposition}\label{prop:uniform}
		Suppose that \eqref{eq:background-bounds} and \eqref{eq:coeff-small}
		hold and that $g$ satisfies
		\eqref{eq:compat-left}--\eqref{eq:linear-outlet-vanishing},
		with $\kappa_g<\infty$ as defined in \eqref{eq:original-data-bound}.
		There exist positive constants $\tau_0,\varepsilon_0$, independent of
		$L$, such that, for $0<\varepsilon<\varepsilon_0$ and
		$\tau\le\tau_0$, the boundary value problem \eqref{eq:bvp} has a
		unique variational solution $\widetilde w_\varepsilon\in H^1(D_\varepsilon)$
		satisfying
		\begin{equation}\label{eq:uniform-statement}
			\|\widetilde w_\varepsilon\|_{H^1(D_\varepsilon)}
			+\|\widetilde w_\varepsilon\|_{1,\alpha;(\beta);\widehat D_\varepsilon,\mathrm{pm}}^{(-\alpha)}
			\le C\kappa_g,
		\end{equation}
		where $C$ is independent of $\varepsilon$ and $L$.
	\end{proposition}
	
	The proof uses three comparison estimates, which we establish first.
	In the next three lemmas, assume \eqref{eq:background-bounds} and
	\eqref{eq:coeff-small}, and let $\widetilde w=\widetilde w_\varepsilon
	\in H^1(D_\varepsilon)$ be a solution of \eqref{eq:bvp}, with
	$0<\varepsilon<\varepsilon_0$ and $0<\tau\le\tau_0$.
	The thresholds $\varepsilon_0,\tau_0>0$ will be chosen in the proofs,
	independently of $L$, and decreased when necessary.

	\begin{lemma}\label{lem:far-barrier}
		There exists a constant $C>0$, independent of $L$ and
		$\varepsilon$, such that
		\begin{equation}\label{eq:far-w-decay}
			|\widetilde w(y)|\le C\kappa_g(1+y_1)^{-\beta}
			\quad\text{in }D_\varepsilon.
		\end{equation}
	\end{lemma}
	\begin{proof}
		The weak maximum principle \cite[Theorem 8.1]{GT} and
		\eqref{eq:uniform-boundary-size} give
		\begin{equation}\label{eq:global-w-supremum}
			\|\widetilde w\|_{L^\infty(D_\varepsilon)}
			\le\|g_\varepsilon\|_{L^\infty(\partial D_\varepsilon)}
			\le C\kappa_g.
		\end{equation}
		To obtain the downstream decay, set $L_0=\ell+\tfrac12$ and
		choose $\varepsilon_0<\tfrac12$. Then
		$\mathcal R=(L_0,L)\times(-m_-,m_+)$ lies in $D_\varepsilon$, and
		\[
		d(y)\ge\tfrac12\qquad\text{for }y\in\mathcal R.
		\]
		After decreasing $\tau_0\le1$ if necessary, the coefficient assumptions
		yield constants $\lambda,\Lambda>0$, independent of $L,\varepsilon$, such that
		\begin{equation}\label{eq:barrier-ellipticity}
			\lambda|\xi|^2\le\xi^{\mathsf T}A(y)\xi\le\Lambda|\xi|^2
			\quad\text{in }D_\varepsilon^\pm,\qquad \xi\in\mathbb R^2.
		\end{equation}
		The lower bound for $d(y)$, $y\in\mathcal R$, also gives a uniform constant
		$M_1$ such that
		\[
		\left|\sum_{i=1}^2\partial_{y_i}a^{ij}\right|\le M_1
		\quad\text{in }\mathcal R\cap D_\varepsilon^\pm,\qquad j=1,2.
		\]
		
		With $\overline m=\max\{m_+,m_-\}$, consider the comparison function
		\begin{equation}\label{eq:far-exponential-barrier}
			V_\infty(y)
			=\frac{e^{C_1\overline m}-e^{C_1|y_2|}+1}
			{(C_2+y_1)^\beta},
			\qquad C_1>0,\quad C_2\ge1.
		\end{equation}
		Its numerator, $\phi(y_2)=e^{C_1\overline m}-e^{C_1|y_2|}+1$,
		satisfies
		\begin{equation}\label{eq:tmp1}
			1\le\phi(y_2)\le e^{C_1\overline m}.
		\end{equation}
		With $\sigma=y_2/|y_2|$ for $y_2\ne0$, differentiation gives
		\[
		(V_\infty)_{y_1}=-\frac{\beta}{C_2+y_1}V_\infty,\qquad
		(V_\infty)_{y_2}=-\frac{\sigma C_1e^{C_1|y_2|}}{\phi(y_2)}V_\infty.
		\]
		Differentiating once more gives
		\[
		(V_\infty)_{y_1y_1}=\frac{\beta(\beta+1)}{(C_2+y_1)^2}V_\infty,\quad
		(V_\infty)_{y_1y_2}=\frac{\sigma\beta C_1e^{C_1|y_2|}}{(C_2+y_1)\phi(y_2)}V_\infty,\quad
		(V_\infty)_{y_2y_2}=-\frac{C_1^2e^{C_1|y_2|}}{\phi(y_2)}V_\infty.
		\]
		Substituting these derivatives and using the symmetry of $A$, we find
		\[
		\begin{aligned}
			\mathscr L V_\infty
			={}&\biggl[-\frac{C_1e^{C_1|y_2|}}{\phi(y_2)}
			\left(C_1a^{22}+\sigma\partial_{y_i}a^{i2}
			-\frac{2\sigma\beta a^{12}}{C_2+y_1}\right)
			+\frac{\beta(\beta+1)a^{11}}{(C_2+y_1)^2}
			-\frac{\beta\partial_{y_i}a^{i1}}{C_2+y_1}\biggr]V_\infty.
		\end{aligned}
		\]
		Choose $C_1>0$ so that $M_1+2\beta\Lambda\le\lambda C_1/2$.
		For this choice, the term in parentheses satisfies, for every $C_2\ge1$,
		\begin{equation}\label{eq:far-axial-choice1}
			C_1a^{22}+\sigma\partial_{y_i}a^{i2}
			-\frac{2\sigma\beta a^{12}}{C_2+y_1}
			\ge\lambda C_1-M_1-2\beta\Lambda
			\ge\frac{\lambda C_1}{2}.
		\end{equation}
		With $C_1$ fixed, choose $C_2$ large enough that
		\begin{equation}\label{eq:far-axial-choice}
			\frac{\Lambda\beta(\beta+1)}{C_2^2}
			+\frac{\beta M_1}{C_2}
			\le\frac{\lambda C_1^2}{4e^{C_1\overline m}}.
		\end{equation}
		Combining~\eqref{eq:tmp1} and \eqref{eq:far-axial-choice1} with the choice of $C_2$
		in \eqref{eq:far-axial-choice} gives
		\begin{equation}\label{eq:far-operator-sign}
			\mathscr L V_\infty
			\le-\frac{\lambda C_1^2}{4e^{C_1\overline m}}V_\infty<0
			\qquad y\in\mathcal R,\quad y_2\ne0.
		\end{equation}
		
		To extend this inequality across $\Gamma_\varepsilon^*$ in the weak
		sense, we also need the sign of the conormal jump. On this line,
		\[
		[(A\nabla V_\infty)_2]
		=\left(-\frac{C_1(a^{22,+}+a^{22,-})}{e^{C_1\overline m}}
		-\frac{\beta[a^{21}]}{C_2+y_1}\right)V_\infty.
		\]
		The contribution of $[a^{21}]$ can be controlled by increasing $C_2$
		so that
		\[
		\left|\frac{\beta[a^{21}]}{C_2+y_1}\right|
		\le\frac{2\beta\Lambda}{C_2}
		\le\lambda C_1e^{-C_1\overline m}.
		\]
		Ellipticity also gives
		\[
		a^{22,+}+a^{22,-}\ge2\lambda.
		\]
		The two bounds imply
		\[
		[(A\nabla V_\infty)_2]
		\le-\lambda C_1e^{-C_1\overline m}V_\infty<0.
		\]
		This sign and \eqref{eq:far-operator-sign} show that $V_\infty$ is a
		weak supersolution on $\mathcal R$.
		
		We now turn to the estimates on $\partial\mathcal R$. On $y_1=L_0$, \eqref{eq:global-w-supremum} gives
		\[
		|\widetilde w|\le C\kappa_g.
		\]
		On the walls and at $y_1=L$, \eqref{eq:uniform-boundary-size}
		already gives the required decay:
		\[
		|g_\varepsilon|\le C\kappa_g(1+y_1)^{-\beta}.
		\]
		The bounds
		\[
		\phi\ge1,\qquad \frac{C_2+y_1}{1+y_1}\le C_2
		\]
		show that choosing
		\[
		C_3\ge C\max\{(C_2+L_0)^\beta,C_2^\beta\}
		\]
		makes $C_3\kappa_g V_\infty$ an upper bound for
		$|\widetilde w|$ on $\partial\mathcal R$.
		The weak comparison principle then yields
		\[
		|\widetilde w(y)|\le C_3\kappa_g V_\infty(y)
		\le C\kappa_g(1+y_1)^{-\beta}
		\quad\text{in }\mathcal R.
		\]
		Combining this estimate with \eqref{eq:global-w-supremum} on the region
		$y_1\le L_0$ proves \eqref{eq:far-w-decay}.
	\end{proof}
	
	\begin{lemma}\label{lem:tail-barrier}
		There exists $r_0>0$, independent of
		$L$ and $\varepsilon$, such that
		\begin{equation}\label{eq:barrier-tail-bound}
			|\widetilde w(y)-k|\le C\kappa_g|y-e_*|^\alpha
			\quad\text{in }D_\varepsilon\cap B_{r_0}(e_*),
		\end{equation}
		where $C$ is independent of $L$ and $\varepsilon$.
	\end{lemma}
	\begin{proof}
		For $\sigma\in\{+,-\}$, write
		$A_b^\sigma=A_b|_{I^\sigma}=(a_b^{ij,\sigma})$ and denote its
		one-sided limit at $y_2=0$ by $A_b^\sigma(0)$. To construct the
		barrier in $D_\varepsilon^\sigma$ near $e_*$, introduce
		\[
		\mu_\sigma
		=\sqrt{\frac{a_b^{11,\sigma}(0)}{a_b^{22,\sigma}(0)}},\qquad
		r_\sigma=\bigl((y_1-\ell)^2+\mu_\sigma^2y_2^2\bigr)^{1/2},\qquad
		\theta_\sigma=\arg\bigl(y_1-\ell+i\mu_\sigma y_2\bigr).
		\]
		Here $\arg$ denotes the principal argument, so that
		\[
		\theta_+\in(0,\pi)\quad\text{in }D_\varepsilon^+,\qquad
		\theta_-\in(-\pi,0)\quad\text{in }D_\varepsilon^-.
		\]
		The ellipticity bounds in \eqref{eq:background-bounds} give
		\[
		c_b r\le r_\sigma\le C_b r,\qquad r=|y-e_*|,
		\]
		where $c_b,C_b>0$ depend only on $\lambda_b,\Lambda_b$.
		Choose
		\[
		0<\alpha<\alpha_1<\tfrac12,
		\qquad \alpha_1\pi<\theta_0<\tfrac\pi2,
		\]
		and define
		\begin{equation}\label{eq:tail-barrier}
			V_*(y)
			=r_\sigma^\alpha
			\sin(\theta_0-\alpha_1|\theta_\sigma|),\qquad y\in D_\varepsilon^\sigma.
		\end{equation}
		The choice $\theta_0>\alpha_1\pi$ keeps the angular factor uniformly
		positive. Together with the bounds for $r_\sigma$, this gives
		\begin{equation}\label{eq:tail-barrier-bounds}
			c r^\alpha\le V_*\le C r^\alpha,
			\qquad |\nabla V_*|\le C r^{\alpha-1},
			\qquad |D^2V_*|\le C r^{\alpha-2}.
		\end{equation}
		The choice of $\mu_\sigma$ gives
		$a_b^{22,\sigma}(0)\mu_\sigma^2=a_b^{11,\sigma}(0)$, and hence
		\[
		\bigl(a_b^{11,\sigma}(0)\partial_{y_1}^2
		+a_b^{22,\sigma}(0)\partial_{y_2}^2\bigr)
		\bigl(r_\sigma^\alpha f(\theta_\sigma)\bigr)
		=a_b^{11,\sigma}(0)r_\sigma^{\alpha-2}(\alpha^2f+f'').
		\]
		For $f(\theta)=\sin(\theta_0-\alpha_1|\theta|)$, we have
		$f''=-\alpha_1^2f$ on both $(0,\pi)$ and $(-\pi,0)$.
		Substituting this identity gives
		\begin{align*}
			\operatorname{div}(A_b^\sigma(0)\nabla V_*)
			&=a_b^{11,\sigma}(0)(\alpha^2-\alpha_1^2)
			r_\sigma^{\alpha-2}
			\sin(\theta_0-\alpha_1|\theta_\sigma|)\le-c_1 r^{\alpha-2}\qquad\text{in }D_\varepsilon^\sigma.
		\end{align*}
		Here $c_1>0$ depends only on the background bounds and the chosen
		exponents and angles.
		
		To estimate the difference between the full operator and this
		constant-coefficient operator, write $\delta A=A-A_b=(\delta a^{ij})$,
		where
		\[
		\delta a^{ij}(y_1,y_2)
		:=a^{ij}(y_1,y_2)-a_b^{ij}(y_2),
		\qquad i,j\in\{1,2\},\quad (y_1,y_2)\in D_\varepsilon^\pm.
		\]
		Since $A_b$ is diagonal, $\delta a^{12}=a^{12}$ and
		$\delta a^{21}=a^{21}$. For $0<r=|y-e_*|<1$,
		\eqref{eq:coeff-small} and the background bounds give
		\[
		|\delta A|\le C\tau,\qquad
		|\nabla \delta A|\le C\tau r^{\alpha-1},\qquad
		|A_b(y_2)-A_b^\sigma(0)|\le C K_b r,\qquad
		|\nabla A_b|\le C K_b.
		\]
		The corresponding decomposition of the operator in $D_\varepsilon^\sigma$ is
		\begin{align*}
			\mathscr L V_*={}&A_b^\sigma(0):D^2V_*
			+(A_b(y_2)-A_b^\sigma(0)+\delta A):D^2V_* +\partial_{y_i}(a_b^{ij}+\delta a^{ij})
			\partial_{y_j}V_*.
		\end{align*}
		Here $M:N:=\sum_{i,j=1}^2 M_{ij}N_{ij}$ for $M,N\in\mathbb R^{2\times2}$.
		We infer from \eqref{eq:tail-barrier-bounds} that
		\[
		\big|\mathscr L V_*-
		\operatorname{div}(A_b^\sigma(0)\nabla V_*)\big|
		\le C(K_b r+\tau+\tau r^\alpha)r^{\alpha-2}.
		\]
		Choose $r_0\in(0,1)$, independently of $L$ and $\varepsilon$, such that
		\[
		B_{r_0}(e_*)\cap\partial D_\varepsilon
		\subset \Gamma_{\mathrm a,\varepsilon}^+
		\cup\Gamma_{\mathrm a,\varepsilon}^-\cup\Sigma_\varepsilon.
		\]
		By decreasing $r_0$ and then $\tau_0$ if necessary, we can also ensure that
		\[
		C\bigl(K_b r_0+\tau_0+\tau_0r_0^\alpha\bigr)\le c_1/2.
		\]
		Hence, for $y\in D_\varepsilon^\pm\cap B_{r_0}(e_*)$, we have
		\begin{equation}\label{eq:tail-local-super}
			\mathscr L V_*\le-\tfrac12 c_1r^{\alpha-2}.
		\end{equation}
		
		To check the supersolution inequality across $\Gamma_\varepsilon^*$,
		consider the traces at $y=(\ell+r,0)$, where $\varepsilon<r<r_0$:
		\[
		\begin{gathered}
			V_*^+=V_*^-=r^\alpha\sin\theta_0,\qquad
			\partial_{y_1}V_*^+=\partial_{y_1}V_*^-=
			\alpha\sin\theta_0\,r^{\alpha-1},\\
			\partial_{y_2}V_*^+=-\mu_+\alpha_1\cos\theta_0\,r^{\alpha-1},\qquad
			\partial_{y_2}V_*^-=\mu_-\alpha_1\cos\theta_0\,r^{\alpha-1}.
		\end{gathered}
		\]
		Using these traces and the diagonal form of $A_b^\pm(0)$, we obtain
		\[
		\begin{aligned}
			[(A_b\nabla V_*)_2]
			&=\bigl(A_b^+(0)\nabla V_*^+\bigr)_2
			-\bigl(A_b^-(0)\nabla V_*^-\bigr)_2\\
			&=-\alpha_1\cos\theta_0
			\bigl(a_b^{22,+}(0)\mu_++a_b^{22,-}(0)\mu_-\bigr)
			r^{\alpha-1}\\
			&=-c_\Gamma r^{\alpha-1},
		\end{aligned}
		\]
		where
		\[
		c_\Gamma:=\alpha_1\cos\theta_0
		\sum_{\sigma\in\{+,-\}}
		\sqrt{a_b^{11,\sigma}(0)a_b^{22,\sigma}(0)}>0.
		\]
		Replacing the background traces by those of $A$ contributes the
		correction
		\[
		\begin{split}
			\mathcal R_*(r)
			&:=[(A\nabla V_*)_2]-[(A_b\nabla V_*)_2]\\
			&=\Bigl(\alpha\sin\theta_0[\delta a^{21}]
			-\alpha_1\cos\theta_0
			(\mu_+\delta a^{22,+}+\mu_-\delta a^{22,-})\Bigr)
			r^{\alpha-1}.
		\end{split}
		\]
		Here the one-sided coefficient differences are
		\[
		\delta a^{ij,\sigma}(\ell+r,0)
		=a^{ij,\sigma}(\ell+r,0)-a_b^{ij,\sigma}(0),\qquad
		[\delta a^{21}]=\delta a^{21,+}-\delta a^{21,-},
		\]
		with all traces taken at $(\ell+r,0)$. The coefficient bound therefore yields
		\[
		|\mathcal R_*(r)|\le C\tau r^{\alpha-1}.
		\]
		Thus the jump remains negative after decreasing $\tau_0$ so that
		$C\tau_0\le c_\Gamma/2$:
		\begin{equation}\label{eq:tail-flux-negative}
			[(A\nabla V_*)_2]
			=-c_\Gamma r^{\alpha-1}+\mathcal R_*(r)
			\le-\frac{c_\Gamma}{2}r^{\alpha-1}<0.
		\end{equation}
		Since $V_*^+=V_*^-$ on $\Gamma_\varepsilon^*$,
		\eqref{eq:tail-local-super} and \eqref{eq:tail-flux-negative} imply
		\[
		\operatorname{div}(A\nabla V_*)\le0
		\]
		in the weak sense.
		
		It remains to compare $V_*$ with the boundary values of
		$|\widetilde w-k|$. On $\Gamma_{\mathrm a,\varepsilon}^+\cap B_{r_0}(e_*)$,
		the condition $g_\varepsilon(\ell,\varepsilon)=k$ and the boundary norm bound give
		\begin{equation}\label{eq:g-tail-modulus}
			\begin{aligned}
				|g_\varepsilon(y_1,\varepsilon)-k|
				&=\left|\int_{y_1}^{\ell}
				\partial_{y_1}g_\varepsilon(\xi,\varepsilon)\,d\xi\right|\\
				&\le C\kappa_g\int_{y_1}^{\ell}
				\bigl((\ell-\xi)^2+\varepsilon^2\bigr)^{(\alpha-1)/2}\,d\xi\\
				&\le C\kappa_g\int_{y_1}^{\ell}(\ell-\xi)^{\alpha-1}\,d\xi
				=\frac{C\kappa_g}{\alpha}(\ell-y_1)^\alpha
				\le C\kappa_g|y-e_*|^\alpha.
			\end{aligned}
		\end{equation}
		The same estimate holds on $\Gamma_{\mathrm a,\varepsilon}^-$, whereas
		$\widetilde w-k=0$ on $\Sigma_\varepsilon$. On the remaining circular
		boundary, \eqref{eq:uniform-boundary-size} and
		\eqref{eq:global-w-supremum} give
		\[
		|k|\le C\kappa_g,\qquad
		|\widetilde w-k|\le C\kappa_g
		\quad\text{on }\partial B_{r_0}(e_*)\cap D_\varepsilon.
		\]
		These boundary bounds, together with \eqref{eq:tail-barrier-bounds}
		and \eqref{eq:g-tail-modulus}, allow us to choose $C_0$ independently
		of $L,\varepsilon$ so that
		\[
		C_0\kappa_g V_*\ge|\widetilde w-k|
		\quad\text{on }\partial\bigl(D_\varepsilon\cap B_{r_0}(e_*)\bigr).
		\]
		The weak comparison principle applied to
		$C_0\kappa_g V_*-(\widetilde w-k)$ and
		$C_0\kappa_g V_*+(\widetilde w-k)$ now gives
		\[
		|\widetilde w(y)-k|\le C_0\kappa_g V_*(y)
		\le C\kappa_g|y-e_*|^\alpha,
		\]
		which proves \eqref{eq:barrier-tail-bound}.
	\end{proof}
	
	The inlet and outlet corners form the set
	\[
	\mathscr S_\varepsilon
	=\{(0,m_+),(0,-m_-),(0,\varepsilon),(0,-\varepsilon),(L,m_+),(L,-m_-)\}.
	\]
	We index these corners by $c$, writing $P_c\in\mathscr S_\varepsilon$.
	Let $\sigma(c)\in\{+,-\}$ be determined by
	$P_c\in\partial D_\varepsilon^{\sigma(c)}$, and let $J_{c,1},J_{c,2}$
	be the two exterior boundary segments meeting at $P_c$. Thus
	$P_c\in\overline{J_{c,1}}\cap\overline{J_{c,2}}$ and
	$J_{c,j}\subset\partial D_\varepsilon^{\sigma(c)}$.

	Fix $R_{\mathrm{loc}}>0$ sufficiently small,
	independently of $L$ and $\varepsilon$, so that
	$D_\varepsilon^{\sigma(c)}\cap B_{R_{\mathrm{loc}}}(P_c)$ is a quadrant
	sector for every $P_c\in\mathscr S_\varepsilon$ and stays at least a fixed positive
	distance from $\{e_*\}$.  For a solution $\widetilde w$ of \eqref{eq:bvp} and
	$P_c\in\mathscr S_\varepsilon$, set
	\[
	\mathcal M_c=
	\|\widetilde w\|_{L^\infty(D_\varepsilon^{\sigma(c)}\cap B_{R_{\mathrm{loc}}}(P_c))}
	+\sum_{j=1}^2
	\|g_\varepsilon|_{J_{c,j}}\|_{1,\alpha;J_{c,j}\cap B_{R_{\mathrm{loc}}}(P_c)}.
	\]
	Let $\tau_{c,j}$ be the unit tangent to $J_{c,j}$ pointing away from
	$P_c$, and let $g_\varepsilon(P_c)$ denote the common endpoint value
	of the two boundary traces. For $j=1,2$, define
	\[
	\partial_{\tau_{c,j}}g_\varepsilon(P_c)
	:=\lim_{s\to0^+}
	\frac{(g_\varepsilon|_{J_{c,j}})(P_c+s\tau_{c,j})
		-g_\varepsilon(P_c)}{s}.
	\]
	The tangents $\tau_{c,1},\tau_{c,2}$ are orthonormal. We match the
	common corner value and the two tangential derivatives with the affine
	function
	\begin{equation}\label{eq:inlet-jet}
		b_c:=\sum_{j=1}^2\partial_{\tau_{c,j}}g_\varepsilon(P_c)\,\tau_{c,j},
		\qquad \ell_c(y):=g_\varepsilon(P_c)+b_c\cdot(y-P_c).
	\end{equation}
	
	\begin{lemma}\label{lem:corner-barrier}
		There exists $r_0\in(0,R_{\mathrm{loc}})$,
		independent of $L$ and $\varepsilon$, such that
		\begin{equation}\label{eq:inlet-jet-decay}
			|\widetilde w(y)-\ell_c(y)|\le C\mathcal M_c|y-P_c|^{1+\alpha}
			\quad\text{in }D_\varepsilon^{\sigma(c)}\cap B_{r_0}(P_c),
		\end{equation}
		for every $P_c\in\mathscr S_\varepsilon$, where $C$ is independent of
		$L$ and $\varepsilon$.
	\end{lemma}
	\begin{proof}
		Fix $P_c\in\mathscr S_\varepsilon$, write $\sigma=\sigma(c)$, and put
		\[
		Q_c=D_\varepsilon^\sigma\cap B_{R_{\mathrm{loc}}}(P_c).
		\]
		Since $Q_c$ lies at a fixed positive distance from $e_*$, the coefficient
		assumption gives
		\[
		\|A\|_{1,\alpha;\overline{Q_c}}\le C K_b
		\qquad\text{for }\tau\le1,
		\]
		where $C$ depends only on the fixed geometry.
		
		By \eqref{eq:inlet-jet}, the coefficients of $\ell_c$ satisfy
		\[
		|g_\varepsilon(P_c)|+|b_c|\le C\mathcal M_c.
		\]
		Subtracting $\ell_c$, we set $v=\widetilde w-\ell_c$ and obtain on
		either boundary segment
		\[
		v(P_c+s\tau_{c,j})=\int_0^s
		\bigl(\partial_{\tau_{c,j}}g_\varepsilon(P_c+\xi\tau_{c,j})
		-\partial_{\tau_{c,j}}g_\varepsilon(P_c)\bigr)\,d\xi.
		\]
		The H\"older continuity of the tangential derivatives then gives, for
		$0<s<R_{\mathrm{loc}}$,
		\begin{equation}\label{eq:inlet-data-remainder}
			\begin{gathered}
				|v(P_c+s\tau_{c,j})|\le C\mathcal M_cs^{1+\alpha},\qquad
				|\partial_{\mathrm{tan}}v(P_c+s\tau_{c,j})|\le C\mathcal M_cs^\alpha,\\
				[\partial_{\mathrm{tan}}v]_{\alpha;J_{c,j}\cap B_{R_{\mathrm{loc}}}(P_c)}
				\le C\mathcal M_c,\qquad j=1,2.
			\end{gathered}
		\end{equation}
		Since $\ell_c$ has constant gradient $b_c$, the equation for $v$ is
		\[
		\mathscr L v=f_c:=-\partial_{y_i}(a^{ij}b_c^j)
		\quad\text{in }Q_c,
		\]
		where
		\[
		\|f_c\|_{L^\infty(Q_c)}\le CK_b\mathcal M_c.
		\]
		
		To control $v$ in $Q_c$, we construct a barrier adapted to the quadrant.
		For this fixed corner, introduce the weighted polar radius as in Lemma~\ref{lem:tail-barrier}:
		\[
		\mu_\sigma=\sqrt{\frac{a_b^{11,\sigma}(P_{c,2})}
			{a_b^{22,\sigma}(P_{c,2})}},\qquad
		r_\sigma=\bigl((y_1-P_{c,1})^2
		+\mu_\sigma^2(y_2-P_{c,2})^2\bigr)^{1/2}.
		\]
		Let $\theta_c\in\{\pm\pi/4,\pm3\pi/4\}$ be the angle bisecting
		the quadrant at $P_c$, and choose the branch
		\[
		\theta_\sigma=\arg\bigl(y_1-P_{c,1}
		+i\mu_\sigma(y_2-P_{c,2})\bigr)
		\in(\theta_c-\pi/4,\theta_c+\pi/4).
		\]
		As before, the background bounds give
		\[
		c_b|y-P_c|\le r_\sigma\le C_b|y-P_c|.
		\]
		Following the angular construction in \cite[Section~4]{CCF},
		choose $1+\alpha<\alpha_2<2$ and define
		\begin{equation}\label{eq:inlet-barrier}
			Z_c(y)=r_\sigma^{1+\alpha}
			\cos\bigl(\alpha_2(\theta_\sigma-\theta_c)\bigr).
		\end{equation}
		The restriction $\alpha_2<2$ ensures that the angular factor is
		positive throughout the quadrant, and gives
		\[
		\begin{gathered}
			\cos(\alpha_2\pi/4)r_\sigma^{1+\alpha}
			\le Z_c\le r_\sigma^{1+\alpha},\\
			|\nabla Z_c|\le Cr_\sigma^\alpha,\qquad
			|D^2Z_c|\le Cr_\sigma^{\alpha-1}.
		\end{gathered}
		\]
		Using $a_b^{22,\sigma}(P_{c,2})\mu_\sigma^2=a_b^{11,\sigma}(P_{c,2})$, we compute
		\[
		\begin{aligned}
			A_b^\sigma(P_{c,2}):D^2Z_c
			&=a_b^{11,\sigma}(P_{c,2})\bigl((1+\alpha)^2-\alpha_2^2\bigr)
			r_\sigma^{\alpha-1}
			\cos\bigl(\alpha_2(\theta_\sigma-\theta_c)\bigr)\\
			&\le-c_*r_\sigma^{\alpha-1}\qquad\text{in }Q_c,
		\end{aligned}
		\]
		where
		\[
		c_*:=\lambda_b\bigl(\alpha_2^2-(1+\alpha)^2\bigr)
		\cos(\alpha_2\pi/4)>0.
		\]
		The difference between this background operator and $\mathscr L$
		is controlled by the coefficient bound
		\[
		|A(y)-A_b^\sigma(P_{c,2})|\le C(\tau+K_b r_\sigma)
		\]
		which, together with the derivative bounds for $Z_c$, yields
		\[
		\begin{split}
			|\mathscr LZ_c-A_b^\sigma(P_{c,2}):D^2Z_c|
			&\le |A-A_b^\sigma(P_{c,2})|\,|D^2Z_c|
			+C|\nabla A|\,|\nabla Z_c|\\
			&\le C(\tau+K_b r_\sigma)r_\sigma^{\alpha-1}.
		\end{split}
		\]
		Choose $R_{\mathrm{loc}}$ at the outset so that
		$C_bR_{\mathrm{loc}}\le1$ and
		$C K_b C_bR_{\mathrm{loc}}\le c_*/4$.
		After decreasing $\tau_0$ to ensure $C\tau_0\le c_*/4$, the perturbation
		is absorbed by the negative background term, leaving
		\[
		\mathscr LZ_c\le-\frac{c_*}{2}r_\sigma^{\alpha-1}
		\quad\text{in }Q_c.
		\]
		
		For the boundary comparison, \eqref{eq:inlet-data-remainder} gives
		on the two straight segments of $\partial Q_c$
		\[
		|v|\le C\mathcal M_cr_\sigma^{1+\alpha}.
		\]
		On the circular part, the definition of $\mathcal M_c$ and
		the lower bound for $Z_c$ give
		\[
		|v|\le C\mathcal M_c,\qquad
		Z_c\ge cR_{\mathrm{loc}}^{1+\alpha}.
		\]
		These bounds allow us to choose a fixed constant
		$C_4\ge C(1+K_b+R_{\mathrm{loc}}^{-1-\alpha})$ large enough that
		$W_c=C_4\mathcal M_cZ_c$ satisfies
		\[
		W_c\ge|v|\quad\text{on }\partial Q_c,\qquad
		\mathscr LW_c\le-|f_c|\quad\text{in }Q_c.
		\]
		For the second inequality, the bound on $f_c$ is combined with
		\[
		r_\sigma^{\alpha-1}\ge1\qquad\text{for }0<r_\sigma\le1.
		\]
		The functions $W_c-v$ and $W_c+v$ therefore satisfy
		\[
		\begin{aligned}
			\mathscr L(W_c-v)&=\mathscr LW_c-f_c\le0,& W_c-v&\ge0\quad\text{on }\partial Q_c,\\
			\mathscr L(W_c+v)&=\mathscr LW_c+f_c\le0,& W_c+v&\ge0\quad\text{on }\partial Q_c.
		\end{aligned}
		\]
		The weak comparison principle gives
		\[
		|\widetilde w(y)-\ell_c(y)|\le W_c(y)
		\le C\mathcal M_c|y-P_c|^{1+\alpha}
		\quad\text{in }Q_c.
		\]
		Since the geometry and coefficient bounds are uniform in
		$c,L,\varepsilon$, we may choose a single
		$0<r_0<R_{\mathrm{loc}}$ for which \eqref{eq:inlet-jet-decay} holds at
		every corner.
	\end{proof}
	
	\Needspace{12\baselineskip}
	The preceding comparison estimates provide the pointwise control
	needed for the weighted H\"older estimate in Proposition~\ref{prop:uniform}.
	
	\begin{proof}[Proof of Proposition~\ref{prop:uniform}]
		\par\medskip\Needspace{3\baselineskip}\noindent\textbf{Step 1. Existence and energy bound.}
		To obtain an energy bound independent of $L$ and $\varepsilon$, we
		first extend $g_\varepsilon$ into $D_\varepsilon$. Fix
		$0<\rho<\tfrac18\min\{1,\ell,m_+,m_-\}$ and a smooth cutoff
		$\chi_*$ supported in $B_{2\rho}(e_*)$, equal to one on $B_\rho(e_*)$.
		Decrease $\varepsilon_0$ so that $\varepsilon_0<\rho/4$.
		The term $k\chi_*$ extends the cap value. For the remaining data near
		the trailing edge, set
		\[
		q_\varepsilon^\sigma(t):=\varphi_\varepsilon^\sigma(\ell-t)-k,
		\qquad 0<t<\rho,\quad \sigma\in\{+,-\}.
		\]
		By the definition of $\varphi_\varepsilon^\sigma$ and
		\eqref{eq:uniform-boundary-size},
		\[
		q_\varepsilon^\sigma(t)=0\quad(0<t\le\varepsilon),\qquad
		|q_\varepsilon^\sigma(t)|+t|(q_\varepsilon^\sigma)'(t)|
		\le C\kappa_g t^\alpha\quad(0<t<\rho).
		\]
		Choose fixed smooth cutoffs $\theta,\psi$ on $[0,\infty)$ with
		$\theta=1$ on $[0,\rho/4]$, $\theta=0$ on $[\rho/2,\infty)$,
		$\psi=1$ near zero and $\psi=0$ on $[1,\infty)$. Write
		\[
		t=\ell-y_1,\qquad n^+=y_2-\varepsilon,\qquad n^-=-y_2-\varepsilon,
		\]
		and define on $D_\varepsilon$
		\[
		E_\varepsilon^\sigma(y)=
		\begin{cases}
			\theta(t)q_\varepsilon^\sigma(t)\psi(n^\sigma/t),
			&0<t<\rho,\ n^\sigma>0,\\
			0, &\text{otherwise}.
		\end{cases}
		\]
		Differentiation and integration over $0<n^\sigma<t$ give
		\[
		\|E_\varepsilon^\sigma\|_{H^1(D_\varepsilon)}^2
		\le C\kappa_g^2\int_0^{\rho/2}t^{2\alpha-1}\,dt
		\le C\kappa_g^2,
		\]
		since $\alpha>0$.
		
		The residual boundary data
		\[
		b_\varepsilon:=g_\varepsilon-
		\left.(k\chi_*+E_\varepsilon^++E_\varepsilon^-)\right|_{\partial D_\varepsilon}
		\]
		vanish on $\partial D_\varepsilon\cap B_{\rho/4}(e_*)$.
		The matching endpoint values allow $b_\varepsilon$ to be extended
		through fixed-width inward strips and fixed-size corner patches
		with uniformly bounded overlap, treating the two airfoil boundaries
		and their inlet endpoints separately on each side. By
		\eqref{eq:uniform-boundary-size}, the resulting extension
		$F_\varepsilon$ satisfies
		\[
		\|F_\varepsilon\|_{H^1(D_\varepsilon)}^2
		\le C\kappa_g^2\left(1+\int_0^\infty(1+s)^{-2\beta}\,ds\right)
		\le C\kappa_g^2.
		\]
		Thus an extension of $g_\varepsilon$ with the required bound is
		\[
		\Psi_\varepsilon=k\chi_*+E_\varepsilon^++E_\varepsilon^-+F_\varepsilon,
		\qquad \|\Psi_\varepsilon\|_{H^1(D_\varepsilon)}\le C\kappa_g.
		\]
		
		Since $D_\varepsilon$ is bounded and Lipschitz and $A$ satisfies
		\eqref{eq:barrier-ellipticity}, the variational Dirichlet theorem
		\cite[Theorem 8.3]{GT} yields a unique
		$\widetilde w=\widetilde w_\varepsilon\in H^1(D_\varepsilon)$
		with trace $g_\varepsilon$, satisfying
		\begin{equation}\label{eq:weak}
			\int_{D_\varepsilon}A\nabla\widetilde w\cdot\nabla\varphi\,dy=0,
			\qquad \varphi\in H_0^1(D_\varepsilon).
		\end{equation}
		Testing with $\widetilde w-\Psi_\varepsilon$ and using ellipticity gives
		\[
		\|\nabla\widetilde w\|_{L^2(D_\varepsilon)}
		\le C\|\nabla\Psi_\varepsilon\|_{L^2(D_\varepsilon)}
		\le C\kappa_g.
		\]
		Together with the decay estimate in Lemma~\ref{lem:far-barrier}, this yields
		\[
		\|\widetilde w\|_{H^1(D_\varepsilon)}\le C\kappa_g.
		\]
		\par\medskip\Needspace{3\baselineskip}\noindent\textbf{Step 2. Estimates near the trailing edge.}
		Choose $0<r_0<\rho$, independently of $L$ and $\varepsilon$,
		no larger than the radius in Lemma~\ref{lem:tail-barrier}.
		Fix $\sigma\in\{+,-\}$ and a point $y\in \widehat D_\varepsilon^\sigma$ with
		\[
		2\varepsilon<r:=|y-e_*|<r_0/2,
		\]
		and set $s=r/32$ and $v=\widetilde w-k$.
		For $a>0$ and $z\in\overline{D_\varepsilon^\sigma}$, write
		$B_a^\sigma(z)=B_a(z)\cap D_\varepsilon^\sigma$.
		By \eqref{eq:semicircular-distance},
		\[
		\operatorname{dist}(y,\Sigma_\varepsilon)
		\ge r-\varepsilon>\frac r2=16s.
		\]
		Hence $B_{8s}(y)=B_{r/4}(y)$ avoids the cap and its junctions with
		the shifted airfoil boundaries and the contact line. The choice of
		$r_0$ also keeps this ball away from the inlet, outlet and nozzle walls.
		Consequently,
		\[
		\partial D_\varepsilon^\sigma\cap B_{8s}(y)
		\subset\Gamma_{\mathrm a,\varepsilon}^\sigma
		\cup\Gamma_\varepsilon^*.
		\]
		Moreover, $B_{8s}(y)\subset B_{r_0}(e_*)$ and
		\[
		\frac{3r}{4}\le d(y')=\min\{1,|y'-e_*|\}=|y'-e_*|\le\frac{5r}{4}
		\qquad\bigl(y'\in B_{8s}(y)\bigr).
		\]
		Lemma~\ref{lem:tail-barrier} therefore gives
		\[
		\|v\|_{L^\infty(B_{8s}(y)\cap D_\varepsilon)}
		\le C\kappa_g r^\alpha.
		\]
		We next prove that
		\begin{equation}\label{eq:tail-pointwise}
			|\nabla \widetilde w(y)|\le C\kappa_g r^{\alpha-1}, \qquad
			[\nabla \widetilde w]_{\alpha;B_{r/64}^\sigma(y)}
			\le C\kappa_g r^{-1},
		\end{equation}
		with $C$ independent of $y,L$ and $\varepsilon$.
		To obtain these bounds, we estimate
		$\|v\|_{1,\alpha;B_{s/2}^\sigma(y)}'$ separately in the interior,
		near a shifted airfoil boundary and near the contact line.
		
		\smallskip\noindent\emph{(i) Interior points.}
		If $\operatorname{dist}(y,\partial D_\varepsilon^\sigma)\ge s$,
		then $B_s(y)\subset D_\varepsilon^\sigma$, so
		$B_{s/2}^\sigma(y)=B_{s/2}(y)$. Since $s=r/32$,
		\eqref{eq:coeff-restriction} gives
		\[
		\|A\|_{0,\alpha;B_s(y)}
		\le C\bigl(K_b+\tau r^{\alpha-1}s^{1-\alpha}\bigr)\le C.
		\]
		The interior Schauder estimate in its scaled form
		\cite[Corollary 8.36]{GT} therefore gives
		\[
		\|v\|_{1,\alpha;B_{s/2}(y)}'
		\le C\|v\|_{L^\infty(B_s(y))}
		\le C\kappa_g r^\alpha.
		\]
		
		\smallskip\noindent\emph{(ii) Points near a shifted airfoil boundary.}
		Suppose $\operatorname{dist}(y,\Gamma_{\mathrm a,\varepsilon}^\sigma)<s$.
		Without loss of generality, take $y\in D_\varepsilon^+$.
		Its orthogonal projection onto the shifted airfoil boundary is
		$z=(y_1,\varepsilon)\in\Gamma_{\mathrm a,\varepsilon}^+$.
		The half-ball $B_{4s}^+(z)$ lies in
		$B_{8s}(y)$ and has flat boundary
		$J=\Gamma_{\mathrm a,\varepsilon}^+\cap B_{4s}(z)$.
		The boundary bounds \eqref{eq:uniform-boundary-size} and
		\eqref{eq:g-tail-modulus}, with $s=r/32$, give
		\[
		\|g_\varepsilon-k\|_{1,\alpha;J}'\le C\kappa_g r^\alpha.
		\]
		The coefficient bound in case (i) also holds on $B_{4s}^+(z)$.
		Applying the Dirichlet boundary estimate in its scaled form
		\cite[Corollary 8.36 and the paragraph following it]{GT}
		and restricting to $B_{s/2}^+(y)\subset B_{2s}^+(z)$, we obtain
		\[
		\|v\|_{1,\alpha;B_{s/2}^+(y)}'
		\le C\bigl(\|v\|_{L^\infty(B_{4s}^+(z))}
		+\|g_\varepsilon-k\|_{1,\alpha;J}'\bigr).
		\]
		For $y\in D_\varepsilon^-$, the same argument applies with
		$z=(y_1,-\varepsilon)$. Thus, on either side,
		\[
		\|v\|_{1,\alpha;B_{s/2}^\sigma(y)}'\le C\kappa_g r^\alpha.
		\]
		
		\smallskip\noindent\emph{(iii) Points near the contact line.}
		Suppose $\operatorname{dist}(y,\Gamma_\varepsilon^*)<s$.
		Let $z=(y_1,0)\in\Gamma_\varepsilon^*$ be the orthogonal projection
		of $y$. Then $B_{4s}(z)\subset D_\varepsilon\cap B_{8s}(y)$.
		Under the change of variables $y'=z+4s\xi$, set
		\[
		\widehat v(\xi)=v(z+4s\xi),\qquad \widehat A(\xi)=A(z+4s\xi).
		\]
		By \eqref{eq:weak}, $\widehat v$ satisfies the weak equation
		\[
		\operatorname{div}_\xi(\widehat A\nabla_\xi\widehat v)=0
		\quad\text{in }B_1(0).
		\]
		The contact line is mapped to $\xi_2=0$.
		Since $d(y')\ge3r/4$ on $B_{4s}(z)$ and $s=r/32$,
		\eqref{eq:coeff-restriction} gives
		\begin{equation}\label{eq:scaled-coefficients}
			\|\widehat A\|_{0,\alpha;B_1(0),\mathrm{pm}}\le C,
		\end{equation}
		where the piecewise norm is taken over the two half-balls and
		$C$ is independent of $y,L$ and $\varepsilon$.
		Applying the piecewise Schauder estimate of Chen--Deng
		\cite[Lemma~2.1]{CD} to $\widehat v$ and returning to the original
		variables, we obtain on $B_{s/2}^\sigma(y)\subset B_{2s}^\sigma(z)$
		\[
		\|v\|_{1,\alpha;B_{s/2}^\sigma(y)}'
		\le C\|v\|_{L^\infty(B_{4s}(z))}
		\le C\kappa_g r^\alpha.
		\]
		
		The boundary inclusion above shows that these three cases cover
		all possibilities. Since $s=r/32$ and $\nabla v=\nabla\widetilde w$,
		they prove \eqref{eq:tail-pointwise}. For points $y,y'$ on the
		same side, the local H\"older bound applies when
		$|y-y'|<d(y,y')/64$; the pointwise gradient bound controls the
		remaining pairs. Together with \eqref{eq:global-w-supremum} and
		the boundedness of $\Delta_y$ near $e_*$, this gives
		\[
		\|\widetilde w\|_{1,\alpha;(\beta);\widehat D_\varepsilon\cap B_{r_0/2}(e_*),\mathrm{pm}}^{(-\alpha)}
		\le C\kappa_g.
		\]
		
		\par\medskip\Needspace{3\baselineskip}\noindent\textbf{Step 3. Estimates at the Dirichlet corners.}
		If necessary, decrease $r_0$ to the common radius in
		Lemma~\ref{lem:corner-barrier}, which is also independent of $L$
		and $\varepsilon$.
		Fix $P_c\in\mathscr S_\varepsilon$, put $\sigma=\sigma(c)$, and set
		$v=\widetilde w-\ell_c$. The neighborhood $B_{r_0}^\sigma(P_c)$
		stays a fixed positive distance from $e_*$, so the coefficients
		have uniformly bounded $C^{1,\alpha}$ norms there. Subtracting the
		affine function $\ell_c$ from the equation for $\widetilde w$ gives
		\[
		\mathscr Lv=f_c:=-\partial_{y_i}(a^{ij}b_c^j)
		\quad\text{in }B_{r_0}^\sigma(P_c),
		\]
		where
		\begin{equation}\label{eq:corner-source-bound}
			\|f_c\|_{L^\infty(B_{r_0}^\sigma(P_c))}\le CK_b\mathcal M_c.
		\end{equation}
		Lemma~\ref{lem:corner-barrier} gives
		\[
		|v(y')|\le C\mathcal M_c|y'-P_c|^{1+\alpha}
		\qquad\bigl(y'\in B_{r_0}^\sigma(P_c)\bigr).
		\]
		
		For $y\in B_{r_0/2}^\sigma(P_c)$, put $r=|y-P_c|$ and $s=r/32$.
		The ball $B_{8s}(y)$ avoids $P_c$ and meets at most one of the two
		boundary segments. If
		$\operatorname{dist}(y,\partial D_\varepsilon^\sigma)\ge s$,
		we set
		\[
		\widehat v(\xi)=v(y+s\xi),\qquad \widehat A(\xi)=A(y+s\xi).
		\]
		The rescaled equation is
		\[
		\operatorname{div}_\xi(\widehat A\nabla_\xi\widehat v)
		=s^2f_c(y+s\xi)\quad\text{in }B_1(0).
		\]
		The interior estimate \cite[Corollary 8.36]{GT} gives, in the original variables,
		\[
		\|v\|_{1,\alpha;B_{s/2}(y)}'
		\le C\bigl(\|v\|_{L^\infty(B_s(y))}
		+s^2\|f_c\|_{L^\infty(B_s(y))}\bigr).
		\]
		Otherwise, let $z$ be the orthogonal projection of $y$ onto the
		nearby boundary segment. The half-ball $B_{4s}^\sigma(z)$ lies in
		$B_{8s}(y)$ and has flat boundary
		$J=\partial D_\varepsilon^\sigma\cap B_{4s}(z)$.
		By \eqref{eq:inlet-data-remainder},
		\[
		\|g_\varepsilon-\ell_c\|_{1,\alpha;J}'
		\le C\mathcal M_cr^{1+\alpha}.
		\]
		On the unit half-ball $\widehat G=(B_{4s}^\sigma(z)-z)/(4s)$, set
		\[
		\widehat v(\xi)=v(z+4s\xi),\qquad \widehat A(\xi)=A(z+4s\xi).
		\]
		Then
		\[
		\operatorname{div}_\xi(\widehat A\nabla_\xi\widehat v)
		=16s^2f_c(z+4s\xi)\quad\text{in }\widehat G.
		\]
		The Dirichlet boundary estimate \cite[Corollary 8.36 and the paragraph following
		it]{GT} gives, in the original variables on $B_{s/2}^\sigma(y)\subset
		B_{2s}^\sigma(z)$,
		\[
		\|v\|_{1,\alpha;B_{s/2}^\sigma(y)}'
		\le C\bigl(\|v\|_{L^\infty(B_{4s}^\sigma(z))}
		+\|g_\varepsilon-\ell_c\|_{1,\alpha;J}'+s^2\|f_c\|_{L^\infty(B_{4s}^\sigma(z))}\bigr).
		\]
		
		Using \eqref{eq:inlet-jet-decay}, \eqref{eq:inlet-data-remainder}
		and \eqref{eq:corner-source-bound}, together with $r^2\le r^{1+\alpha}$,
		we obtain in both cases
		\begin{equation}\label{eq:corner-local-scaled}
			\|v\|_{1,\alpha;B_{r/64}^\sigma(y)}'
			\le C\mathcal M_c(r^{1+\alpha}+K_br^2)
			\le C\mathcal M_cr^{1+\alpha}.
		\end{equation}
		The gradient decay and local H\"older bounds in
		\eqref{eq:corner-local-scaled} extend $\nabla v$ continuously to
		$P_c$, with value zero. Adding back the affine function $\ell_c$
		gives
		\[
		\|\widetilde w\|_{1,\alpha;\overline{B_{r_0/4}^\sigma(P_c)}}\le C\mathcal M_c.
		\]
		Since $\Delta_y$ is comparable to $1+P_{c,1}$ on this neighborhood,
		Lemma~\ref{lem:far-barrier} and \eqref{eq:uniform-boundary-size} yield
		\begin{equation}\label{eq:corner-unweighted-bound}
			\|\widetilde w\|_{1,\alpha;(\beta);B_{r_0/4}^\sigma(P_c)}^{(-\alpha)}
			\le C(1+P_{c,1})^\beta\mathcal M_c\le C\kappa_g.
		\end{equation}

		\par\medskip\Needspace{3\baselineskip}\noindent\textbf{Step 4. The endpoint at the outlet.}
		Choose $0<r_{\mathrm{out}}<\min\{r_0,\delta_{\mathrm{out}}/4\}$, independently
		of $L$ and $\varepsilon$. Since the outlet data vanish on
		$B_{2r_{\mathrm{out}}}(e_L)\cap O$, we extend the solution by odd reflection
		across $y_1=L$. With
		\[
		\mathsf R=\operatorname{diag}(-1,1),\qquad
		\mathscr R(y)=e_L+\mathsf R(y-e_L)=(2L-y_1,y_2),
		\]
		we extend $\widetilde w$ and $A$ to $B_{2r_{\mathrm{out}}}(e_L)\cap\{y_1>L\}$,
		retaining the same notation, by
		\[
		\widetilde w(y)=-\widetilde w(\mathscr R(y)),\qquad A(y)=\mathsf R A(\mathscr R(y))\mathsf R.
		\]
		Then $\widetilde w\in H^1(B_{2r_{\mathrm{out}}}(e_L))$ satisfies
		\[
		\operatorname{div}(A\nabla \widetilde w)=0
		\quad\text{in }B_{2r_{\mathrm{out}}}(e_L)
		\]
		in the weak sense. The vanishing of the off-diagonal coefficients on the outlet in
		\eqref{eq:coeff-small} preserves the piecewise $C^{0,\alpha}$ coefficient bounds
		under this reflection. Since $\Delta_y$ is comparable to $1+L$ on
		$B_{r_{\mathrm{out}}}(e_L)\cap D_\varepsilon$, the estimate of Chen--Deng
		\cite[Lemma~2.1]{CD} and Lemma~\ref{lem:far-barrier} give
		\[
		\|\widetilde w\|_{1,\alpha;(\beta);B_{r_{\mathrm{out}}}(e_L)\cap D_\varepsilon,\mathrm{pm}}^{(-\alpha)}
		\le C(1+L)^\beta\|\widetilde w\|_{L^\infty(B_{2r_{\mathrm{out}}}(e_L)\cap D_\varepsilon)}
		\le C\kappa_g.
		\]
		The constant $C$ is independent of $L$ and $\varepsilon$.
		
		\par\medskip\Needspace{3\baselineskip}\noindent\textbf{Step 5. The global weighted estimate.}
		Fix $\sigma\in\{+,-\}$ and let $G_\varepsilon^\sigma$ be the set
		of points $y\in \widehat D_\varepsilon^\sigma$ satisfying
		\[
		|y-e_*|\ge r_0/4,\qquad |y-e_L|\ge r_{\mathrm{out}}/2,\qquad
		|y-P_c|\ge r_0/8\quad\text{whenever }\sigma(c)=\sigma.
		\]
		On neighborhoods of a fixed radius centered in $G_\varepsilon^\sigma$,
		the coefficients have uniform $C^{0,\alpha}$ bounds on either side
		and the downstream weight is comparable to its value at the center.
		We apply the interior and Dirichlet boundary estimates
		\cite[Corollary 8.36 and the paragraph following it]{GT}, using
		the piecewise estimate of Chen--Deng \cite[Lemma~2.1]{CD} where
		a neighborhood crosses the contact line. With the solution and
		boundary bounds from Lemma~\ref{lem:far-barrier} and
		\eqref{eq:uniform-boundary-size}, these estimates yield
		\[
		\|\widetilde w\|_{1,\alpha;(\beta);G_\varepsilon^\sigma}^{(-\alpha)}
		\le C\kappa_g,
		\]
		with $C$ independent of $L$ and $\varepsilon$.
		
		We may decrease $\varepsilon_0$ so that $2\varepsilon_0<r_0/4$.
		Together with $G_\varepsilon^\sigma$, the neighborhoods in
		Steps~2--4 cover $\widehat D_\varepsilon^\sigma$; they overlap
		$G_\varepsilon^\sigma$ along the annular portions
		\[
		r_0/4<|y-e_*|<r_0/2,\qquad
		r_0/8<|y-P_c|<r_0/4,\qquad
		r_{\mathrm{out}}/2<|y-e_L|<r_{\mathrm{out}}.
		\]
		The radial widths of these overlaps are at least
		$\min\{r_0/8,r_{\mathrm{out}}/2\}>0$, independently of $L$ and $\varepsilon$.
		Close pairs are controlled by the local estimates, using Step~2 near
		$e_*$; for the remaining pairs, the weighted pointwise gradient bounds
		control the difference quotient.
		Combining the estimate on $G_\varepsilon^\sigma$ with the local
		weighted estimates in Steps~2--4 therefore gives
		\begin{equation}\label{eq:weighted-proof-conclusion}
			\|\widetilde w\|_{1,\alpha;(\beta);\widehat D_\varepsilon,\mathrm{pm}}^{(-\alpha)}\le C\kappa_g.
		\end{equation}
		Together with the $H^1$ bound in Step~1, this proves
		\eqref{eq:uniform-statement}.
	\end{proof}
	
	\subsection{Passage to the limit as \texorpdfstring{$\varepsilon\to0$}{epsilon tends to zero}}
	\label{subsec:epsilon-limit}
	
	Keeping $L$, $A$ and $g$ fixed, we pass to the limit as $\varepsilon\to0$.
	The domains $D_\varepsilon$ exhaust $D_L$, while the boundary data
	constructed above converge to $g$ away from $e_*$.
	
	\begin{proposition}\label{prop:slit-limit}
		Suppose that \eqref{eq:background-bounds} and
		\eqref{eq:coeff-small} hold with $\tau\le\tau_0$, and that $g$ satisfies
		\eqref{eq:compat-left}, \eqref{eq:original-boundary-structure} and
		\eqref{eq:linear-outlet-vanishing}, with finite norm
		\eqref{eq:original-data-bound}. The boundary value problem
		\eqref{eq:slit-problem} has a unique
		solution $\widetilde w\in H^1(D_L)\cap
		C^{1,\alpha;(\beta)}_{(-\alpha),\mathrm{pm}}(D_L)$ satisfying
		$\widetilde w(e_*)=k$ and
		\begin{equation}\label{eq:finite-W-bound}
			\|\widetilde w\|_{H^1(D_L)}
			+\|\widetilde w\|_{1,\alpha;(\beta);D_L,\mathrm{pm}}^{(-\alpha)}
			\le C\kappa_g,
		\end{equation}
		where $C$ is independent of $L$.
	\end{proposition}
	\begin{proof}
		\noindent\textbf{Step 1. Uniform estimates.}
		Applying Proposition~\ref{prop:uniform} to the restrictions of $A$
		and the constructed boundary data gives
		\[
		\|\widetilde w_\varepsilon\|_{H^1(D_\varepsilon)}
		+\|\widetilde w_\varepsilon\|_{1,\alpha;(\beta);\widehat D_\varepsilon,\mathrm{pm}}^{(-\alpha)}
		\le C\kappa_g.
		\]
		
		\par\medskip\Needspace{3\baselineskip}\noindent\textbf{Step 2. Passage to the limit.}
		To treat the moving slit boundaries on fixed sets, choose
		$0<a<\ell$ and $0<b<\frac12\min\{m_+,m_-\}$, and set
		\[
		Q^+=(0,a)\times(0,b),\qquad Q^-=(0,a)\times(-b,0).
		\]
		On these rectangles, define the pullbacks
		\[
		v_\varepsilon^\sigma=\widetilde w_\varepsilon\circ T_\varepsilon^\sigma,
		\qquad T_\varepsilon^\pm(y)=(y_1,y_2\pm\varepsilon).
		\]
		For sufficiently small $\varepsilon$, the translated rectangles lie
		in $\widehat D_\varepsilon^\sigma$, and Step~1 gives uniform $C^{1,\alpha}$
		bounds for $v_\varepsilon^\sigma$ on $\overline{Q^\sigma}$.
		The pullbacks $v_\varepsilon^\sigma$ are used only on $Q^\sigma$.
		On fixed neighborhoods away from the closed slit $[0,\ell]\times\{0\}$,
		the original functions $\widetilde w_\varepsilon$ have the same one-sided
		bounds for all sufficiently small $\varepsilon$.
		A countable family of the sets $\overline{Q^\sigma}$ with $a\to\ell$,
		together with such fixed neighborhoods, covers
		$\overline{D_L^\sigma}\setminus\{e_*\}$. Applying Arzel\`a--Ascoli
		diagonally to the pullbacks near the slit and to
		$\widetilde w_\varepsilon$ on the remaining neighborhoods gives local
		$C^{1,\alpha'}$ limits on each side for every $0<\alpha'<\alpha$.
		
		These limits agree on overlaps. Indeed, any compact subset $K$
		of the overlap between $Q^\sigma$ and an unshifted neighborhood
		stays a positive distance from the closed slit. The uniform gradient and
		H\"older bounds therefore give
		\[
		\|v_\varepsilon^\sigma-\widetilde w_\varepsilon\|_{C^1(K)}
		\le C\kappa_g(\varepsilon+\varepsilon^\alpha)\longrightarrow0.
		\]
		Here $C$ may depend on $K$, but is independent of $\varepsilon$.
		Thus there are $\varepsilon_j\to0$ and a single limit $\widetilde w$
		such that, for each $\sigma\in\{+,-\}$,
		\begin{equation}\label{eq:limit-local-convergence}
			\widetilde w_{\varepsilon_j}\longrightarrow\widetilde w
			\quad\text{in }C^{1,\alpha'}(K),\qquad 0<\alpha'<\alpha,
		\end{equation}
		for every compact $K\subset\overline{D_L^\sigma}$ disjoint from the
		closed slit, and
		\[
		v_{\varepsilon_j}^\sigma\longrightarrow\widetilde w
		\quad\text{in }C^{1,\alpha'}(\overline{Q^\sigma}).
		\]
		The pulled-back coefficients $A\circ T_\varepsilon^\sigma$ converge
		locally uniformly to $A$. By \eqref{eq:assembled-boundary-data},
		the corresponding pullbacks of the boundary data converge to $g$
		in local $C^{1,\alpha'}$ norms on each boundary segment away from
		$e_*$. Hence the limit solves the equation in $D_L^+$ and $D_L^-$
		and has the prescribed boundary values.
		Continuity of the approximating solutions across
		$\Gamma_{\varepsilon_j}^*$, together with
		\eqref{eq:limit-local-convergence} on compact subsets of $\Gamma_L^*$,
		also gives
		\[
		\widetilde w^+=\widetilde w^-\quad\text{on }\Gamma_L^*.
		\]
		The weighted pointwise and two-point inequalities pass to the limit
		with exponent $\alpha$; in particular, the comparison estimate at
		the trailing edge yields
		\[
		|\widetilde w(y)-k|\le C\kappa_g|y-e_*|^\alpha.
		\]
		
		Every open set $G\Subset D_L$ lies in $D_{\varepsilon_j}$ for all
		sufficiently large $j$. The uniform energy bound therefore yields,
		after extraction, weak $H^1(G)$ convergence to the same locally
		uniform limit. By lower semicontinuity and exhaustion of $D_L$,
		\[
		\|\widetilde w\|_{H^1(D_L)}^2
		\le\liminf_{j\to\infty}
		\|\widetilde w_{\varepsilon_j}\|_{H^1(D_{\varepsilon_j})}^2
		\le C\kappa_g^2.
		\]
		
		\par\medskip\Needspace{3\baselineskip}\noindent\textbf{Step 3. The weak equation and uniqueness.}
		Since every $\varphi\in C_c^\infty(D_L)$ is an admissible test function
		on $D_\varepsilon$ for all sufficiently small $\varepsilon$, we may
		pass to the limit in \eqref{eq:weak}. Density then gives
		\[
		\int_{D_L}A\nabla \widetilde w\cdot\nabla \varphi\,dy=0,
		\qquad\varphi\in H_0^1(D_L).
		\]
		For uniqueness, let $\widetilde w_1,\widetilde w_2$ be two solutions
		in the stated class and set $z=\widetilde w_1-\widetilde w_2$.
		Then $z\in H^1(D_L)$ has zero trace on every exterior boundary
		segment, including both sides of the slit. To show that
		$z\in H_0^1(D_L)$, we first remove a neighborhood of the trailing edge using a radial cutoff $\chi_\delta$ which
		vanishes on $B_\delta(e_*)$ and equals one outside $B_{2\delta}(e_*)$,
		with $|\nabla\chi_\delta|\le C/\delta$.
		The weighted derivative bounds and the common value at $e_*$ give
		\[
		|z(y)|\le C|y-e_*|^\alpha.
		\]
		Consequently,
		\[
		\|z\nabla\chi_\delta\|_{L^2(D_L)}^2
		\le C\delta^{-2}\int_{B_{2\delta}(e_*)}|y-e_*|^{2\alpha}\,dy
		\le C\delta^{2\alpha}\longrightarrow0.
		\]
		Away from $e_*$, we use a partition of unity to localize near the
		outer boundary and separately in half-rectangles along the two slit
		faces and quadrants at their inlet endpoints. Each localized function
		has zero trace on its Lipschitz patch and vanishes near the artificial
		patch boundary, so it can be approximated in $H^1$ by smooth functions
		compactly supported there. Extending these approximants by zero and
		combining them with the interior approximations gives
		$\chi_\delta z\in H_0^1(D_L)$. The preceding bound, together with the
		absolute continuity of the integrals of $|z|^2$ and $|\nabla z|^2$,
		yields $\chi_\delta z\to z$ in $H^1(D_L)$, so $z\in H_0^1(D_L)$.
		We may thus subtract the weak equations and test with $z$ to obtain
		\[
		\lambda\|\nabla z\|_{L^2(D_L)}^2
		\le\int_{D_L}A\nabla z\cdot\nabla z\,dy=0.
		\]
		The zero trace then implies $z=0$. Uniqueness shows that the
		convergences in Step~2 hold for the whole family, and the limiting
		estimates give \eqref{eq:finite-W-bound}.
	\end{proof}
	
	\subsection{Pressure recovery and proof of Theorem~\ref{thm:finite-linear}}
	\label{subsec:pressure}
	
	We recover the pressure associated with $\widetilde w$ from
	\eqref{eq:linear-first-order}, with its value at $e_L$ prescribed by
	\eqref{eq:finite-linear-pressure-normalization}.
	
	\begin{proposition}\label{prop:pressure-recovery}
		Let $\widetilde w$ be given by Proposition~\ref{prop:slit-limit}.
		There is a unique pressure $\widetilde p$, continuous across
		$\Gamma_L^*$, with
		$\widetilde p-p_0\in C^{1,\alpha;(\beta-1,\beta)}_{(-\alpha),\mathrm{pm}}(D_L)$,
		satisfying \eqref{eq:linear-first-order} and the pressure condition
		\eqref{eq:finite-linear-pressure-normalization}.
		Moreover,
		\begin{equation}\label{eq:pressure-regularity}
			\|\widetilde p-p_0\|_{1,\alpha;(\beta-1,\beta);D_L,\mathrm{pm}}^{(-\alpha)}
			+\|\nabla\widetilde p\|_{L^2(D_L)}\le C\kappa_g,
		\end{equation}
		where $C$ depends only on $\alpha,\beta,A_b,\ell,m_\pm$
		and is independent of $L$.
	\end{proposition}
	\begin{proof}
		\noindent\textbf{Step 1. The path integral.}
		Set $\boldsymbol F=(F_1,F_2)=A\nabla\widetilde w$.
		The identity $\operatorname{div}\boldsymbol F=0$ allows us to recover
		the pressure by integration. Starting at $e_L$, we follow the outlet
		to $(L,y_2)$ and then the horizontal segment to $(y_1,y_2)$, and define
		\begin{equation}\label{eq:pressure-path}
			\widetilde p(y_1,y_2)=p_0
			-\int_0^{y_2}F_1(L,\xi)\,d\xi
			-\int_{y_1}^LF_2(s,y_2)\,ds,\qquad y_2\in I.
		\end{equation}
		Write $N=\|\widetilde w\|_{1,\alpha;(\beta);D_L,\mathrm{pm}}^{(-\alpha)}$.
		The coefficient bound and \eqref{eq:finite-W-bound} give
		\begin{equation}\label{eq:pressure-flux-point}
			|\boldsymbol F(z)|\le CN\Delta_z^{-\beta}d(z)^{\alpha-1}.
		\end{equation}
		The distance weight satisfies, for almost every $s$,
		\[
		d\bigl((s,y_2)\bigr)^{\alpha-1}\le1+|s-\ell|^{\alpha-1}.
		\]
		The singularity at $s=\ell$ is integrable because $\alpha>0$.
		For $|s-\ell|\ge1$, we have $d((s,y_2))=1$.
		Combining this identity with $\alpha>0$ and $\beta>1$ gives
		\[
		\sup_{y_2\in I^+\cup I^-}
		\int_{y_1}^L(1+s)^{-\beta}d\bigl((s,y_2)\bigr)^{\alpha-1}\,ds
		\le C(1+y_1)^{1-\beta},\qquad 0\le y_1<L.
		\]
		Since $L-\ell>1$, we have $d((L,\xi))=1$ for $-m_-<\xi<m_+$,
		and \eqref{eq:pressure-flux-point} yields
		\[
		\int_{-m_-}^{m_+}|\boldsymbol F(L,\xi)|\,d\xi\le CN(1+L)^{-\beta}.
		\]
		These bounds show that \eqref{eq:pressure-path} is well defined and
		yields
		\begin{equation}\label{eq:pressure-supremum}
			|\widetilde p(y)-p_0|
			\le CN\left((1+L)^{-\beta}
			+\int_{y_1}^L(1+s)^{-\beta}d\bigl((s,y_2)\bigr)^{\alpha-1}\,ds\right)
			\le C N(1+y_1)^{1-\beta}.
		\end{equation}
		
		\par\medskip\noindent\textbf{Step 2. Verification of \eqref{eq:linear-first-order}.}
		Differentiating \eqref{eq:pressure-path} in $y_1$ gives
		$\widetilde p_{y_1}=F_2$. For the second identity,
		$\widetilde p_{y_2}=-F_1$, take $\xi_a<\xi_b$ in the same interval
		$I^\sigma$ and integrate $\operatorname{div}\boldsymbol F=0$ over
		$(y_1,L)\times(\xi_a,\xi_b)$. By first using interior rectangles
		and then passing to the outlet, we obtain
		\[
		\int_{y_1}^L(F_2(s,\xi_b)-F_2(s,\xi_a))\,ds
		=-\int_{\xi_a}^{\xi_b}(F_1(L,\xi)-F_1(y_1,\xi))\,d\xi.
		\]
		Substituting this identity into the difference of the path formulas
		at $\xi_a$ and $\xi_b$ gives
		\[
		\widetilde p(y_1,\xi_b)-\widetilde p(y_1,\xi_a)
		=-\int_{\xi_a}^{\xi_b}F_1(y_1,\xi)\,d\xi.
		\]
		Thus $\widetilde p_{y_2}=-F_1$, and
		$(\widetilde w,\widetilde p)$ satisfies \eqref{eq:linear-first-order}
		in $D_L^+$ and $D_L^-$.
		
		\par\medskip\noindent\textbf{Step 3. Continuity across the contact line and regularity.}
		The weak identity $\operatorname{div}\boldsymbol F=0$ on $D_L$ and
		the one-sided regularity of $\boldsymbol F$ imply $[F_2]=0$ on
		$\Gamma_L^*$. Taking the two traces in \eqref{eq:pressure-path}
		therefore gives, for $y_1>\ell$,
		\[
		[\widetilde p](y_1,0)
		=-\int_{y_1}^{L}[F_2](s,0)\,ds=0.
		\]
		For fixed $L$, the value bound and $\boldsymbol F\in L^2(D_L)$,
		together with this continuity, give $\widetilde p\in H^1(D_L)$.
		Consequently, \eqref{eq:linear-first-order} holds in the weak sense throughout $D_L$.
		
		Near $e_*$, two points in the same half-annulus of radius $r$ can
		be joined within it by a path of length at most $Cr$, so
		\eqref{eq:pressure-flux-point} bounds the pressure oscillation by
		$CNr^\alpha$. Summation over overlapping dyadic half-annuli then
		gives a limit in each of $D_L^+$ and $D_L^-$; the two limits agree
		by continuity across $\Gamma_L^*$. Away from $e_*$,
		$\nabla \widetilde p=(F_2,-F_1)$ gives the asserted piecewise local
		$C^{1,\alpha}$ regularity, including the Dirichlet corners and $e_L$.
		
		By \eqref{eq:background-bounds} and \eqref{eq:coeff-small},
		integration of the weighted gradient bound along line segments in each
		$D_L^\sigma$ gives
		\[
		\|A\|_{0,\alpha;D_L,\mathrm{pm}}\le C,
		\]
		with $C$ independent of $L$.
		The H\"older product estimate applied to \eqref{eq:linear-first-order},
		together with \eqref{eq:pressure-supremum}, then yields
		\[
		\|\widetilde p-p_0\|_{1,\alpha;(\beta-1,\beta);D_L,\mathrm{pm}}^{(-\alpha)}
		\le CN.
		\]
		Moreover,
		\[
		\|\nabla\widetilde p\|_{L^2(D_L)}
		\le\|A\|_{L^\infty(D_L)}\|\nabla\widetilde w\|_{L^2(D_L)}.
		\]
		Combining these bounds with \eqref{eq:finite-W-bound} gives
		\eqref{eq:pressure-regularity}.
		
		\par\medskip\noindent\textbf{Step 4. Uniqueness.}
		Equation~\eqref{eq:linear-first-order} determines the pressure up to
		one constant in each of $D_L^+$ and $D_L^-$. Continuity across $\Gamma_L^*$ and the prescribed value at $e_L$ fix these
		constants.
	\end{proof}
	
	With the pressure recovered, we can complete the proof of
	Theorem~\ref{thm:finite-linear}.
	
	\begin{proof}[Proof of Theorem~\ref{thm:finite-linear}]
		Proposition~\ref{prop:slit-limit} gives $\widetilde w$, and
		Proposition~\ref{prop:pressure-recovery} gives the pressure $\widetilde p$ with the
		prescribed value at $e_L$. Together they yield a solution of the linear problem
		with the stated regularity and estimates. Any other solution in the stated class
		satisfies \eqref{eq:linear-first-order} almost everywhere in $D_L$.
		Since $\widetilde p\in H^1(D_L)$, taking the divergence shows that
		$\widetilde w$ solves \eqref{eq:slit-problem} weakly. For the difference
		$z$ of two such solutions, the boundary data and continuity give
		$z(y)\to0$ as $y\to e_*$. With the radial cutoff from
		Proposition~\ref{prop:slit-limit},
		\[
		\|z\nabla\chi_\delta\|_{L^2(D_L)}^2
		\le C\sup_{D_L\cap B_{2\delta}(e_*)}|z|^2\longrightarrow0.
		\]
		The same cutoff argument therefore gives $z\in H_0^1(D_L)$.
		Testing the difference equation with $z$
		proves uniqueness of $\widetilde w$. The pressure is then uniquely
		determined by Proposition~\ref{prop:pressure-recovery}.
	\end{proof}
	
	\section{The nonlinear problem in a finite nozzle}
	\label{sec:finite-nonlinear}
	
	We now consider Problem~\ref{prob:lagrangian} on the truncated domain
	$D_L$, where $L>\ell+1$. The inlet and solid-boundary data are retained,
	and an additional condition is prescribed at the artificial outlet.
	Choose a smooth function
	$\chi:\mathbb R\to[0,1]$ with
	$\chi(t)=0$ for $t\le\tfrac14$ and $\chi(t)=1$ for
	$t\ge\tfrac34$, and set
	\[
	\chi_{\mathrm{out}}^\pm(y_2)=\chi\!\left(\frac{\pm y_2}{m_\pm}\right).
	\]
	The outlet value of $w$ is chosen to match the wall slopes and
	vanish near $y_2=0$:
	\begin{equation}\label{eq:artificial-outlet-data}
		g_{\mathrm{out}}(y_2)
		=(W^+)'(L)\chi_{\mathrm{out}}^+(y_2)
		+(W^-)'(L)\chi_{\mathrm{out}}^-(y_2).
	\end{equation}
	Together with the prescribed inlet and solid-boundary values, this
	defines the boundary function
	\begin{equation}\label{eq:finite-nonlinear-boundary-data}
		g(y_1,y_2)=
		\begin{cases}
			(W^\sigma)'(y_1),&(y_1,y_2)\in\Gamma_{\mathrm w}^\sigma,\\
			(\zeta^\sigma)'(y_1),&(y_1,y_2)\in\Gamma_{\mathrm a}^\sigma,\\
			w_{\mathrm{in}}(y_2),&(y_1,y_2)\in I_{\mathrm{in}}^+\cup I_{\mathrm{in}}^-,\\
			g_{\mathrm{out}}(y_2),&(y_1,y_2)\in O_L.
		\end{cases}
	\end{equation}
	The boundary conditions for the finite-nozzle problem are
	\begin{equation}\label{eq:finite-nonlinear-outlet}
		w=g\quad\text{on }\partial D_L,\qquad p(e_L)=p_0.
	\end{equation}
	
	\Needspace{11\baselineskip}
	\begin{problem}[Nonlinear problem in a finite nozzle]\label{prob:finite-nonlinear}
		Find $(w,p)$ on $D_L$, with $p>0$, such that:
		\begin{enumerate}
			\item The pair satisfies \eqref{eq:nonlinear-rotated-system} in
			$D_L^+$ and $D_L^-$, with $A=A(S,B,w,p)$ given by
			\eqref{eq:chen-matrix}, where $S,B$ are determined by
			\eqref{eq:inlet-transported-profiles}.
			\item The contact discontinuity conditions $[w]=[p]=0$ hold on
			$\Gamma_L^*$, and the boundary conditions
			\eqref{eq:finite-nonlinear-outlet} hold.
			\item The conditions \eqref{eq:wp-admissibility} hold with
			$D_\infty$ replaced by $D_L$.
		\end{enumerate}
	\end{problem}
	
	The reference state $\widehat U_0$ in
	\eqref{eq:lagrangian-reference-state} gives
	$A_b(y_2)=\operatorname{diag}(a_b^{11}(y_2),a_b^{22}(y_2))$, which satisfies
	\begin{equation}\label{eq:finite-background-bounds}
		\lambda\le a_b^{11},a_b^{22}\le\Lambda,\qquad
		\|A_b\|_{1,\alpha;I,\mathrm{pm}}\le\Lambda,
	\end{equation}
	for fixed positive constants $\lambda,\Lambda$.
	The size of the prescribed data is measured by
	\begin{equation}\label{eq:finite-data-size}
		\hat\kappa=\|\widehat u_{2,\mathrm{in}}\|_{1,\alpha;I,\mathrm{pm}}
		+\sum_{\sigma\in\{+,-\}}\left(
		\|(W^\sigma)'\|_{1,\alpha;(\beta);[0,\infty)}
		+\|(\zeta^\sigma)'\|_{1,\alpha;(0,\ell)}^{(-\alpha;\{\ell\})}\right).
	\end{equation}
	
	\begin{theorem}[Solvability of the nonlinear problem]\label{thm:finite-nonlinear}
		Suppose that $\alpha,\beta$ and $\widehat U_0$ are as in
		Theorem~\ref{thm:infinite-existence}.
		Then there exist positive constants $\hat\kappa_0,\sigma_0,C$,
		independent of $L$, such that, if $\hat\kappa\le\hat\kappa_0$,
		Problem~\ref{prob:finite-nonlinear} has a solution $(w,p)$ in the
		piecewise weighted $C^{1,\alpha}$ class satisfying
		\begin{equation}\label{eq:nonlinear-final-bound}
			\|w\|_{1,\alpha;(\beta);D_L,\mathrm{pm}}^{(-\alpha)}
			+\|p-p_0\|_{1,\alpha;(\beta-1,\beta);D_L,\mathrm{pm}}^{(-\alpha)}
			\le C\hat{\kappa}.
		\end{equation}
		The solution is unique among solutions of Problem~\ref{prob:finite-nonlinear}
		for which the left-hand side of \eqref{eq:nonlinear-final-bound} is at
		most $\sigma_0$. The constants depend only on
		$\widehat U_0,\ell,m_\pm,\gamma,\alpha,\beta$.
	\end{theorem}
	
	\subsection{Preliminary estimates}
	\label{subsec:finite-preliminaries}
	
	The fixed point construction requires a compact set of inputs and
	uniform control of the coefficients they determine. A further estimate
	for differences will give uniqueness. For fixed $L$, the input space is
	$Y=C^0_{\mathrm{pm}}(\overline{D_L})\times C^0_{\mathrm{pm}}(\overline{D_L})$,
	with norm
	\[
	\|(v,h)\|_{Y}
	=\|v\|_{0,0;\overline{D_L},\mathrm{pm}}
	+\|h\|_{0,0;\overline{D_L},\mathrm{pm}}.
	\]
	To include the estimates needed later on the semi-infinite nozzle,
	write $D_R=D_L$ for $R=L$ and $D_R=D_\infty$ for $R=\infty$.
	The first lemma records the compactness and integrability consequences
	of the weighted bounds.
	
	\begin{lemma}\label{lem:weighted-calculus}
		For $R=L$ or $R=\infty$ and $1\le s<2/(1-\alpha)$,
		\begin{equation}\label{eq:uniform-weighted-integrability}
			\begin{aligned}
				\|v\|_{W^{1,s}_{\mathrm{pm}}(D_R)}
				&\le C\|v\|_{1,\alpha;(\beta);D_R,\mathrm{pm}}^{(-\alpha)},\\
				\|\nabla h\|_{L^s_{\mathrm{pm}}(D_R)}
				&\le C\|h\|_{1,\alpha;(\beta-1,\beta);D_R,\mathrm{pm}}^{(-\alpha)}.
			\end{aligned}
		\end{equation}
		The constant $C$ may depend on $s$, but is independent of $R$.
		For fixed finite $L$ and $K_0>0$,
		the pairs $(v,h)\in Y$ satisfying
		\[
		\|v\|_{1,\alpha;(\beta);D_L,\mathrm{pm}}^{(-\alpha)}
		+\|h\|_{1,\alpha;(\beta-1,\beta);D_L,\mathrm{pm}}^{(-\alpha)}\le K_0
		\]
		form a compact convex subset of $Y$.
	\end{lemma}
	\begin{proof}
		The derivative bounds in the weighted norms are integrable because
		$s(1-\alpha)<2$ and $\beta s>1$. Indeed,
		\[
		\int_{D_R}\Delta_y^{-\beta s}d(y)^{-s(1-\alpha)}\,dy
		\le C\left(\int_0^1 r^{1-s(1-\alpha)}\,dr
		+\int_0^R(1+\xi)^{-\beta s}\,d\xi\right)\le C.
		\]
		The second integral also controls the function value of $v$, giving
		\eqref{eq:uniform-weighted-integrability}.
		
		For compactness, fix $L<\infty$. On either convex half-rectangle
		$G=D_L^\sigma$, the weighted gradient bound gives uniform H\"older
		control up to the boundary. More precisely, for
		$N=\|f\|_{1,\alpha;G}^{(-\alpha)}$ and $y,z\in G$ with $0<a=|y-z|<1$,
		\[
		|f(y)-f(z)|\le N\int_0^a d(y+\xi\omega)^{\alpha-1}\,d\xi
		\le CNa^\alpha,\qquad \omega=\frac{z-y}{a}.
		\]
		Here $d(y+\xi\omega)^{\alpha-1}\le1+|\xi-\xi_*|^{\alpha-1}$,
		where $\xi_*$ is the projection parameter of $e_*$ onto the
		supporting line; for $a\ge1$, the supremum bound suffices.
		Thus every sequence $(v_j,h_j)$ in the stated norm ball is uniformly
		bounded and equicontinuous on each of $\overline{D_L^+}$ and
		$\overline{D_L^-}$. Arzel\`a--Ascoli gives a subsequence converging
		to $(v,h)$ in $Y$. The local H\"older bounds for the gradients give,
		after a further diagonal extraction, locally uniform convergence
		to $Dv,Dh$ up to each side boundary away from $e_*$.
		Passing to the limit in the weighted pointwise and two-point
		inequalities yields
		\[
		\|v\|_{1,\alpha;(\beta);D_L,\mathrm{pm}}^{(-\alpha)}+\|h\|_{1,\alpha;(\beta-1,\beta);D_L,\mathrm{pm}}^{(-\alpha)}
		\le\liminf_{j\to\infty}
		\bigl(\|v_j\|_{1,\alpha;(\beta);D_L,\mathrm{pm}}^{(-\alpha)}
		+\|h_j\|_{1,\alpha;(\beta-1,\beta);D_L,\mathrm{pm}}^{(-\alpha)}\bigr)\le K_0.
		\]
		The norm ball is therefore compact in $Y$, and the triangle
		inequality gives convexity.
	\end{proof}
	To apply the linear theory to an input pair $(w,p)$, the associated
	matrix must remain close to $A_b$. The prescribed inlet data first give,
	by \eqref{eq:finite-data-size} and the positive lower bound for
	$\widehat u_{1,\mathrm{in}}$,
	\begin{equation}\label{eq:inlet-pullback-smallness}
		\|\widehat u_{2,\mathrm{in}}\|_{1,\alpha;I,\mathrm{pm}}
		+\|w_{\mathrm{in}}\|_{1,\alpha;I,\mathrm{pm}}\le C\hat{\kappa}.
	\end{equation}
	At the reference state $\widehat U_0$, the identities
	$w=0$, $p=p_0$, $q=\widehat u_{1,\mathrm{in}}$ and
	$c^2=\gamma p_0/\widehat\rho_{\mathrm{in}}$ reduce
	\eqref{eq:chen-matrix} to
	\begin{equation}\label{eq:reference-matrix}
		A_b(y_2)=\operatorname{diag}\!\left(
		\widehat u_{1,\mathrm{in}},
		\frac{\widehat\rho_{\mathrm{in}}^{\,2}\widehat u_{1,\mathrm{in}}^{\,3}}{1-\mathrm{Ma}_0^2}\right),
	\end{equation}
	where $\mathrm{Ma}_0^2=\widehat\rho_{\mathrm{in}}\widehat u_{1,\mathrm{in}}^{\,2}/(\gamma p_0)$
	is the squared Mach number of the reference state. The next lemma
	shows that this ellipticity and the required coefficient bounds persist
	for small inputs, including their dependence on the inlet pressure trace.
	
	\begin{lemma}\label{lem:Euler-coefficients}
		There exist positive constants $\sigma_b,c_b,C_b$, independent of $R$, such that if
		$\hat{\kappa}\le\sigma_b$ and
		\[
		\|w\|_{1,\alpha;D_R,\mathrm{pm}}^{(-\alpha)}
		+\|p-p_0\|_{1,\alpha;D_R,\mathrm{pm}}^{(-\alpha)}\le\sigma_b,
		\]
		then
		\begin{equation}\label{eq:Euler-state-bounds}
			p,\ q,\ c^2-q^2\ge c_b,\qquad \frac{q^2}{c^2}\le1-\frac{\delta_b}{2},
		\end{equation}
		and
		\[
		c_b|\xi|^2\le\xi^{\mathsf T}A(y)\xi\le C_b|\xi|^2.
		\]
		Moreover,
		\begin{equation}\label{eq:Euler-coefficient-smallness}
			\|A-A_b\|_{1,\alpha;D_R,\mathrm{pm}}^{(-\alpha)}
			\le C_b\left(
			\|w\|_{1,\alpha;D_R,\mathrm{pm}}^{(-\alpha)}
			+\|p-p_0\|_{1,\alpha;D_R,\mathrm{pm}}^{(-\alpha)}+\hat{\kappa}\right).
		\end{equation}
		For two pairs satisfying the preceding assumptions and having the
		same inlet data, let $A_1,A_2$ denote the corresponding matrices. Then
		\begin{equation}\label{eq:Euler-coefficient-Lipschitz}
			\|A_1-A_2\|_{L^\infty_{\mathrm{pm}}(D_R)}
			\le C_b\left(
			\|w_1-w_2\|_{0,0;\overline{D_R},\mathrm{pm}}
			+\|p_1-p_2\|_{0,0;\overline{D_R},\mathrm{pm}}\right).
		\end{equation}
		For finite $R=L$, if $w(L,y_2)=0$ for
		$|y_2|\le\delta_{\mathrm{out}}$, then
		$a^{12}(L,y_2)=a^{21}(L,y_2)=0$ there.
	\end{lemma}
	\begin{proof}
		\noindent\textbf{Step 1. The pressure trace and the coefficient functions.}
		The inlet pressure enters the transported profiles in
		\eqref{eq:inlet-transported-profiles}. Since
		$d((0,y_2))\ge\min\{1,\ell\}>0$, its trace satisfies
		\[
		\|p(0,\cdot)-p_0\|_{1,\alpha;I,\mathrm{pm}}
		\le C\|p-p_0\|_{1,\alpha;D_R,\mathrm{pm}}^{(-\alpha)}.
		\]
		Its extension independent of $y_1$ obeys the same corner-weighted
		bound. The lower bound $\widehat\rho_{\mathrm{in}}\ge\rho_*>0$
		from \eqref{eq:physical-background}, together with
		$p_0/2\le p\le3p_0/2$ for $\sigma_b\le p_0/2$, keeps the arguments
		of the coefficient functions in a fixed compact range.
		Differentiating \eqref{eq:coefficient-functions-SB} and
		\eqref{eq:inlet-transported-profiles}, and using
		\eqref{eq:inlet-pullback-smallness}, yields
		\[
		\left\|c^2-\frac{\gamma p_0}{\widehat\rho_{\mathrm{in}}}\right\|_{1,\alpha;D_R,\mathrm{pm}}^{(-\alpha)}
		+\|q^2-\widehat u_{1,\mathrm{in}}^{\,2}\|_{1,\alpha;D_R,\mathrm{pm}}^{(-\alpha)}
		\le C_b\left(\|p-p_0\|_{1,\alpha;D_R,\mathrm{pm}}^{(-\alpha)}+\hat{\kappa}\right).
		\]
		For the H\"older seminorms, integration along line segments in either
		convex side gives
		\[
		[f]_{\alpha;G}\le C\|f\|_{1,\alpha;G}^{(-\alpha)},\qquad
		G=D_R^+\ \text{or}\ D_R^-.
		\]
		The mixed terms in the differentiated products are then controlled by
		\[
		d(y,z)\frac{|f(y)-f(z)|\,|Dg(z)|}{|y-z|^\alpha}
		\le [f]_{\alpha;G}\|g\|_{1,\alpha;G}^{(-\alpha)}
		d(y,z)d(z)^{\alpha-1}
		\le [f]_{\alpha;G}\|g\|_{1,\alpha;G}^{(-\alpha)}.
		\]
		Together with the chain rule, this also controls smooth compositions
		on the compact range above.
		
		\par\medskip\noindent\textbf{Step 2. Ellipticity and smallness.}
		The inequalities in \eqref{eq:physical-background} give
		$\widehat u_{1,\mathrm{in}}\ge u_*>0$ and
		$\widehat\rho_{\mathrm{in}}\widehat u_{1,\mathrm{in}}^2/(\gamma p_0)
		\le1-\delta_b$. The estimates for $c^2$ and $q^2$ in Step~1 therefore
		give \eqref{eq:Euler-state-bounds} after decreasing $\sigma_b$.
		Together with the resulting upper bounds for $p,c,q$, these inequalities
		keep the denominators $c^2$ and $c^2-q^2$ in \eqref{eq:chen-matrix}
		away from zero, uniformly in $R$. The formula for $A$ then gives
		positive lower bounds for $a^{11}$ and $\det A$, so the identity
		\[
		\xi^{\mathsf T}A\xi
		=a^{11}\left(\xi_1+\frac{a^{12}}{a^{11}}\xi_2\right)^2
		+\frac{\det A}{a^{11}}\xi_2^2
		\]
		gives uniform ellipticity.
		To compare $A$ with $A_b$, keep $\widehat\rho_{\mathrm{in}}$ and
		$\widehat u_{1,\mathrm{in}}$ fixed and apply the mean value theorem
		to \eqref{eq:chen-matrix} in the variables
		$w,p,p(0,y_2),\widehat u_{2,\mathrm{in}}$ about $(0,p_0,p_0,0)$.
		The matrix at this reference point is \eqref{eq:reference-matrix}.
		The chain rule and the preceding weighted product bound now give
		\eqref{eq:Euler-coefficient-smallness}.
		
		\par\medskip\noindent\textbf{Step 3. Differences and the outlet.}
		For two inputs with the same inlet data, the trace difference satisfies
		\[
		\|p_1(0,\cdot)-p_2(0,\cdot)\|_{0,0;I,\mathrm{pm}}
		\le\|p_1-p_2\|_{0,0;\overline{D_R},\mathrm{pm}}.
		\]
		Applying the mean value theorem to
		\eqref{eq:inlet-transported-profiles}, \eqref{eq:coefficient-functions-SB}
		and \eqref{eq:chen-matrix}, with this trace bound, proves
		\eqref{eq:Euler-coefficient-Lipschitz}. At the outlet, the factor $w$
		in $a^{12}=a^{21}$ gives the asserted vanishing.
	\end{proof}
	
	The uniqueness argument uses the first-order system for the
	difference of two solutions, with the coefficient difference placed
	in the source term. The following estimate controls both unknowns
	by this source, with a constant independent of the nozzle length.
	
	\begin{lemma}\label{lem:weak-difference}
		Let $R=L$ or $R=\infty$, and let $A$ be a
		bounded measurable symmetric matrix on $D_R$ with fixed ellipticity
		bounds $0<\lambda\le\Lambda$.
		Suppose that $z\in H_0^1(D_R)$ and
		$h\in H^1_{\mathrm{loc}}(D_R)\cap C^0_{\mathrm{pm}}(\overline{D_R})$ satisfy
		\begin{equation}\label{eq:first-order-difference-gradient}
			\nabla_y^\perp h=A\nabla_y z+H,
			\qquad |H(y)|\le N\Delta_y^{-\beta}d(y)^{\alpha-1},
		\end{equation}
		where the equation holds in distributions on $D_R$, the source $H:D_R\to\mathbb R^2$ is measurable, and $N\ge0$ is a constant.
		For finite $R$, assume $h(e_R)=0$, where $e_R=(R,0)$. For $R=\infty$, assume that $h(y_1,\cdot)$
		tends uniformly to zero as $y_1\to\infty$. Then
		\begin{equation}\label{eq:weak-pair-bound}
			\|z\|_{0,0;\overline{D_R},\mathrm{pm}}
			+\|h\|_{0,0;\overline{D_R},\mathrm{pm}}\le CN,
		\end{equation}
		where $C$ depends only on $\ell,m_\pm,\lambda,\Lambda,\alpha,\beta$
		and is independent of $R$.
	\end{lemma}
	\begin{proof}
		The assumptions $\alpha>0$ and $\beta>1$ ensure $H\in L^2(D_R)$.
		Taking the divergence of \eqref{eq:first-order-difference-gradient}
		therefore gives, first for compactly supported smooth test functions
		and then by density,
		\begin{equation}\label{eq:weak-difference-equation}
			\int_{D_R}A\nabla z\cdot\nabla\phi\,dy
			=-\int_{D_R}H\cdot\nabla\phi\,dy,
			\qquad \phi\in H_0^1(D_R).
		\end{equation}
		\noindent\textbf{Step 1. Localized energy.}
		For $a\in[0,R]$, or $a\ge0$ when $R=\infty$, set
		$\chi_a(y_1)=e^{-\delta|y_1-a|}$, where $\delta>0$ is a fixed
		constant to be chosen below. Write $\|\cdot\|_2=\|\cdot\|_{L^2(D_R)}$.
		The one-dimensional Poincar\'e inequality between the nozzle walls
		applies to $\chi_a z$, since $\chi_a$ depends only on $y_1$ and
		$z$ has zero wall traces. Hence
		\[
		\|\chi_a z\|_2\le C_\perp\|\chi_a\partial_{y_2}z\|_2
		\le C_\perp\|\chi_a\nabla z\|_2,
		\]
		where $C_\perp$ depends only on $m_+,m_-$. The function
		$\chi_a^2z$ belongs to $H_0^1(D_R)$, and
		\[
		|\chi_a'|\le\delta\chi_a.
		\]
		Testing \eqref{eq:weak-difference-equation} with $\chi_a^2z$ therefore yields
		\[
		\begin{aligned}
			\lambda\|\chi_a\nabla z\|_2^2
			&\le 2\Lambda C_\perp\delta\|\chi_a\nabla z\|_2^2\\
			&\quad+(1+2C_\perp\delta)
			\|\chi_aH\|_2\|\chi_a\nabla z\|_2.
		\end{aligned}
		\]
		Fix $\delta$ so that $2\Lambda C_\perp\delta\le\lambda/2$.
		Then
		\[
		\|\chi_a\nabla z\|_2\le C\|\chi_aH\|_2.
		\]
		
		To estimate the weighted source, split the integral at $B_1(e_*)$,
		where $d(y)=|y-e_*|$; outside this disk, $d(y)=1$. The inequality
		\[
		1+a\le(1+y_1)(1+|y_1-a|)
		\]
		then separates the dependence on $a$ from the integrable tail:
		\[
		\begin{aligned}
			\|\chi_aH\|_2^2
			&\le N^2(1+a)^{-2\beta}
			\int_{D_R}e^{-2\delta|y_1-a|}(1+|y_1-a|)^{2\beta}
			d(y)^{2\alpha-2}\,dy\\
			&\le CN^2(1+a)^{-2\beta}
			\left(\int_0^1 r^{2\alpha-1}\,dr
			+\int_{\mathbb R}e^{-2\delta|\xi|}
			(1+|\xi|)^{2\beta}\,d\xi\right)\\
			&\le CN^2(1+a)^{-2\beta}.
		\end{aligned}
		\]
		Here the bounded width of $D_R$ controls the second integral. Thus
		\begin{equation}\label{eq:localized-difference-energy}
			\|\chi_a\nabla z\|_2
			\le C\|\chi_aH\|_2\le CN(1+a)^{-\beta},
		\end{equation}
		with $C$ independent of $a,R$.
		
		\par\medskip\noindent\textbf{Step 2. Local gradient estimates.}
		Fix $0<r_*<\tfrac12$ sufficiently small in terms of $\ell,m_+,m_-$,
		independently of $R$, and set
		\[
		E_a=D_R\cap\{|y_1-a|<r_*\},\qquad
		E_a^*=D_R\cap\{|y_1-a|<2r_*\}.
		\]
		The localized energy bound can be improved to an $L^s$ estimate
		with $s>2$. For $b\in\overline{D_R}$ and small $r$, write
		$G_r(b)=D_R\cap B_r(b)$. Testing
		\eqref{eq:weak-difference-equation} with $\eta^2(z-\bar z)$,
		where $\eta$ is a cutoff between $B_r(b)$ and $B_{2r}(b)$, gives
		\[
		\int_{G_r(b)}|\nabla z|^2\,dy
		\le \frac C{r^2}\int_{G_{2r}(b)}|z-\bar z|^2\,dy
		+C\int_{G_{2r}(b)}|H|^2\,dy.
		\]
		Here $\bar z$ is the average over $B_{2r}(b)$ in an interior ball
		and is zero in a Dirichlet boundary neighborhood. A ball crossing
		the contact line but no exterior boundary is an interior ball.
		
		Fix $s_0\in(1,2)$. The scale-invariant Sobolev--Poincar\'e
		inequality applies on disks, half-disks and quadrants, with mean
		subtraction or zero Dirichlet trace as appropriate. At an interior
		point of the slit, we apply it separately on the two half-disks;
		at the inlet endpoint of the slit, we use the two quadrants. At $e_*$, set
		$B_{2r}^{\sigma}=B_{2r}(e_*)\cap D_R^{\sigma}$ for
		$\sigma\in\{+,-\}$. On each half-disk, $z$ has zero trace on the
		negative horizontal radius, whose relative length is unchanged by
		scaling. Hence
		\[
		\left(\fint_{B_{2r}^{\sigma}}|z|^2\,dy\right)^{1/2}
		\le Cr\left(\fint_{B_{2r}^{\sigma}}
		|\nabla z|^{s_0}\,dy\right)^{1/s_0}.
		\]
		Combining the two inequalities gives the same bound on the entire
		cut disk. This subdivision is used only for Sobolev--Poincar\'e.
		The cutoff test remains on the whole cut disk, with no boundary
		condition imposed on the contact line.
		
		A ball meeting the boundary is reduced to one of these geometries
		by recentering at a boundary point and enlarging its radius by a
		fixed factor. Near a corner or $e_*$, that point is used as the
		center. Consequently, for a fixed geometric constant $c_*>1$,
		\[
		\left(\fint_{G_r(b)}|\nabla z|^2\,dy\right)^{1/2}
		\le C\left(\fint_{G_{c_*r}(b)}|\nabla z|^{s_0}\,dy\right)^{1/s_0}
		+C\left(\fint_{G_{c_*r}(b)}|H|^2\,dy\right)^{1/2}.
		\]
		The area bounds are uniform:
		\[
		c r^2\le|G_r(b)|\le Cr^2.
		\]
		Extend $|\nabla z|$ and $|H|$ by zero outside $D_R$.
		Every sufficiently small ambient ball meeting $D_R$ is contained in a
		ball of comparable radius centered in $\overline{D_R}$. The area bounds
		therefore give the same reverse H\"older inequality for the zero extensions,
		with a fixed enlargement of the ambient ball.
		Choose fixed nested cubes and a radius $\rho>0$ so that all required
		enlargements of balls centered in the inner cube lie in the outer
		cube $Q$. Let $M$ be the Hardy--Littlewood maximal operator and let
		$M_\rho$ denote its restriction to radii at most $\rho$. Set
		$q=2/s_0$, $g=|\nabla z|^{s_0}\mathbf1_Q$, and
		$f=[M(|H|^2\mathbf1_Q)]^{s_0/2}\mathbf1_Q$.
		The reverse H\"older inequality gives, on the inner cube,
		\[
		M_\rho(g^q)\le C(Mg)^q+C f^q\le C[M(g+f)]^q,
		\]
		where the last inequality follows from Lebesgue differentiation.
		This is the hypothesis of the inhomogeneous self-improvement
		theorem~\cite[p.~409]{Stredulinsky}. For $s>2$, the boundedness
		of $M$ on $L^{s/2}$ gives
		\[
		\|f\|_{L^{s/s_0}(Q)}\le C\|H\|_{L^s(Q\cap D_R)}^{s_0}.
		\]
		The resulting local Meyers-type estimate~\cite{Meyers} has an exponent
		$s_*>2$ determined by the ellipticity bounds and these fixed geometries.
		A cover of $E_a$ by balls of a fixed radius can be chosen so that
		their required enlargements meet $D_R$ only inside $E_a^*$.
		Combining the local estimates over this cover gives
		\begin{equation}\label{eq:local-difference-Meyers}
			\|\nabla z\|_{L^s(E_a)}
			\le C\bigl(\|\nabla z\|_{L^2(E_a^*)}
			+\|H\|_{L^s(E_a^*)}\bigr),
			\qquad 2<s<\min\{s_*,2/(1-\alpha)\}.
		\end{equation}
		Here $C$ may depend on $s$, but is independent of $a$ and $R$.
		The number of balls is uniformly bounded because the strip width
		and $r_*$ are fixed.
		
		Fix $s$ in this range. The assumed pointwise bound on $H$ is
		integrable to this power near $e_*$, since
		\[
		\int_{E_a^*}d(y)^{-s(1-\alpha)}\,dy
		\le |E_a^*|+2\pi\int_0^1 r^{1-s(1-\alpha)}\,dr
		\le C.
		\]
		On $E_a^*$ the ratio $(1+y_1)/(1+a)$ has fixed positive bounds,
		and
		\[
		\chi_a\ge e^{-2\delta r_*}.
		\]
		Hence the source bound and
		\eqref{eq:localized-difference-energy} give
		\[
		\|H\|_{L^s(E_a^*)}+\|\nabla z\|_{L^2(E_a^*)}
		\le CN(1+a)^{-\beta}.
		\]
		The same estimate controls $\nabla h$ through
		\eqref{eq:first-order-difference-gradient}, which gives
		\[
		|\nabla h|\le\Lambda|\nabla z|+|H|.
		\]
		Together with \eqref{eq:local-difference-Meyers}, this yields
		\begin{equation}\label{eq:local-pair-gradient}
			\|\nabla z\|_{L^s(E_a)}+\|\nabla h\|_{L^s(E_a)}
			\le CN(1+a)^{-\beta}.
		\end{equation}
		
		\par\medskip\noindent\textbf{Step 3. Function values and the pressure constant.}
		The gradient estimate gives function-value bounds on each side.
		Indeed, $E_a^\sigma=E_a\cap D_R^\sigma$ is a rectangle of fixed
		height and width between $r_*$ and $2r_*$, so its Poincar\'e and
		Morrey constants are uniform.
		Write $\operatorname{osc}_Gh=\sup_{y,y'\in G}|h(y)-h(y')|$.
		The zero outer-wall trace of $z$ and
		\eqref{eq:local-pair-gradient} imply
		\[
		\|z\|_{L^\infty(E_a^\sigma)}
		\le C\|\nabla z\|_{L^s(E_a^\sigma)}
		\le CN(1+a)^{-\beta}
		\]
		and
		\[
		\operatorname{osc}_{\overline{E_a^\sigma}}h
		\le C\|\nabla h\|_{L^s(E_a^\sigma)}
		\le CN(1+a)^{-\beta}.
		\]
		In particular, $z$ has a continuous representative on each of
		$\overline{D_R^+}$ and $\overline{D_R^-}$.
		The remaining pressure constant is determined by its prescribed
		outlet value when $R<\infty$, or by its limit at infinity.
		
		For finite $R$, the $H^1_{\mathrm{loc}}(D_R)$ regularity of $h$ gives equal traces
		on the interior contact line. Their one-sided continuity extends this equality to
		$e_R$, where the prescribed value gives $h^+(e_R)=h^-(e_R)=0$. For $y\in D_R^\sigma$, choose $n\ge1$
		with $(n-1)r_*<R-y_1\le nr_*$ and set
		\[
		a_j=y_1+jr_*\quad(0\le j<n),\qquad a_n=R.
		\]
		The horizontal chain stays in the chosen side. Telescoping from $y$
		to the outlet gives
		\[
		h(y)=\sum_{j=0}^{n-1}
		\bigl(h(a_j,y_2)-h(a_{j+1},y_2)\bigr)
		+h(R,y_2)-h(R,0).
		\]
		The two points in each difference lie in the closure of one of the
		rectangles above. Thus
		\[
		\begin{aligned}
			|h(y)|
			&\le CN\left(\sum_{j=0}^{\infty}(1+y_1+jr_*)^{-\beta}
			+(1+R)^{-\beta}\right)\\
			&\le CN\left(2(1+y_1)^{-\beta}
			+\frac1{r_*}\int_{y_1}^{\infty}(1+\xi)^{-\beta}\,d\xi\right)
			\le CN(1+y_1)^{1-\beta}.
		\end{aligned}
		\]
		For $R=\infty$, end the same chain at $a_n=y_1+nr_*$. The remaining value
		$h(a_n,y_2)$ tends uniformly to zero by hypothesis. Letting $n\to\infty$ gives
		the same estimate. Together with the bound for $z$, this
		proves~\eqref{eq:weak-pair-bound}.
	\end{proof}
	
	\subsection{Existence}
	\label{subsec:finite-existence}
	
	A candidate pair $(w,p)$ determines the coefficients of a linear
	problem, whose solution defines the fixed point map. The estimates of
	Theorem~\ref{thm:finite-linear} and Lemma~\ref{lem:Euler-coefficients}
	keep this map in a compact convex subset of $Y$. It remains to choose
	that set and prove continuity of the map.
	
	\begin{proof}[Proof of Theorem~\ref{thm:finite-nonlinear}: existence]
		\noindent\textbf{Step 1. The set of inputs.}
		For the boundary data $g$ in
		\eqref{eq:finite-nonlinear-boundary-data} and a constant $K_0\ge1$
		to be chosen below, set
		\begin{equation}\label{eq:Schauder-set}
			\begin{split}
				K=\{(w,p)\in Y:\ &
				w\in C^{1,\alpha;(\beta)}_{(-\alpha),\mathrm{pm}}(D_L),\quad
				p-p_0\in C^{1,\alpha;(\beta-1,\beta)}_{(-\alpha),\mathrm{pm}}(D_L),\\
				&[w]=[p]=0\text{ on }\Gamma_L^*,\\
				&w=g\text{ on }\partial D_L,\quad p(e_L)=p_0,\\
				&\|w\|_{1,\alpha;(\beta);D_L,\mathrm{pm}}^{(-\alpha)}
				+\|p-p_0\|_{1,\alpha;(\beta-1,\beta);D_L,\mathrm{pm}}^{(-\alpha)}
				\le K_0\hat{\kappa}\}.
			\end{split}
		\end{equation}
		The boundary and contact conditions, including the prescribed pressure
		at $e_L$, are closed affine constraints in $Y$. Lemma~\ref{lem:weighted-calculus},
		applied to $(w,p-p_0)$, therefore shows that $K$ is compact and convex.
		
		To find an element of $K$, first note that $g$ satisfies the boundary
		assumptions of the linear problem. The inlet compatibility and common
		trailing-edge value follow from \eqref{eq:lagrangian-compatibility},
		while \eqref{eq:artificial-outlet-data} matches the wall values at the
		outlet corners and vanishes for $|y_2|\le\delta_{\mathrm{out}}$.
		These are precisely \eqref{eq:compat-left}--\eqref{eq:linear-outlet-vanishing}.
		The estimates \eqref{eq:inlet-pullback-smallness} and
		\eqref{eq:artificial-outlet-data}, together with the solid-boundary
		terms in \eqref{eq:finite-data-size}, also give
		\begin{equation}\label{eq:finite-boundary-data-bound}
			\kappa_g=\|g\|_{1,\alpha;(\beta);\partial D_L,\mathrm{pm}}^{(-\alpha)}
			\le C\hat{\kappa},
		\end{equation}
		where $C$ is independent of $L$.
		
		Proposition~\ref{prop:slit-limit} with matrix $A_b$ and boundary data
		$g$ now gives a solution $w^{(0)}$ satisfying
		\[
		\|w^{(0)}\|_{1,\alpha;(\beta);D_L,\mathrm{pm}}^{(-\alpha)}
		\le C_{\mathrm{init}}\hat{\kappa},
		\]
		with $C_{\mathrm{init}}$ independent of $L$. Thus
		$(w^{(0)},p_0)\in K$ for $K_0\ge C_{\mathrm{init}}$.
		
		\par\medskip\noindent\textbf{Step 2. The solution map.}
		Fix $0<\tau\le\tau_0$ as in Proposition~\ref{prop:uniform}.
		For coefficients satisfying the corresponding perturbation bound,
		Propositions~\ref{prop:slit-limit} and~\ref{prop:pressure-recovery},
		together with \eqref{eq:finite-boundary-data-bound}, bound the linear
		solution pair by $C_{\mathrm{map}}\hat{\kappa}$ independently of $L$.
		To make these estimates applicable to every input in $K$, choose
		\[
		K_0>\max\{1,C_{\mathrm{init}},C_{\mathrm{map}}\},\qquad
		0<\sigma_0\le\sigma_b,\qquad C_b\sigma_0\le\tau/2,
		\]
		and then choose $\hat{\kappa}_0>0$ such that
		\begin{equation}\label{eq:fixed-point-smallness}
			K_0\hat{\kappa}_0\le\sigma_0,\qquad C_b\hat{\kappa}_0\le\tau/2,\qquad
			\hat{\kappa}_0\le1.
		\end{equation}
		The choices of $K_0,\sigma_0,\hat{\kappa}_0$ are independent of $L$.
		
		For $(w,p)\in K$ with $\hat{\kappa}\le\hat{\kappa}_0$, the inlet trace
		$p(0,\cdot)$ fixes $S,B$ through \eqref{eq:inlet-transported-profiles}.
		The quantities $c,q,A$ are then evaluated from these profiles and
		$(w,p)$. Since $\Delta_y\ge1$ and $\beta>1$,
		\[
		\|w\|_{1,\alpha;D_L,\mathrm{pm}}^{(-\alpha)}
		+\|p-p_0\|_{1,\alpha;D_L,\mathrm{pm}}^{(-\alpha)}\le K_0\hat{\kappa}.
		\]
		The coefficient estimate \eqref{eq:Euler-coefficient-smallness} then gives
		\[
		\|A-A_b\|_{1,\alpha;D_L,\mathrm{pm}}^{(-\alpha)}
		\le C_b(K_0+1)\hat{\kappa}
		\le C_b(\sigma_0+\hat{\kappa}_0)\le\tau.
		\]
		Since $w=g_{\mathrm{out}}=0$ on the outlet near $e_L$,
		Lemma~\ref{lem:Euler-coefficients} also gives
		\[
		a^{12}(L,y_2)=a^{21}(L,y_2)=0,
		\qquad |y_2|\le\delta_{\mathrm{out}}.
		\]
		The coefficients therefore satisfy all the hypotheses of the linear
		problem. Propositions~\ref{prop:slit-limit} and~\ref{prop:pressure-recovery}
		give a pair $(\widetilde w,\widetilde p)$ satisfying
		\[
		\nabla_y^\perp\widetilde p=A\nabla_y\widetilde w,
		\qquad \widetilde w=g\text{ on }\partial D_L,
		\qquad \widetilde p(e_L)=p_0.
		\]
		Define
		\[
		\mathcal S(w,p)=(\widetilde w,\widetilde p).
		\]
		The linear estimates give
		\[
		\|\widetilde w\|_{1,\alpha;(\beta);D_L,\mathrm{pm}}^{(-\alpha)}
		+\|\widetilde p-p_0\|_{1,\alpha;(\beta-1,\beta);D_L,\mathrm{pm}}^{(-\alpha)}
		\le C_{\mathrm{map}}\hat{\kappa}\le K_0\hat{\kappa}.
		\]
		The boundary and contact conditions in \eqref{eq:Schauder-set} are
		preserved by the linear problem, so $\mathcal S(K)\subset K$.
		
		\par\medskip\noindent\textbf{Step 3. Continuity of the solution map.}
		Let $(w_n,p_n)\to(w,p)$ in $Y$, with all
		pairs in $K$. Write $A_n,A$ for the associated
		matrices, and set
		\[
		(\widetilde w_n,\widetilde p_n)=\mathcal S(w_n,p_n),
		\qquad
		(\widetilde w,\widetilde p)=\mathcal S(w,p).
		\]
		By \eqref{eq:Euler-coefficient-Lipschitz},
		$\|A_n-A\|_{L^\infty_{\mathrm{pm}}(D_L)}\to0$.
		The common boundary values of $\widetilde w_n$ and $\widetilde w$
		allow the zero-trace argument of Proposition~\ref{prop:slit-limit}
		to be applied to their difference, which therefore belongs to $H_0^1(D_L)$.
		Subtracting the weak equations gives
		\[
		\int_{D_L}A_n\nabla(\widetilde w_n-\widetilde w)
		\cdot\nabla\varphi\,dy
		=\int_{D_L}(A-A_n)\nabla\widetilde w
		\cdot\nabla\varphi\,dy,
		\qquad \varphi\in H_0^1(D_L).
		\]
		Testing with $\widetilde w_n-\widetilde w$ and using the vertical
		Poincar\'e inequality yields
		\begin{equation}\label{eq:map-energy-continuity}
			\|\widetilde w_n-\widetilde w\|_{H^1(D_L)}
			\le C\|A_n-A\|_{L^\infty_{\mathrm{pm}}(D_L)}
			\|\nabla\widetilde w\|_{L^2(D_L)}\longrightarrow0.
		\end{equation}
		The first-order equations then imply
		\[
		\begin{aligned}
			\|\nabla\widetilde p_n-\nabla\widetilde p\|_{L^2(D_L)}
			&\le C\|\nabla\widetilde w_n-\nabla\widetilde w\|_{L^2(D_L)}\\
			&\quad+\|A_n-A\|_{L^\infty_{\mathrm{pm}}(D_L)}
			\|\nabla\widetilde w\|_{L^2(D_L)}
			\longrightarrow0.
		\end{aligned}
		\]
		Since the outputs lie in the compact set $K$, every subsequence has
		a further subsequence converging in $Y$. The energy convergence in
		\eqref{eq:map-energy-continuity} identifies its first component as
		$\widetilde w$, while the preceding gradient convergence determines
		the second component up to one constant on each side. Continuity
		across $\Gamma_L^*$ and the common value $p_0$ at $e_L$ fix these
		constants. All subsequential limits are consequently
		$(\widetilde w,\widetilde p)$, which proves continuity of $\mathcal S$ in $Y$.
		
		\par\medskip\noindent\textbf{Step 4. The fixed point.}
		Schauder's fixed point theorem now gives $(w,p)=\mathcal S(w,p)$.
		At this pair, \eqref{eq:linear-first-order} is exactly
		\eqref{eq:nonlinear-rotated-system}, and the boundary and contact
		conditions follow from membership in $K$. It remains to check
		\eqref{eq:wp-admissibility}. The definition of $K$ and
		\eqref{eq:fixed-point-smallness} give $|w|\le\sigma_0$.
		By \eqref{eq:inlet-transported-profiles} and
		$\widehat\rho_{\mathrm{in}}\ge\rho_*$,
		$S\le(p_0+\sigma_0)/\rho_*^\gamma$.
		Together with the lower bounds for $p,q$ in
		\eqref{eq:Euler-state-bounds}, these inequalities give
		\[
		\rho u_1=\left(\frac pS\right)^{1/\gamma}
		\frac{q}{\sqrt{1+w^2}}\ge C^{-1}>0.
		\]
		Equation~\eqref{eq:Euler-state-bounds} also gives $p\ge c_b$ and
		$q^2/c^2\le1-\delta_b/2$. The pair therefore solves
		Problem~\ref{prob:finite-nonlinear}, and Step~2 gives
		\eqref{eq:nonlinear-final-bound}.
	\end{proof}
	
	\subsection{Uniqueness}
	\label{subsec:finite-uniqueness}
	
	The difference of two nonlinear solutions satisfies the first-order
	system in Lemma~\ref{lem:weak-difference}, with a source arising from
	the difference of their coefficient matrices. Estimate
	\eqref{eq:Euler-coefficient-Lipschitz} controls this source and allows
	the lemma to yield uniqueness when $\sigma_0$ is sufficiently small.
	
	\begin{proof}[Proof of Theorem~\ref{thm:finite-nonlinear}: uniqueness]
		Take $\hat{\kappa}_0,\sigma_0\le\sigma_b$, where $\sigma_b$ is given by
		Lemma~\ref{lem:Euler-coefficients}. Let $(w_1,p_1)$ and $(w_2,p_2)$
		be two solutions in the stated class, with coefficient matrices
		$A_1,A_2$. Define
		\begin{equation}\label{eq:nonlinear-differences}
			z=w_1-w_2,\qquad h=p_1-p_2,\qquad
			H=(A_1-A_2)\nabla w_2,
		\end{equation}
		and set
		$\mathfrak d=\|z\|_{0,0;\overline{D_L},\mathrm{pm}}
		+\|h\|_{0,0;\overline{D_L},\mathrm{pm}}$.
		The common boundary data and the zero-trace argument in
		Proposition~\ref{prop:slit-limit} give $z\in H_0^1(D_L)$.
		The assumed weighted bounds, \eqref{eq:uniform-weighted-integrability}
		with $s=2$, and $[h]=0$ give $h\in H^1(D_L)$.
		Subtracting the two first-order equations gives
		\begin{equation}\label{eq:exact-nonlinear-difference}
			\nabla_y^\perp h=A_1\nabla_y z+H,\qquad h(e_L)=0.
		\end{equation}
		Lemma~\ref{lem:Euler-coefficients} ensures uniform ellipticity of
		$A_1$ and gives
		\[
		\|A_1-A_2\|_{L^\infty_{\mathrm{pm}}(D_L)}\le C_b\mathfrak d.
		\]
		Combining this with the weighted bound for $w_2$ yields
		\[
		|H(y)|\le C_b\sigma_0\mathfrak d\,
		\Delta_y^{-\beta}d(y)^{\alpha-1}.
		\]
		By Lemma~\ref{lem:weak-difference},
		\[
		\mathfrak d\le C_U\sigma_0\mathfrak d,
		\]
		where $C_U$ is independent of $L$. Decreasing $\sigma_0$ so that
		$C_U\sigma_0<1$ forces $\mathfrak d=0$. After reducing $\hat{\kappa}_0$
		to preserve \eqref{eq:fixed-point-smallness}, the solution constructed
		above remains in this uniqueness class, since
		\[
		C_{\mathrm{map}}\hat{\kappa}\le K_0\hat{\kappa}\le\sigma_0.
		\]
	\end{proof}
	
	\section{Proof of Theorem~\ref{thm:infinite-existence}}
	\label{sec:infinite}
	
	The estimates in Theorem~\ref{thm:finite-nonlinear} are independent
	of $L$, so we can let the artificial outlet recede to infinity and
	recover a solution of Problem~\ref{prob:lagrangian}. Write $(w_L,p_L)$
	for the finite-nozzle solutions and $A_L$ for their coefficient matrices.
	
	\begin{proof}
		\noindent\textbf{Step 1. Passage to the limit.}
		Choose $\hat{\kappa}_0$ below the threshold of
		Theorem~\ref{thm:finite-nonlinear}. For fixed $R>\ell$, restrict
		the solutions with $L>R$ to $D_R$; the bound
		\eqref{eq:nonlinear-final-bound} remains uniform in $L$.
		The segment estimate in the proof of Lemma~\ref{lem:weighted-calculus}
		gives uniform $C^{0,\alpha}$ bounds for $w_L,p_L$ on
		$\overline{D_R^+}$ and $\overline{D_R^-}$. The weighted
		norms also give uniform $C^{1,\alpha}$ bounds on compact subsets away
		from $e_*$. A diagonal argument therefore gives a sequence $L_j\to\infty$ such that
		\[
		(w_{L_j},p_{L_j})\longrightarrow(w,p)
		\quad\text{in }C^0_{\mathrm{pm}}(\overline{D_R})\times
		C^0_{\mathrm{pm}}(\overline{D_R})
		\]
		for every fixed $R>\ell$, and in $C^{1,\alpha'}$ on compact subsets
		of each $\overline{D_\infty^\sigma}\setminus\mathcal E$, for
		$0<\alpha'<\alpha$. In particular, the inlet pressure traces converge
		uniformly. The formulas \eqref{eq:inlet-transported-profiles},
		\eqref{eq:coefficient-functions-SB} and \eqref{eq:chen-matrix}
		therefore imply
		\[
		\|A_{L_j}-A\|_{0,0;\overline{D_R},\mathrm{pm}}\longrightarrow0,
		\]
		where $A$ is the matrix determined by $(w,p)$ and its inlet pressure
		trace. The inequalities \eqref{eq:Euler-state-bounds} hold for
		$(p_{L_j},q_{L_j},c_{L_j})$ with constants independent of $L_j$.
		By \eqref{eq:coefficient-functions-SB}, these quantities converge
		locally uniformly, so the same inequalities hold for $(p,q,c)$.
		
		Each weighted pointwise and two-point inequality in
		\eqref{eq:nonlinear-final-bound} passes to the limit and gives
		\eqref{eq:lagrangian-main-estimate}.
		The first-order elliptic equations \eqref{eq:nonlinear-rotated-system} pass to the limit in $D_\infty^+$ and
		$D_\infty^-$. Uniform convergence up to the boundary preserves the
		inlet, wall and airfoil values, the common trailing-edge value $k$,
		and $[w]=[p]=0$ on $\Gamma_\infty^*$.
		Moreover, \eqref{eq:uniform-weighted-integrability} with $s=2$ gives
		\[
		\|w_{L_j}\|_{H^1(D_{L_j})}
		+\|\nabla p_{L_j}\|_{L^2(D_{L_j})}\le C\hat{\kappa}.
		\]
		Local weak convergence and exhaustion yield
		$w\in H^1(D_\infty)$ and $\nabla p\in L^2(D_\infty)$.
		The continuity of $w$ and $p$ across $\Gamma_\infty^*$ makes the
		limiting first-order relation valid on $D_\infty$ in the
		distributional sense. In particular,
		\[
		\int_{D_\infty}A\nabla_yw\cdot\nabla_y\phi\,dy=0,
		\qquad \phi\in H_0^1(D_\infty).
		\]

		The value weights in \eqref{eq:lagrangian-main-estimate} give
		\[
		\sup_{y_2\in I}|w(y_1,y_2)|\le C\hat{\kappa}(1+y_1)^{-\beta},
		\qquad
		\sup_{y_2\in I}|p(y_1,y_2)-p_0|
		\le C\hat{\kappa}(1+y_1)^{1-\beta}.
		\]
		Since $\beta>1$, the pressure tends uniformly to $p_0$ downstream.
		The bounds in \eqref{eq:Euler-state-bounds}, together with
		$|w|\le\sigma_0$ and $S\le(p_0+\sigma_0)/\rho_*^\gamma$,
		also give
		\begin{equation*}
			\rho u_1=\left(\frac pS\right)^{1/\gamma}\frac{q}{\sqrt{1+w^2}}\ge C^{-1}.
		\end{equation*}
		Thus $(w,p)$ solves Problem~\ref{prob:lagrangian}. We decrease
		$\hat{\kappa}_0$ if needed so that $C\hat{\kappa}_0\le\sigma_0$.
		
		\par\medskip\noindent\textbf{Step 2. Uniqueness.}
		Let $(w_1,p_1),(w_2,p_2)$ belong to the uniqueness class in the theorem,
		and let $A_1,A_2$ be their matrices. Define $z,h,H$ by
		\eqref{eq:nonlinear-differences} on $D_\infty$ and put
		$\mathfrak d=\|z\|_{0,0;\overline{D_\infty},\mathrm{pm}}
		+\|h\|_{0,0;\overline{D_\infty},\mathrm{pm}}$.
		The weighted bounds give $z\in H^1(D_\infty)$ and
		$\nabla h\in L^2(D_\infty)$. The zero boundary traces, a cutoff at
		infinity and the zero-trace argument at $e_*$ in
		Proposition~\ref{prop:slit-limit} give $z\in H_0^1(D_\infty)$.
		Moreover,
		\[
		\nabla_y^\perp h=A_1\nabla_yz+H,\qquad
		|H(y)|\le C_b\sigma_0\mathfrak d\,
		\Delta_y^{-\beta}d(y)^{\alpha-1},
		\]
		and $h(y_1,\cdot)\to0$ uniformly as $y_1\to\infty$.
		Lemma~\ref{lem:weak-difference} gives
		\[
		\mathfrak d\le C_U\sigma_0\mathfrak d.
		\]
		The choice $C_U\sigma_0<1$ gives $z=h=0$.
	\end{proof}
	
	\Needspace{9\baselineskip}
	\section{Recovery of the flow and proof of Theorem~\ref{thm:physical-main}}\label{sec:physical-return}
	
	The pair $(w,p)$ from Theorem~\ref{thm:infinite-existence} determines
	the density and velocity through the entropy and Bernoulli relations.
	We recover the inverse coordinate map by integrating the flow direction $w$ from
	the inlet map. Its boundary values identify the prescribed nozzle walls and
	airfoil surfaces and determine the contact curve.
	
	\begin{proof}
		The mass fluxes $m_\pm$ in \eqref{eq:mass-inlet-maps} and the inverse
		inlet maps $s_\pm$ in \eqref{eq:inverse-inlet-map} depend only on
		$\rho_{\mathrm{in}},u_{1,\mathrm{in}}$ and $H^\pm$.
		The positivity and $C^{1,\alpha}$ regularity of
		$\rho_{\mathrm{in}}u_{1,\mathrm{in}}$ give
		$s_\pm\in C^{2,\alpha}(\overline{I^\pm})$.
		Thus $\widehat\rho_{\mathrm{in}},\widehat u_{1,\mathrm{in}}$ have
		$C^{1,\alpha}_{\mathrm{pm}}(I)$ bounds determined by these data and
		satisfy \eqref{eq:physical-background} in the transformed coordinates.
		For the transverse inlet velocity, composition with $s_\pm$ gives
		\[
		\|\widehat u_{2,\mathrm{in}}\|_{1,\alpha;I,\mathrm{pm}}
		\le C_{\mathrm{EL}}\|u_{2,\mathrm{in}}\|_{1,\alpha;I_0,\mathrm{pm}}.
		\]
		Since $y_1=x_1$, the wall and airfoil terms in
		\eqref{eq:physical-kappa} and \eqref{eq:lagrangian-data-size} coincide.
		Consequently,
		\[
		\hat\kappa\le C_{\mathrm{EL}}\kappa,
		\]
		where $C_{\mathrm{EL}}\ge1$ depends only on
		$\rho_{\mathrm{in}},u_{1,\mathrm{in}},p_0,\ell,H^\pm,\gamma,\alpha,\beta$.
		We choose $C_{\mathrm{EL}}\kappa_0\le\hat{\kappa}_0$, where
		$\hat{\kappa}_0$ is the threshold in Theorem~\ref{thm:infinite-existence}.
		The inverse inlet maps satisfy
		\[
		s_+(m_+)=H^+,\qquad s_-(-m_-)=H^-,\qquad s_\pm(0)=0.
		\]
		Using these endpoint values in \eqref{eq:lagrangian-inlet-angle}, we obtain the
		first two lines of \eqref{eq:lagrangian-compatibility} from
		\eqref{eq:physical-compatibility}. The last line of \eqref{eq:lagrangian-compatibility} is the cusp condition
		\eqref{eq:airfoil-cusp}, unchanged because $y_1=x_1$.
		Theorem~\ref{thm:infinite-existence} therefore provides $(w,p)$
		satisfying \eqref{eq:lagrangian-main-estimate}, with
		$C\hat\kappa\le CC_{\mathrm{EL}}\kappa$ on the right-hand side.
		
		\par\medskip\noindent\textbf{Step 1. Recovery in Lagrangian coordinates.}
		With $S$ and $B$ given by \eqref{eq:inlet-transported-profiles},
		we recover the remaining flow variables by
		\begin{equation}\label{eq:chen-reconstruction}
			\rho=\left(\frac pS\right)^{1/\gamma},\qquad
			u_1=\frac{q}{\sqrt{1+w^2}},\qquad u_2=wu_1,
		\end{equation}
		where $c,q$ are defined by \eqref{eq:coefficient-functions-SB}.
		Consequently,
		\[
		c^2=\frac{\gamma p}{\rho},\qquad
		q^2=|u|^2=2B-\frac{2\gamma p}{(\gamma-1)\rho}.
		\]
		By \eqref{eq:lagrangian-main-estimate},
		\[|w(y)|+|p(y)-p_0|+|p(0,y_2)-p_0|\le C\kappa.\]
		Substituting these bounds and
		$\|\widehat u_{2,\mathrm{in}}\|_{0,0;I,\mathrm{pm}}\le C\kappa$
		into \eqref{eq:inlet-transported-profiles} and
		\eqref{eq:chen-reconstruction} gives
		\begin{equation}\label{eq:recovered-state-closeness}
			\|U-\widehat U_0\|_{0,0;D_\infty,\mathrm{pm}}\le C\kappa.
		\end{equation}
		In view of \eqref{eq:physical-background}, decreasing $\kappa_0$ ensures
		\begin{equation}\label{eq:recovered-state-bounds}
			p\ge\frac{p_0}{2},\qquad \rho\ge\frac{\rho_*}{2},\qquad
			u_1\ge\frac{u_*}{2},\qquad \frac{q^2}{c^2}\le1-\frac{\delta_b}{2}.
		\end{equation}
		Using $w=u_2/u_1$ and $c^2=\gamma p/\rho$ in
		\eqref{eq:chen-matrix} gives
		\begin{equation}\label{eq:Euler-matrix-explicit}
			a^{11}=\frac{u_1(c^2-u_1^2)}{c^2-q^2},\qquad
			a^{12}=a^{21}=-\frac{c^2\rho u_1u_2}{c^2-q^2},\qquad
			a^{22}=\frac{c^2\rho^{\,2}u_1q^2}{c^2-q^2}.
		\end{equation}
		Since $S$ and $B$ depend only on $y_2$, differentiation of the
		recovered variables gives
		\[
		\rho_{y_1}=\frac{p_{y_1}}{c^2},\qquad
		q_{y_1}=-\frac{p_{y_1}}{\rho q}.
		\]
		Hence
		\[
		\begin{aligned}
			\partial_{y_1}\!\left(\frac1{\rho u_1}\right)
			&=\frac{w}{\rho u_1(1+w^2)}w_{y_1}
			+\frac{c^2-q^2}{\rho^2c^2u_1q^2}p_{y_1},\\
			(u_2)_{y_1}
			&=\frac{u_1}{1+w^2}w_{y_1}
			-\frac{w}{\rho u_1(1+w^2)}p_{y_1}.
		\end{aligned}
		\]
		Substituting \eqref{eq:nonlinear-rotated-system} with the coefficients
		\eqref{eq:Euler-matrix-explicit} into these identities yields
		\begin{equation}\label{eq:kinematic-system}
			\partial_{y_1}\!\left(\frac1{\rho u_1}\right)=w_{y_2},\qquad
			(u_2)_{y_1}=-p_{y_2}.
		\end{equation}
		The identity for $q_{y_1}$ also gives
		\[
		\partial_{y_1}\!\left(\frac{|u|^2}{2}\right)
		=u_1(u_1)_{y_1}+u_2(u_2)_{y_1}=-\frac{p_{y_1}}\rho.
		\]
		Together with \eqref{eq:kinematic-system}, this yields
		\[
		\partial_{y_1}\!\left(u_1+\frac{p}{\rho u_1}\right)-\partial_{y_2}(pw)
		=\frac{u_1(u_1)_{y_1}+u_2(u_2)_{y_1}+\rho^{-1}p_{y_1}}{u_1}=0.
		\]
		This is the second equation of \eqref{eq:lagrangian-conservation}.
		The last equation is $B_{y_1}=0$. Hence the recovered variables solve
		\eqref{eq:lagrangian-conservation} in $D_\infty^+$ and $D_\infty^-$.
		At $y_1=0$, \eqref{eq:coefficient-functions-SB} and
		\eqref{eq:inlet-transported-profiles}, together with the inlet
		condition \eqref{eq:lagrangian-inlet-angle}, give
		\[
		\rho(0,y_2)=\widehat\rho_{\mathrm{in}}(y_2),\qquad
		u_i(0,y_2)=\widehat u_{i,\mathrm{in}}(y_2),\qquad i=1,2.
		\]
		Hence 
		\[(\rho u_1)^{-1}(0,y_2)=s'(y_2).\]
		
		Since $d((0,y_2))\ge\min\{1,\ell\}>0$,
		\eqref{eq:lagrangian-main-estimate} gives
		\[
		\|p(0,\cdot)-p_0\|_{1,\alpha;I,\mathrm{pm}}\le C\kappa.
		\]
		Applying this trace estimate to \eqref{eq:inlet-transported-profiles} and then
		using \eqref{eq:chen-reconstruction} and
		\eqref{eq:recovered-state-closeness}--\eqref{eq:recovered-state-bounds}, we
		obtain
		\begin{equation}\label{eq:recovered-weighted-bounds}
			\|\rho\|_{1,\alpha;D_\infty,\mathrm{pm}}^{(-\alpha)}
			+\|u_1\|_{1,\alpha;D_\infty,\mathrm{pm}}^{(-\alpha)}
			+\|u_1^{-1}\|_{1,\alpha;D_\infty,\mathrm{pm}}^{(-\alpha)}\le C.
		\end{equation}
		
		Since $p\to p_0$ and $w\to0$ uniformly downstream, the reconstruction
		formulas give the limiting profiles in these coordinates
		\begin{equation}\label{eq:downstream-states}
			\widehat\rho_\infty(y_2)=\left(\frac{p_0}{S(y_2)}\right)^{1/\gamma},\qquad
			\widehat u_{1,\infty}(y_2)=\sqrt{2B(y_2)-\frac{2\gamma p_0}{(\gamma-1)\widehat\rho_\infty(y_2)}}.
		\end{equation}
		The inlet pressure estimate and \eqref{eq:inlet-transported-profiles}
		show that these are positive $C^{1,\alpha}_{\mathrm{pm}}(I)$ profiles,
		with norms bounded in terms of the data listed in
		Theorem~\ref{thm:physical-main}. To estimate convergence to these profiles,
		we use the identities
		\[
		\rho-\widehat\rho_\infty
		=\widehat\rho_\infty\left[\left(\frac p{p_0}\right)^{1/\gamma}-1\right]
		\]
		and
		\[
		q^2-\widehat u_{1,\infty}^2
		=-\frac{2\gamma p_0}{(\gamma-1)\widehat\rho_\infty}
		\left[\left(\frac p{p_0}\right)^{(\gamma-1)/\gamma}-1\right].
		\]
		Together with $u_1=q/\sqrt{1+w^2}$, these identities and the weighted
		product estimates yield
		\begin{equation}\label{eq:downstream-state-estimate}
			\|\rho-\widehat\rho_\infty\|_{1,\alpha;(\beta-1);D_\infty,\mathrm{pm}}^{(-\alpha)}
			+\|u_1-\widehat u_{1,\infty}\|_{1,\alpha;(\beta-1);D_\infty,\mathrm{pm}}^{(-\alpha)}
			\le C\kappa.
		\end{equation}
		The weight $\beta-1$ also applies to the transverse derivatives:
		differentiating the limiting profiles produces terms proportional to
		$p-p_0$, whose decay is $(1+y_1)^{1-\beta}$.
		
		\par\medskip\noindent\textbf{Step 2. The nozzle and the contact discontinuity.}
		We define
		\begin{equation}\label{eq:inverse-map}
			\varphi(y_1,y_2)=s(y_2)+\int_0^{y_1}w(\xi,y_2)\,d\xi,
			\qquad \Phi(y_1,y_2)=(y_1,\varphi(y_1,y_2)).
		\end{equation}
		For $y_2\ne0$, we differentiate \eqref{eq:inverse-map} and use
		\eqref{eq:kinematic-system} to obtain
		\[
		\varphi_{y_1}=w,\qquad
		\varphi_{y_2}=\frac1{\widehat\rho_{\mathrm{in}}(y_2)\widehat u_{1,\mathrm{in}}(y_2)}+
		\int_0^{y_1}\partial_{y_1}\!\left(\frac1{\rho u_1}\right)(\xi,y_2)\,d\xi
		=\frac1{\rho u_1}(y_1,y_2)>0.
		\]
		Equation~\eqref{eq:lagrangian-main-estimate} gives the bound
		\[
		|w_{y_2}(\xi,y_2)|\le C\kappa(1+\xi)^{-\beta}
		\bigl(1+|\xi-\ell|^{\alpha-1}\bigr),\qquad \xi\ne\ell.
		\]
		The right-hand side is integrable in $\xi$ because $\alpha>0$ and
		$\beta>1$. Dominated convergence therefore extends the identities
		for $\varphi_{y_1},\varphi_{y_2}$ to the one-sided limits at $y_2=0$.
		The contact graph and its downstream limit are
		\begin{equation}\label{eq:asymptotic-map}
			\eta(x_1)=\zeta_*+\int_\ell^{x_1}w(\xi,0)\,d\xi,\qquad
			\eta_\infty=\zeta_*+\int_\ell^\infty w(\xi,0)\,d\xi.
		\end{equation}
		The endpoint identities $s(m_+)=H^+$ and $s(-m_-)=H^-$,
		together with the wall conditions for $w$, give
		$\varphi(y_1,m_+)=W^+(y_1)$ and
		$\varphi(y_1,-m_-)=W^-(y_1)$.
		Similarly, $s_\pm(0)=0$ and the airfoil boundary data give
		$\varphi^+(y_1,0)=\zeta^+(y_1)$ and
		$\varphi^-(y_1,0)=\zeta^-(y_1)$ for $0\le y_1\le\ell$.
		The equality $\zeta^+(\ell)=\zeta^-(\ell)=\zeta_*$ and continuity of
		$w$ across $\Gamma_\infty^*$ give the single contact graph
		\eqref{eq:asymptotic-map}.
		The bounds \eqref{eq:physical-kappa} and
		\eqref{eq:lagrangian-main-estimate}, together with $w(e_*)=k$, give,
		for $0<s<\min\{1,\ell\}$,
		\[
		\begin{aligned}
			|(\zeta^\sigma)'(\ell-s)-k|+|\eta'(\ell+s)-k|
			&\le C\kappa s^\alpha,\\
			|\zeta^\sigma(\ell-s)-\zeta_*+ks|
			+|\eta(\ell+s)-\zeta_*-ks|
			&\le C\kappa s^{1+\alpha}.
		\end{aligned}
		\]
		Thus each airfoil boundary joins the contact graph with
		$C^{1,\alpha}$ regularity. To complete the construction of the
		coordinate map, we use \eqref{eq:recovered-state-closeness}--\eqref{eq:recovered-weighted-bounds}
		and \eqref{eq:lagrangian-main-estimate} to obtain
		\[
		C^{-1}\le\varphi_{y_2}=\frac1{\rho u_1}\le C,\qquad
		\|D_y\Phi\|_{0,\alpha;\overline{D_\infty},\mathrm{pm}}\le C.
		\]
		These estimates and the boundary identities show that each restriction
		$\Phi:D_\infty^\sigma\to\Omega^\sigma$ is onto and uniformly
		bi-Lipschitz.
		
		The decay of $w$ also gives a uniform limit
		$\varphi_\infty(y_2)=\lim_{y_1\to\infty}\varphi(y_1,y_2)$.
		Passing to the limit in $\varphi_{y_2}=1/(\rho u_1)$ gives
		\[
		\varphi_\infty(y_2)=\eta_\infty+
		\int_0^{y_2}\frac{dt}{\widehat\rho_\infty(t)\widehat u_{1,\infty}(t)},\qquad
		\varphi_\infty'(y_2)=\frac1{\widehat\rho_\infty(y_2)\widehat u_{1,\infty}(y_2)}>0.
		\]
		Since $\beta>1$, the wall profiles have finite limits. The limiting
		nozzle cross-section is
		\[
		I_\infty:=\left(\lim_{x_1\to\infty}W^-(x_1),
		\lim_{x_1\to\infty}W^+(x_1)\right).
		\]
		The boundary identities give
		\[
		\varphi_\infty(\pm m_\pm)=\lim_{x_1\to\infty}W^\pm(x_1),\qquad
		\varphi_\infty(0)=\eta_\infty.
		\]
		Thus $\eta_\infty\in I_\infty$, and $\varphi_\infty$ maps $I^-$ and
		$I^+$ onto the portions of $I_\infty$ below and above $\eta_\infty$,
		respectively, as $C^{2,\alpha}$ diffeomorphisms up to the endpoints.
		The profiles in \eqref{eq:downstream-states} therefore define the
		physical asymptotic state by
		\[
		U_\infty(\varphi_\infty(y_2))
		=(\widehat\rho_\infty(y_2),\widehat u_{1,\infty}(y_2),0,p_0),\qquad y_2\in I^\pm,
		\]
		which belongs to $C^{1,\alpha}_{\mathrm{pm}}(I_\infty)$.
		Since $y_2$ is constant along each streamline, Step~1 gives convergence along
		streamlines to this physical state.
		
		\par\medskip\noindent\textbf{Step 3. The Euler equations and estimates.}
		We define $U(x)=U(\Phi^{-1}(x))$, retaining the notation for the
		flow variables. The inverse map satisfies
		$\psi_{x_2}=\rho u_1$ and $\psi_{x_1}=-\rho u_2$.
		Applying this change of variables to the Lagrangian equations
		verified in Step~1 yields \eqref{eq:physical-Euler} in each $\Omega^\sigma$.
		The boundary identities from Step~2 give
		\eqref{eq:physical-slip} and \eqref{eq:physical-inlet}.
		On the contact curve, they give \eqref{eq:physical-contact}, which
		ensures the Rankine--Hugoniot conditions. The bounded fluxes have integrals tending
		to zero on shrinking circles about $E_*$, so no point source remains.
		Since $x_1=y_1$, the uniform downstream pressure limit $p\to p_0$
		is also preserved.
		
		To transfer the estimates to the physical coordinates, we use the
		fact that the map preserves the first coordinate and satisfies
		\[
		c|y-e_*|\le|\Phi(y)-E_*|\le C|y-e_*|,
		\]
		with fixed positive constants $c,C$. For a scalar $v$,
		$D(v\circ\Phi^{-1})=(Dv\circ\Phi^{-1})D\Phi^{-1}$.
		Let $y,y'\in D_\infty^\sigma$ be distinct, and write $x=\Phi(y)$ and
		$x'=\Phi(y')$. The H\"older bound for $D\Phi^{-1}$ controls the
		additional term in the weighted gradient seminorm:
		\[
		\Delta_{y,y'}^\beta d(y,y')|Dv(y')|
		\frac{|D\Phi^{-1}(x)-D\Phi^{-1}(x')|}{|x-x'|^\alpha}
		\le C\|Dv\|_{0,\alpha;(\beta);D_\infty^\sigma}^{(1-\alpha)}
		d(y,y')^\alpha.
		\]
		Since $x_1=y_1$, the function-value weights are unchanged. Combining
		the preceding gradient estimate with the distance comparison gives
		\begin{equation}\label{eq:physical-coordinate-estimate}
			\|v\circ\Phi^{-1}\|_{1,\alpha;(\boldsymbol\beta);\Omega,\mathrm{pm}}^{(-\alpha)}
			\le C\|v\|_{1,\alpha;(\boldsymbol\beta);D_\infty,\mathrm{pm}}^{(-\alpha)}
		\end{equation}
		for $\boldsymbol\beta=(\beta,\beta)$ or $(\beta-1,\beta)$.
		The same argument without the factors $\Delta$ gives
		\eqref{eq:physical-coordinate-estimate} for $\boldsymbol\beta=(0,0)$.
		
		Since $u_2=u_1w$, the weighted product estimate, together with
		\eqref{eq:lagrangian-main-estimate} and \eqref{eq:recovered-weighted-bounds}, gives
		\[
		\|u_2\|_{1,\alpha;(\beta);D_\infty,\mathrm{pm}}^{(-\alpha)}
		\le C\|u_1\|_{1,\alpha;D_\infty,\mathrm{pm}}^{(-\alpha)}
		\|w\|_{1,\alpha;(\beta);D_\infty,\mathrm{pm}}^{(-\alpha)}
		\le C\kappa.
		\]
		Transferring this bound and the pressure estimate to the physical
		coordinates controls the first two terms in \eqref{eq:physical-main-estimate}.
		Since $\eta'=w(\cdot,0)$, taking the trace in
		\eqref{eq:lagrangian-main-estimate} gives
		\[
		\|\eta'\|_{1,\alpha;(\beta);(\ell,\infty)}
		^{(-\alpha;\{\ell\})}\le C\kappa.
		\]
		Moreover,
		\[
		\eta(x_1)-\eta_\infty=-\int_{x_1}^\infty w(t,0)\,dt,
		\]
		so $\beta>1$ yields
		\[
		\sup_{x_1>\ell}(1+x_1)^{\beta-1}|\eta(x_1)-\eta_\infty|\le C\kappa.
		\]
		Together with the bounds for $u_2$ and $p$, these estimates yield
		\eqref{eq:physical-main-estimate}.
		
		\par\medskip\noindent\textbf{Step 4. Uniqueness.}
		Fix $C$ in \eqref{eq:physical-main-estimate} large enough for the
		existence estimate above, and consider any solution satisfying that
		bound. In particular, $p-p_0$ and $u_2$ are bounded by $C\kappa$
		in the weighted norms displayed there. Along a streamline,
		the entropy and Bernoulli relations give
		\[
		\rho=\left(\frac pS\right)^{1/\gamma},\qquad
		u_1^2=2B-\frac{2\gamma p}{(\gamma-1)\rho}-u_2^2.
		\]
		Using \eqref{eq:physical-background},
		\eqref{eq:inlet-transported-profiles} and
		$\|p-p_0\|_{L^\infty(\Omega)}+\|u_2\|_{L^\infty(\Omega)}\le C\kappa$
		in these identities gives, for sufficiently small $\kappa_0$,
		\[
		C^{-1}\le\rho,u_1\le C.
		\]
		Hence the stream-function map and its inverse are uniformly
		Lipschitz on each side. The inlet trace of $p$ and
		\eqref{eq:inlet-transported-profiles} give
		$S,B\in C^{1,\alpha}_{\mathrm{pm}}(I)$ with bounded norms.
		Differentiating the two identities above and using
		\eqref{eq:physical-main-estimate} gives
		\[
		|\nabla_y\rho|+|\nabla_yu_1|+|\nabla_y(u_1^{-1})|
		\le C\bigl(1+\kappa\Delta_y^{-\beta}d(y)^{\alpha-1}\bigr).
		\]
		Integrating along line segments in $D_\infty^\sigma$ gives uniform
		$C^{0,\alpha}$ bounds for $\rho,u_1,u_1^{-1}$ and hence for $D\Phi$.
		The argument proving \eqref{eq:physical-coordinate-estimate} therefore
		also applies to composition with $\Phi$.
		The reconstruction formulas then give \eqref{eq:recovered-weighted-bounds}
		for these variables.
		Finally, the product estimate for $w=u_2/u_1$ gives
		\[
		\|w\|_{1,\alpha;(\beta);D_\infty,\mathrm{pm}}^{(-\alpha)}
		+\|p-p_0\|_{1,\alpha;(\beta-1,\beta);D_\infty,\mathrm{pm}}^{(-\alpha)}
		\le C_*\kappa.
		\]
		Here $C_*$ depends only on the data listed in
		Theorem~\ref{thm:physical-main} and the chosen constant
		in \eqref{eq:physical-main-estimate}. We take $C_*\kappa_0\le\sigma_0$,
		where $\sigma_0$ is the uniqueness threshold in
		Theorem~\ref{thm:infinite-existence}. Any two such solutions therefore have the same $(w,p)$. Their common inlet data and
		pressure trace determine the same $S,B$ and flow variables, and
		\eqref{eq:inverse-map} determines the same inverse map and contact graph.
	\end{proof}
	
	\section*{Acknowledgments}
	The research of Jun Chen was supported in part by the
	Hubei Provincial Natural Science Foundation of China (Grant No.~2025AFB530).
	The work of Xiaoguang You was partially supported by the National
	Natural Science Foundation of China (Grant No.~12561038) and the Jiangxi
	Provincial Natural Science Foundation (Grant No.~20242BAB25008).
	Aibin Zang is supported by the National Natural Science Foundation
	of China (Grant No.~12261093) and the Jiangxi Provincial Natural Science
	Foundation (Grant No.~20224ACB201004).
	
	\section*{Declaration of competing interests}
	The authors declare that they have no competing interests.
	
	\section*{Declaration of generative AI use}
	The authors used ChatGPT (OpenAI) to polish the language and help identify
	possible errors in the manuscript. The core mathematical ideas, methods,
	derivations, and calculations are the authors' own. The authors reviewed
	the AI-assisted revisions and take full responsibility for the content
	of the manuscript.

\end{document}